\documentclass[EJP]{ejpecp} 

\usepackage{epsfig}
\usepackage{bbm}

\AUTHORS{%
  Martin~Hutzenthaler\footnote{University of Munich, Germany.
    \EMAIL{hutzenthaler@bio.lmu.de}}
  }

\SHORTTITLE{Interacting diffusions and trees of excursions}

\TITLE{Interacting diffusions and trees of excursions:\\ convergence and comparison}

\KEYWORDS{Island model; virgin island model; mean field model; McKean-Vlasov limit; extinction; excursion measure; many-demes-limit; propagation of chaos}

\AMSSUBJ{60K35} 
\AMSSUBJSECONDARY{60E15; 92D25} 

\SUBMITTED{April 6, 2011} 
\ACCEPTED{August 17, 2012} 

\ARXIVID{1104.1060v1}  

\VOLUME{0}
\YEAR{2012}
\PAPERNUM{0}
\DOI{vVOL-PID}

\ABSTRACT{
  We consider systems of interacting diffusions
  with local population regulation
  representing populations on countably many islands.
  Our main result shows that the total mass process of such a system
  is bounded above by the total mass process of a tree of excursions with
  appropriate drift and diffusion coefficients.
  As a corollary, this entails a sufficient, explicit condition
  for extinction of the total mass as time tends to infinity.
  On the way to our comparison result, we establish that systems of interacting diffusions
  with uniform migration between finitely many islands converge to a tree of excursions
  as the number of islands tends to infinity.
  In the special case of logistic branching, this leads to a duality between
  a tree of excursions
  and the solution of a McKean-Vlasov equation.
}

\theoremstyle{ejpecpbodyit}
\newtheorem{assumption}{Assumption}[section]

\providecommand{\dup}{{\ensuremath{{\operatorname{d}}}}}
\providecommand{\Leb}{{\ensuremath{{\operatorname{Leb}}}}}

 \providecommand{\C}{{\ensuremath{\textup{C}}}}

 \providecommand{\E}{{\ensuremath{\mathbf{E}}}}

 \providecommand{\SetF}{{\ensuremath{\mathbbm{F}}}}

 \providecommand{\N}{{\ensuremath{\mathbbm{N}}}}

 \renewcommand{\P}{{\ensuremath{\mathbf{P}}}}
 
 \providecommand{\P}{{\ensuremath{\mathbf{P}}}}
 
 \providecommand{\R}{{\ensuremath{\mathbbm{R}}}}

 \providecommand{\1}{{\ensuremath{\mathbbm{1}}}}

 \providecommand{\MCA}{{\ensuremath{\mathcal A}}}

 \providecommand{\MCD}{{\ensuremath{\mathcal D}}}
 \providecommand{\MCE}{{\ensuremath{\mathcal E}}}
 \providecommand{\MCF}{{\ensuremath{\mathcal F}}}
 \providecommand{\MCG}{{\ensuremath{\mathcal G}}}

 \providecommand{\MCM}{{\ensuremath{\mathcal M}}}

 \providecommand{\MCV}{{\ensuremath{\mathcal V}}}

 \providecommand{\th} {{\ensuremath{\hat{t}}}}

 \providecommand{\Ft}    {{\ensuremath{\tilde{F}}}}

 \providecommand{\Mt}    {{\ensuremath{\tilde{M}}}}

 \providecommand{\Xt}    {{\ensuremath{\tilde{X}}}}
 \providecommand{\Yt}    {{\ensuremath{\tilde{Y}}}}
 \providecommand{\Zt}    {{\ensuremath{\tilde{Z}}}}

 \providecommand{\ct} {{\ensuremath{\tilde{c}}}}

 \providecommand{\ft} {{\ensuremath{\tilde{f}}}}

 \providecommand{\mt} {{\ensuremath{\tilde{m}}}}

 \providecommand{\Bb}    {{\ensuremath{\bar{B}}}}

 \providecommand{\Fb}    {{\ensuremath{\bar{F}}}}

 \providecommand{\Mb}    {{\ensuremath{\bar{M}}}}

 \providecommand{\Xb}    {{\ensuremath{\bar{X}}}}
 \providecommand{\Yb}    {{\ensuremath{\bar{Y}}}}

 \providecommand{\fb} {{\ensuremath{\bar{f}}}}

\providecommand{\mub}  {{\ensuremath{\bar{\mu}}}}
\providecommand{\mut}  {{\ensuremath{\tilde{\mu}}}}

\providecommand{\etab}{{\ensuremath{\bar{\eta}}}}

\providecommand{\dlt}{{\ensuremath{\tilde{\dl}}}}

\providecommand{\dlb}{{\ensuremath{\bar{\dl}}}}

\providecommand{\ldt}  {{\ensuremath{\tilde{\ld}}}}

\providecommand{\Pit}  {{\ensuremath{\tilde{\Pi}}}}

\renewcommand{\Xt}{{\ensuremath{\tilde{X}}}}

\providecommand{\zetah}{\ensuremath{\hat{\zeta}}}%
\providecommand{\zetab}{\ensuremath{\bar{\zeta}}}%
\providecommand{\zetat}{\ensuremath{\tilde{\zeta}}}%

\providecommand{\scal}[2]{{\ensuremath{\langle#1,#2\rangle}}}

\providecommand{\scalb}[2]{{\ensuremath{\bigl\langle#1,#2\bigr\rangle}}}

\providecommand{\ro}[1]    {{{(#1)}}}
\providecommand{\rob}[1]   {{{\bigl(#1\bigr)}}}
\providecommand{\robb}[1]  {{{\biggl(#1\biggr)}}}
\providecommand{\roB}[1]   {{{\Bigl(#1\Bigr)}}}

\providecommand{\ru}[1]    {{{(#1)}}}
\providecommand{\rub}[1]   {{{\bigl(#1\bigr)}}}
\providecommand{\rubb}[1]  {{{\biggl(#1\biggr)}}}
\providecommand{\ruB}[1]   {{{\Bigl(#1\Bigr)}}}

\providecommand{\sqb}[1]  {{{\bigl[#1\bigr]}}}

\providecommand{\sqB}[1]  {{{\Bigl[#1\Bigr]}}}

\providecommand{\eckb}[1]  {{{\bigl[#1\bigr]}}}
\providecommand{\eckbb}[1] {{{\biggl[#1\biggr]}}}
\providecommand{\eckB}[1]  {{{\Bigl[#1\Bigr]}}}

\providecommand{\curl}[1]   {{{\{#1\}}}}
\providecommand{\curlb}[1]  {{{\bigl\{#1\bigr\}}}}
\providecommand{\curlbb}[1] {{{\biggl\{#1\biggr\}}}}
\providecommand{\curlB}[1]  {{{\Bigl\{#1\Bigr\}}}}

\providecommand{\ddt}{{\ensuremath{\frac{d}{dt}}}}

\providecommand{\delzdxidxj}{{\ensuremath{\frac{\del^2}{\del x_{i}\del x_j}}}}

\providecommand{\abs}[1]  {{\ensuremath{|#1|}}}
\providecommand{\absb}[1] {{\ensuremath{\bigl|#1\bigr|}}}
\providecommand{\absbb}[1]{{\ensuremath{\Bigl|#1\Bigr|}}}
\providecommand{\absB}[1] {{\ensuremath{\Bigl|#1\Bigr|}}}

\providecommand{\al}      {{\ensuremath{\alpha}}}

\renewcommand{\th}    {{\ensuremath{\theta}}}
\providecommand{\ld}      {{\ensuremath{\lambda}}}

\providecommand{\eps}     {{\ensuremath{\varepsilon}}}
\providecommand{\dl}      {{\ensuremath{\delta}}}

\providecommand{\ldt}      {{\ensuremath{\tilde{\lambda}}}}

\providecommand{\dlt}      {{\ensuremath{\tilde{\delta}}}}

\providecommand{\taut}     {{\ensuremath{\tilde{\tau}}}}

\providecommand{\limdl}  {{\ensuremath{{\displaystyle \lim_{\delta \ra 0}}}}}
\providecommand{\limdlO} {{\ensuremath{{\displaystyle \lim_{\delta \ra 0}}}}}

\providecommand{\limm}   {{\ensuremath{{\displaystyle \lim_{m \ra \infty}}}}}

\providecommand{\limtO}  {{\ensuremath{{\displaystyle \lim_{t \ra 0}}}}}

\providecommand{\limK}   {{\ensuremath{{\displaystyle \lim_{K \ra \infty}}}}}
\providecommand{\limN}   {{\ensuremath{{\displaystyle \lim_{N \ra \infty}}}}}

\providecommand{\limsupK}  {{\ensuremath{{\displaystyle \limsup_{K \ra \infty}}}}}
\providecommand{\limsupN}  {{\ensuremath{{\displaystyle \limsup_{N \ra \infty}}}}}

\providecommand{\qqasN}  {\ensuremath{\qquad\text{as }N\to\infty}}

\providecommand{\ra}{\rightarrow}

\providecommand{\wlra} {\xrightarrow{w}}

\providecommand{\lra}{\longrightarrow}

\providecommand{\lradl}   {\xrightarrow{\dl  \ra0}}
\providecommand{\lradlO}  {\xrightarrow{\dl  \ra0}}

\providecommand{\lrat} {\xrightarrow{t  \ra\infty}}

\providecommand{\lraN} {\xrightarrow{N  \ra\infty}}

\providecommand{\wlimeps}%
      {{\ensuremath{\stackrel{\eps \rightarrow \infty}%
                            {\Longrightarrow}}}}

\providecommand{\varwlim}  {\xrightarrow{\text{w}}}

\providecommand{\sgn}{\operatorname{sgn}}

\providecommand{\Var}{\operatorname{Var}}

\providecommand{\mal}{\ensuremath{{\displaystyle \cdot}}}

\providecommand{\fa}{\ensuremath{\;\;\forall\;}}

\providecommand{\Gen}{{\ensuremath{\MCG}}}

\providecommand{\GenA}{{\ensuremath{\MCA}}}

\providecommand{\Dom}{{\ensuremath{\mathfrak{D}}}}

\providecommand{\uline}[1]{{\ensuremath{\underline{#1}}}}

\providecommand{\del}{{\ensuremath{\partial}}}

\providecommand{\clearemptydoublepage}%
    {\newpage{\pagestyle{empty}\cleardoublepage}}

\newlength{\mylen}

\newenvironment{beidschub}[1][5mm]{%
  \setlength{\mylen}{\textwidth}%
  \addtolength{\mylen}{-#1}%
  \addtolength{\mylen}{-#1}%
  \par\smallskip\noindent\hspace{#1}\begin{minipage}[t]{\mylen}}
  {\end{minipage}\par\smallskip\noindent}

\providecommand{\eqd}{{\ensuremath{\;\overset{d}{=}\;}}}

\newcommand{\cadlag}{c\`adl\`ag}

\renewcommand{\dup}{d}

\renewcommand{\E}{{\ensuremath{\mathbbm{E}}}}
\renewcommand{\P}{{\ensuremath{\mathbbm{P}}}}

\newcommand{\leqFix}{\ensuremath{\leq_{\SetF_{+-}}}}
\newcommand{\leqFiv}{\ensuremath{\leq_{\SetF_{++}}}}
\newcommand{\leqFst}{\ensuremath{\leq_{\SetF_{+\pm}}}}

\newcommand{\VIM}{\ensuremath{\MCV}}

\begin{document}



\section{Introduction}%
\label{sec:introduction}
In population dynamics and population genetics a prominent role is played by
diffusion processes on $I^G$ with $I=[0,\infty)$ or $I=[0,1]$ driven by stochastic differential
equations (SDEs) of the form 
\begin{equation}  \begin{split}\label{eq:X}
  dX_t(i)= \eckbb{\sum_{j\in G}X_t(j)m(j,i)-X_t(i)+ \mu\rub{X_t(i)}
                 }dt
         +\sqrt{\sigma^2\rub{X_t(i)}}\,dB_t(i),\quad i\in G,
\end  {split}     \end  {equation}
for $t\geq 0$
where $G$ is a finite or countable set and $(B_t(i))_{t\geq 0}$, $i\in G$,
are independent standard Brownian motions
and where $(m(j,i))_{j,i\in G}$ is a stochastic matrix.
We will refer to the solution $(X_t)_{t\geq0}$ of~\eqref{eq:X}
as $(G,m,\mu,\sigma^2)$-process.
In appropriate timescales and for suitable choices of $\mu$ and $\sigma^2$,
the component $X_t(i)$  describes the (rescaled) {\em population size}
on island $i\in G$ at time $t\geq 0$,  or  the relative frequency of a genetic type that is
present on island $i\in G$ at time $t\geq 0$.
The linear interaction term on the 
right-hand side of~\eqref{eq:X} models a {\em mass flow}
between the islands,
which might be caused by migration of individuals or a flow of genes.
Here we will use a picture from population dynamics.
The coefficient $\sigma^2(x)$ then is the infinitesimal variance of
the local population size given its current value  $x\in[0,\infty)$.
A classical case are Feller's branching diffusions where
$\sigma^2(x)= {\rm const}\cdot x$ for $x\in[0,\infty)$.
Moreover the drift term  $\mu$  describes the mean growth rate of a local population
apart from immigration from other islands and emigration.
A prototypical example  is $\mu(x) = x(K-x)$, i.e. logistic growth.

Here we address the question of the maximal effect of the migration matrix
$m(i,j)_{i,j\in G}$ for fixed $\mu$ and $\sigma^2$.
We will show in Theorem~\ref{thm:comparison} that the total mass process $(\sum_{i\in G}X_t(i))_{t\geq0}$
of $(X_t)_{t\geq0}$ is dominated by a tree of excursions,
which has been constructed in \cite{Hutzenthaler2009EJP}.
This dominating process does not depend on $G$
or on the migration matrix.
The intuition which leads to this comparison result is as follows.
Consider a model with supercritical population-size independent branching
and additional deaths due to competition within each island of $G$,
i.e., $\mu$ is a concave function with $\mu(0)=0$ and $\sigma^2$ is a linear function.
Now compare different distributions of individuals over space.
If there is at most one individual  on each island, then there are no deaths
due to competition until the next birth or migration event.
If, however, all individuals are on the same island, then there are deaths due to competition.
The effect of the population-size independent branching is the same in both situations.
We infer that more individuals survive if the distribution of individuals is more uniform
over space.
Now the migration dynamics which distributes mass uniformly over space would be
uniform migration on $G$.
As there is no uniform migration on an infinite set $G$, we approximate $G$ with larger and
larger finite subsets and consider uniform migration on the finite subsets.
This intuition leads to
considering an $N$-island model $(X_t^N)_{t\geq0}$ for $N\in\N$ which is
the solution of~\eqref{eq:X}
with $G:= \{1,\ldots, N\}$ and $m(i,j) := \frac 1N$, $i,j\in G$, that is,
the solution of
\begin{equation}  \begin{split}\label{eq:XN}
  dX_t^N(i)= \eckbb{\tfrac{1}{N}\sum_{j=1}^{N}X_t^N(j)-X_t^N(i)+ \mu\rub{X_t^N(i)}
                 }dt
         +\sqrt{\sigma^2\rub{X_t^N(i)}}\,dB_t(i),\quad i\in \{1,2,\ldots,N\},
\end  {split}     \end  {equation}
for $t\in[0,\infty)$ for every $N\in\N$.
We will refer to this $N$-island model as
$(N,\mu,\sigma^2)$-process.
Now if the initial configuration $X_0^N$ converges in distribution to $X_0$ as $N\to\infty$,
then the above intuition leads to
the assertion that the total mass of
the $(G,m,\mu,\sigma^2)$-process
is dominated by the limit of the total mass process of an $N$-island model
$(X_t^N)_{t\geq0}$ as $N\to\infty$
\begin{equation}   \label{eq:comp_intro}
  \sum_{i\in G}{X_t}(i)\leq\limN\sum_{i=1}^N{X_t^N}(i),\quad t\geq0,
\end{equation}
where the stochastic order being used here will be specified in~\eqref{eq:stochastic.order}
below.
We confirm this intuition in Theorem~\ref{thm:comparison}.
Please note that a comparison result
analogous to~\eqref{eq:comp_intro}
of a $(G,m,\mu,\sigma^2)$-process 
with an $N$-island process cannot be expected in general for finite $N\in\N$.

Next we review part of the literature on the comparison and on the limit
in~\eqref{eq:comp_intro} and begin
with the limits of the $N$-island model as $N\to\infty$.
Starting with Kac (1957)\nocite{Kac1957}
and McKean (1967)\nocite{McKean1967}, the case of exchangeable initial
configurations $X_0^N$, $N\in\N$, (so that $\sum_{j=1}^N X_0^N(j)$ is $O(N)$ as $N\to\infty$)
has been studied intensively
(e.g.~\cite{Gaertner1988, Leonard1986, Oelschlaeger1984} and the references therein).
The most general result assumes the drift and the diffusion coefficient to depend
continuously on $x_i$ and on the measure $\tfrac{1}{N}\sum_{j=1}^N\dl_{x_j}$ for
$x\in I^N$. If, in addition, certain Lyapunov conditions hold and if $\sigma(x)\neq 0$ for all $x\in I$, then
$(X_t^N(i))_{t\geq0}$ converges in distribution as $N\to\infty$ to a limiting process $(M_t)_{t\geq0}$
for every $i\in\N$ which in our case solves the McKean-Vlasov equation
\begin{equation}  \begin{split}\label{eq:M}
  dM_t=\eckB{\E M_t-M_t+\mu(M_t)}\,dt
         +\sqrt{\sigma^2\rub{M_t}}\,dB_t
\end  {split}     \end  {equation}
for $t\in[0,\infty)$, see Theorem 4.1 of~\cite{Gaertner1988}.
In particular, the Lyapunov assumptions of this theorem
are satisfied if the coefficients are locally Lipschitz continuous
and satisfy a linear growth condition,
see Proposition 5.1 of~\cite{Gaertner1988}.
To the best of our knowledge, our case of H\"older-$\tfrac{1}{2}$-continuous
diffusion coefficients which are not strictly positive has not been considered so far.
Proposition~\ref{p:MVL} below fills this gap.
The idea of comparing a 
$(G,m,\mu,\sigma^2)$-process
with a limit of $N$-island processes is due to Hutzenthaler and Wakolbinger (2007)\nocite{HutzenthalerWakolbinger2007}.
Their Proposition 2.2 establishes for the case of an associated initial configuration $X_0$
with identically distributed marginals and concave $\mu$ that $X_t(i)\leq M_t$ in a suitable stochastic order
for every $t\in[0,\infty)$ and every $i\in G$.
As a corollary hereof and of an extinction result for the McKean-Vlasov process~\eqref{eq:M},
Hutzenthaler and Wakolbinger (2007) obtain a sufficient condition for local extinction of the
$(G,m,\mu,\sigma^2)$-process, see Theorem 1 in~\cite{HutzenthalerWakolbinger2007}.

In this paper we focus on initial configurations with bounded total
mass (so that $\sum_{j=1}^N X_0^N(j)$ is $O(1)$ as $N\to\infty$).
The main difficulty here is that the limit of $(X_t^N)_{t\geq0}$ as $N\to\infty$
was unknown so far.
We will establish convergence of $(X_t^N)_{t\geq0}$ as $N\to\infty$
in Theorem \ref{thm:convergence}
below for the case of additive $\sigma^2$.
This is the first convergence result of the $N$-island process in case
of bounded total mass (together with the concurrent preprint~\cite{DawsonGreven2011pre} of Dawson and Greven
who consider interacting Wright-Fisher diffusions with selection and rare mutations).
The limiting process turns out to be a \emph{forest of (mass) excursions} started
in $x_i$, $i\in \N$.
This forest of excursions has been
constructed and analysed in \cite{Hutzenthaler2009EJP}
and has been denoted as \emph{virgin island model},
see also Section~\ref{sec:vim} for a formal definition.
In this model, every migrant populates a new island,
similar to the infinite-alleles model in which every mutation
produces a new allele.
Jean Bertoin turned this idea into
a {\em tree of alleles}
\cite{Bertoin2010SPA}
and a {\em partition into colonies}
\cite{Bertoin2010AAPpre}
(see also Section 7 of Pardoux and Wakolbinger (2010)\nocite{PardouxWakolbinger2010}
for a connection with the virgin island model).
In particular, Bertoin
\cite{Bertoin2010SPA}
shows in case of state-independent branching that the virgin island model
is the diffusion approximation of a discrete mass branching process.

Here is a brief heuristics how a forest of (mass) excursions
emerges from the $(N,\mu,\sigma^2)$-processes as $N\to \infty$.
Due to bounded initial mass, the total mass up to a finite time
is also bounded in $N\in\N$,
see Lemma~\ref{l:second_moment_estimate_X_Z}.
So the total mass that immigrates into a fixed island
is of order $O(1/N)$ as $N\to\infty$.
Among this immigrated mass might be a ``successful'' founder
producing a substantial family.
The probability of this event is of order $O(1/N)$
(this follows from~\eqref{eq:Q}
since $S(\eps)\sim S^{'}(0)\eps$ as $\eps\to\infty$).
Among $N$ islands, there is then a Poisson number of islands having
a founder whose progeny reaches a fixed level $\delta>0$, say.
Moreover all of these founders are on distinct islands. The reason for this
is that the probability to have two successful founders on the same
island is of order $O(1/N^2)$.
Consequently, this does not appear within $N$ islands
in the limit as $N\to\infty$,
see Lemma~\ref{l:vanishing_immigration_weak_process} for the details.
Thus, in the limit as $N\to\infty$, every ``successful'' emigrant populates a previously
unpopulated island.
The evolution of the population size on
every freshly populated island is described by a random (mass) excursion.
These mass excursions are born (densely in time) in a Poissonian manner
on ever new islands with an intensity proportional to the currently extant mass,
and, once born, follow the 
SDE 
\begin{align}\label{eq:Y}
  dY_t = -Y_tdt +\mu(Y_t)dt + \sqrt{\sigma^2(Y_t)}\, dB_t\, .
\end{align}
Formally, this is described by means of the {\em excursion measure}
$Q$ associated with \eqref{eq:Y} in the sense of
Pitman and Yor (1982)\nocite{PitmanYor1982} (see also \cite{Hutzenthaler2009EJP}).
The intensity measure with which a path $(\eta_t)_{t\geq0}$
spawns a ``daughter'' excursion born  at time $t\geq 0$ is $\eta_t \, dt \otimes \, Q$. 
The roots of the forest are random paths
which are independent solutions of \eqref{eq:Y}.
The virgin island model is then a countable family
$\VIM:=\{(s,\chi)\}$ of islands where island $(s,\chi)$
is populated at time $s\geq0$ and carries mass $\chi_{t-s}$
at time $t\geq0$.

\begin{figure}[ht]
\begin{beidschub}[0.0cm]
\begin{minipage}[t]{\linewidth}
  \epsfig{file=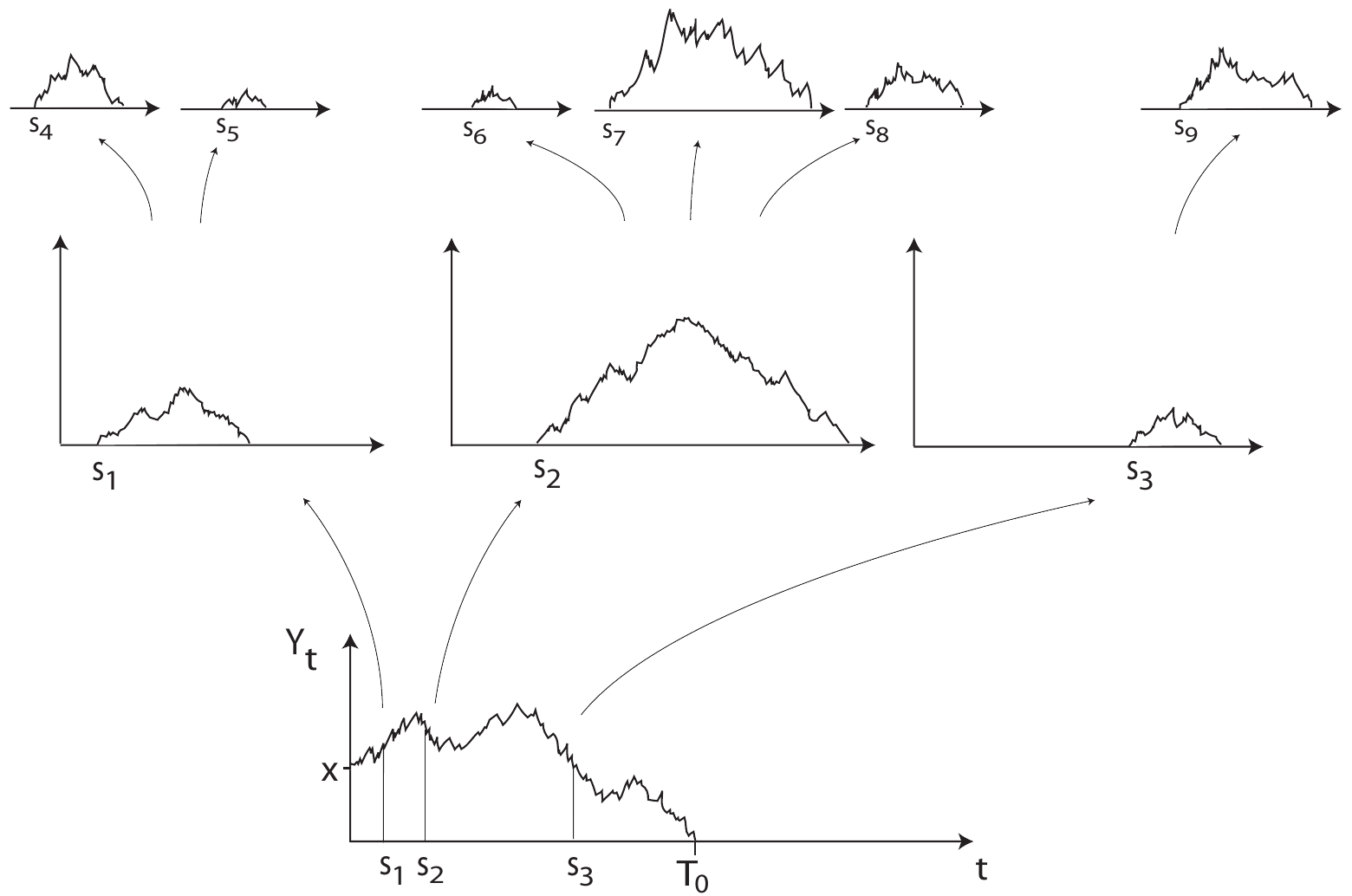, clip=,
      width=\linewidth}
  \caption{\footnotesize Subtree of the Virgin Island Model.
      Only offspring islands with a certain excursion height are drawn.
      Note that infinitely many islands are colonized e.g.\ between
      times $s_1$ and $s_2$.}
   \label{f:excursion_tree}
\end{minipage}
\end{beidschub}
\end{figure}

One tree of excursions in the forest of excursions is illustrated in
Figure~\ref{f:excursion_tree}.
Note that Figure~\ref{f:excursion_tree} does not contain the whole tree.
In fact islands are populated by emigrants densely in time
but only finitely many excursions started by these emigrants reach
a given strictly positive height.
Now a noteworthy observation is that the tree structure provides us with
independence of disjoint subtrees.
Putting this differently, the virgin island model is a branching process in discrete
time in the sense of Ji\v{r}ina (1958) except that there are now infinitely many
types, one for each excursion path.
Due to this branching structure, the virgin island model is easier to study
than the $N$-island process.
Several authors considered analogue of the virgin island 
model in the case of state-independent branching;
see \cite{AbrahamDelmas2009AIHPPS, Bertoin2009,
Bertoin2010AAPpre,
Bertoin2010SPA,
DawsonLi2003,
GriffithPakes1988,Lambert2009,SchinaziSchweinsberg2008, Taib1992}
for a selection of articles.

To state our convergence result of Theorem~\ref{thm:convergence} more formally,
denote the {\em population size spectrum} of the $(N,\mu,\sigma^2)$-process as
$\zeta_t^N := \sum_{i=1}^N\delta_{X^N_t(i)}$ for all $t\in[0,\infty)$ and
for every $N\in\N$.
Furthermore define the  population size spectrum of the virgin island model as
$\zeta_t := \sum_{(s,\chi)\in \VIM} \delta_{\chi_{t-s}}$ for all $t\in[0,\infty)$.
Our convergence result asserts that the population size spectra converge in distribution, i.e.,
$\zeta^N \to \zeta$
in distribution as $N\to \infty$; see Theorem~\ref{thm:convergence} for the precise statement.

Now we state the comparison result~\eqref{eq:comp_intro}
more precisely.
Recall that
$\mu$ is the infinitesimal mean in a non-spatial situation.
We assume $\mu$ to be {\em subadditive}.
Then
a population of size $x$ that is  separated into two islands experiences (in sum)
a larger growth rate than a population of the same size that is concentrated on one island.
Thus the virgin island model should offer in expectation a more prolific evolution of the total
mass than a model \eqref{eq:X} with the same coefficients $\mu$ and $\sigma^2$.
The infinitesimal variance $\sigma^2$ has an impact on
a comparison in distribution. More precisely the stochastic order in~\eqref{eq:comp_intro}
depends on whether $\sigma^2$ is superadditive, additive or subadditive.
In our prototype example of additive $\sigma^2$, the stochastic order is the usual increasing order.
In that case we have that
\begin{equation} \label{eq:comp_f_intro}
  \E\sqB{f\roB{
  \sum_{i\in G}{X_t}(i)}}
  \leq
  \E\sqB{f\roB{\sum_{(s,\chi)\in\VIM}\chi_{t-s}}},\quad t\geq0,
\end{equation}
for every non-decreasing function $f\colon[0,\infty)\to[0,\infty)$.
If $\sigma^2$ is superadditive or subadditive, then we will use concave, non-decreasing functions
and convex, non-decreasing functions, respectively.
In fact Theorem~\ref{thm:comparison} is a comparison result not only on the one-dimensional,
but on the finite dimensional distributions.
Thereby results on the distribution of the total mass of the virgin island model
have an immediate impact on the distribution of interacting diffusions with local population
regulation.

As a very special application of our comparison result, we obtain a sufficient condition
for global extinction 
for interacting diffusions with local population regulation.
Here we speak of global extinction
if the total mass $\sum_{i\in G}{X_t}(i)$ converges to zero in distribution
as $t\to\infty$ whenever $\sum_{i\in G}{X_0}(i)<\infty$.
Theorem 2 of~\cite{Hutzenthaler2009EJP} shows under certain assumptions that global extinction
of the virgin island model with coefficients $\mu$ and $\sigma$
is equivalent to
\begin{equation}   \label{eq:criterion_extinction}
  \int_0^{\infty}\frac{y}{\sigma^2(y)/2}\exp\robb{\int_0^y\frac{-x+\mu(x)}{\sigma^2(x)/2}\,dx}\,dy
  \leq1.
\end{equation}
As a consequence of~\eqref{eq:comp_f_intro}, if condition~\eqref{eq:criterion_extinction}
is satisfied, then
the $(G,m,\mu,\sigma^2)$-process dies out globally no matter what the
migration matrix $m(\cdot,\cdot)$ is.
Here are further implications of our comparison result.
Theorem 2 of~\cite{Hutzenthaler2009EJP} together with~\eqref{eq:comp_f_intro}
implies an upper bound
for the survival probability of $(X_t)_{t\geq0}$.
Theorem 3 of~\cite{Hutzenthaler2009EJP} yields an upper bound
for $\E[\int_0^\infty \sum_{i\in G}X_s(i)ds]$ if the left-hand side of~\eqref{eq:criterion_extinction}
is strictly smaller than one
and an upper bound for the growth rate of
$\int_0^t \E[ \sum_{i\in G}X_s(i)ds]$ as $t\to\infty$
if~\eqref{eq:criterion_extinction} fails to hold.
Moreover Theorem 4 of~\cite{Hutzenthaler2009EJP} implies an upper bound
for the growth rate of $\sum_{i\in G}X_t(i)$ as $t\to\infty$
if~\eqref{eq:criterion_extinction} fails to hold.

Here is a selection of models for which one could think of a similar
comparison result.
Mueller and Tribe (1994)~\nocite{MuellerTribe1994}
investigate a one-dimensional SPDE analog
of interacting Feller branching diffusions with logistic growth.
Bolker and Pacala (1997)~\nocite{BolkerPacala1997} propose a branching random walk
in which the individual mortality rate is increased by a weighted sum
of the entire population.
Etheridge (2004)~\nocite{Etheridge2004} studies two diffusion limits hereof.
The ``stepping stone version of the Bolker-Pacala model'' is a system
of interacting Feller branching diffusions with non-local logistic growth.
The ``superprocess version of the Bolker-Pacala model'' is an analog of this
in continuous space.
Blath, Etheridge and Meredith (2007)\nocite{BlathEtAl2007} study a two-type version hereof, which is
a spatial extension of the classical Lotka-Volterra model.
Fournier and M\'el\'eard (2004)\nocite{FournierMeleard2004} generalize the model of
Bolker and Pacala (1997)\nocite{BolkerPacala1997} by allowing spatial dependence of all rates.
A model in discrete time and discrete space is constructed
in Birkner and Depperschmidt (2007)\nocite{BirknerDepperschmidt2007}.
In that paper an individual has a Poisson number of offspring with mean depending on the
current configuration and, once created, offspring take an independent random
walk step from the location of their mother.

%
%
\section{The virgin island model}%
\label{sec:vim}
The virgin island model without immigration
has been introduced in~\cite{Hutzenthaler2009EJP}.
Here we slightly generalize this model by adding independent
immigration of mass.

The virgin island model is an analog of~\eqref{eq:X} in which every
emigrant populates a new island.
Islands with positive mass at time zero evolve as the one-dimensional
diffusion $(Y_t)_{t\geq0}$ solving~\eqref{eq:Y}.
The following assumption guarantees existence and uniqueness of
a strong $[0,\infty)$-valued solution of
equation~\eqref{eq:Y}, see e.g.\ Theorem IV.3.1
in~\cite{IkedaWatanabe1981}.
\begin{assumption}\label{a:A1}
 The set $I$ is an interval of length $\abs{I}\in(0,\infty]$ which is either
 of the form $[0,\abs{I}]$ if $\abs{I}<\infty$ or of the form
 $[0,\infty)$ if $\abs{I}=\infty$.
 The functions
 $\mu\colon I\to\R$ and
 $\sigma^2\colon I\to[0,\infty)$
 are  locally Lipschitz continuous in $I$ and satisfy
 $\mu(0) = 0=\sigma^2(0)$ and
 if $\abs{I}<\infty$, then $\mu(\abs{I})\leq 0=\sigma^2(\abs{I})$. 
 The function $\sigma^2(\cdot)$ is strictly positive on $(0,\abs{I})$
 and the function $\mu(\cdot)$ is globally upward Lipschitz continuous,
 that is, $\mu(x)-\mu(y)\leq L_\mu\abs{x-y}$ whenever $x>y\in I$ where
 $L_\mu\in[0,\infty)$ is a finite constant.
 Furthermore $\sigma^2$ satisfies
 the growth condition $\sigma^2(y)\leq L_{\sigma}(y+y^2)$
 for all $y\in I$ where $L_\sigma\in[0,\infty)$ is a finite constant.
\end  {assumption}
\noindent
Note that zero is a trap for $(Y_t)_{t\geq0}$, that is, $Y_t=0$
implies $Y_{t+s}=0$ for all $s\geq0$.

Mass emigrates from each island at rate one and colonizes new islands.
A new population should evolve as the process $(Y_t)_{t\geq0}$
and should start from a single individual which has mass zero
due to the diffusion approximation.
Thus we need the law of excursions of $(Y_t)_{t\geq0}$ from
the trap zero.
For this, define the
set of excursions from zero by 
\begin{equation}  \label{eq:U}
  U:=\bigl\{\chi\in\mathbf{C}\rub{(-\infty,\infty),[0,\infty)}\colon
     T_0\in(0,\infty],
     \,\chi_t=0\fa t\in(-\infty,0]\cup[T_0,\infty)\bigr\}
\end  {equation}
where $T_y=T_y(\chi):=\inf\{t>0\colon\chi_t=y\}$ is the first hitting time
of $y\in[0,\infty)$.
The set $U$ is furnished with locally uniform convergence.
For existence of the excursion measure $Q$
and in order to apply the results of~\cite{Hutzenthaler2009EJP},
we need to assume additional properties of $\mu(\cdot)$
and of $\sigma^2(\cdot)$.
For the motivation of these assumptions, we refer the reader
to~\cite{Hutzenthaler2009EJP}.
Assume $\int_0^\eps\tfrac{y}{\sigma^2(y)}dy<\infty$ for some $0<\eps<|I|$.
Then the scale function $S\colon I\to[0,\infty)$ defined through
\begin{equation} \label{eq:S}
  s(z):=\exp\ruB{-\int_0^z\frac{-x+\mu(x)}{\sigma^2(x)/2}\,dx},\quad
  S(y):=\int_0^y s(z)\,dz,\quad
  z,y\in I.
\end  {equation}
is well-defined.
\begin{assumption}      \label{a:Hutzenthaler2009EJP} 
  The functions $\mu(\cdot)$ and $\sigma^2(\cdot)$ satisfy
  \begin{equation}
    \int_{0}^{\eps} \frac{y}{\sigma^2(y)}\,dy<\infty
    \text{ and }
    \int_{\eps}^{\abs{I}} \frac{y}{\sigma^2(y) s(y)}\,dy<\infty
  \end{equation}
  for some $0<\eps<|I|$.
\end{assumption}
\noindent
Under Assumption~\ref{a:Hutzenthaler2009EJP}, the process $(Y_t)_{t\geq0}$
hits zero in finite time almost surely and the expected total
emigration intensity of the virgin island model is finite,
see Lemma 9.5 and Lemma 9.6 in~\cite{Hutzenthaler2009EJP}.
Moreover the scale function $S(\cdot)$ is well-defined and satisfies
$S^{'}(0)=1$.
A generic example which satisfies Assumption~\ref{a:Hutzenthaler2009EJP} is
$\mu(y)=c_1 y^{\kappa_1}-c_2 y^{\kappa_2}$, $\sigma^2(y)=c_3y^{\kappa_3}$
where $c_1,c_2,c_3>0$, $\kappa_2>\kappa_1\geq1$ and $\kappa_3\in[1,2)$.
Assumption~\ref{a:Hutzenthaler2009EJP} is not met by $\sigma^2(y)=y^2$.

The excursion measure $Q$ of the SDE~\eqref{eq:Y} has first been
described
in Pitman and Yor (1982)\nocite{PitmanYor1982}.
Here we use the description of $Q$ as limit
of the law of $(Y_t)_{t\geq0}$ started in $\eps>0$
and rescaled with $S(\eps)$ as $\eps\to0$.
More formally,
under Assuming~\ref{a:A1} and~\ref{a:Hutzenthaler2009EJP},
Theorem 1 in
\cite{Hutzenthaler2009EJP}
implies
that there exists a unique
$\sigma$-finite measure $Q$ on $U$ such that
\begin{equation}  \label{eq:Q}
  \lim_{\eps\to0}\frac{1}{S(\eps)}\E^\eps\sqb{ F\ru{Y}}=\int F(\chi)Q(d\chi)
\end  {equation}
for every bounded continuous function
$F\colon\mathbf{C}\rub{[0,\infty),[0,\infty)}\to\R$
for which there exists a $\dl>0$ such that $F(\chi)=0$ whenever
$\sup_{t\geq0}\chi_t<\dl$.
The reader might want to think of $Q$ as describing the evolution
of a population founded by a single individual.
In the special case $\sigma^2(y)=2\beta y$, $\mu(y)=cy$ with $\beta>0$ and
$c\in\R$, the process $(Y_t)_{t\geq0}$ is Feller's branching diffusion
whose law is infinitely divisible.
In that case the excursion measure coincides with the canonical measure.

Having introduced the excursion measure, we now construct
the virgin island model with constant immigration rate $\th\in[0,\infty)$
and started in $(x_k)_{k\in\N}\subset I$.
Let $\{(Y_{t}^{k,x_k})_{t\geq0}\colon k\in\N\}$
be independent solutions of~\eqref{eq:Y} such that $Y_0^{k,x_k}=x_k$
almost surely.
Moreover let $\Pi^{\th}$ be a Poisson point process
on $[0,\infty)\times U$
with intensity measure
\begin{equation}  \label{eq:intensity_measure_Pi_immi}
  \E\eckb{ \Pi^{\th}(dt\otimes d\psi)}
  =\th\,dt\otimes Q(d\psi).
\end  {equation}
The elements of the Poisson point process $\Pi^{\th}$
are the islands whose founders immigrated into the system.
Next we construct all islands which are colonized from a given
mother island.
Let
$\{\Pi^{(n,s,\chi)}\colon (n,s,\chi)\in\N_0\times[0,\infty)\times\C\big([0,\infty),I\big)\}$
be a set of independent
Poisson point processes on $[0,\infty)\times U$
with intensity measure
\begin{equation}  \label{eq:intensity_measure_Pi}
  \E\eckb{ \Pi^{(n,s,\chi)}(dt\otimes d\psi)}
  =\chi_{t-s}\,dt\otimes Q(d\psi)\quad (n,s,\chi)\in\N_0\times[0,\infty)\times\C\big([0,\infty),I\big).
\end  {equation}
All ingredients are assumed to be independent.
The elements of the Poisson point process $\Pi^{(n,s,\chi)}$
are the islands which descend from an island with population
size trajectory $(\chi_{t-s})_{t\geq0}$ and where the ancestral
lineages of individuals living on these islands have exactly
$n\in\N_0$ migration events.
Now the $0$-th generation is a random $\sigma$-finite measure on $[0,\infty)\times\C\big([0,\infty),I\big)$
defined through
  $\VIM^{(0)}:=\sum_{k\in\N}\dl_{(0,Y^{k,x_k})}+\Pi^{\th}$.
The $(n+1)$-st generation, $n\geq0$,
is the random $\sigma$-finite measure on all islands
which have been colonized from islands of the $n$-th generation, that is,
$\VIM^{(n+1)}:= \sum_{(s,\xi)\in\VIM^{(n)}}\Pi^{(n,s,\xi)}$ for all $n\in\N_0$.
The virgin island model $\VIM$ is then the sum of all of these measures
\begin{equation} \label{eq:nth_island}
  \VIM:=\sum_{n\in\N_{\geq0}}\VIM^{(n)}.
\end  {equation}
We call $\VIM$
the virgin island model with immigration rate $\th$
and initial configuration $(x_k)_{k\in\N}$.

%
%

\section{Main results}%
\label{sec:main results}
We begin with convergence of the
$N$-island process.
In this convergence, we allow the drift function $\mu_N$ and
the diffusion function $\sigma_N$
to depend on $N$ in order to include the case of weak immigration.
For example, one could be interested in an $N$-island model
with logistic branching and weak immigration at rate $\tfrac{\th}{N}$
on each island. In that case, one would set
$\mu_N(x)=\tfrac{\th}{N}+\gamma x(K-x)$ and $\sigma_N^2(x):=x$ for $x\in I$.
The equation of the $N$-island process now reads as
\begin{equation}  \begin{split}\label{eq:XN_mueN}
  dX_t^N(i)=\eckbb{\frac{1}{N}\sum_{j=1}^{N}X_t^N(j)-X_t^N(i)+ \mu_N\rub{X_t^N(i)}
                 }dt
         +\sqrt{ \sigma_N^2\rub{X_t^N(i)}}\,dB_t(i)
\end  {split}     \end  {equation}
where $i=1,\ldots,N$ and where 
$(B_t(i))_{t\geq0}$, $i\in \N$, are independent standard Brownian motions.
The idea to include weak immigration into a convergence result
is due to \cite{DawsonGreven2010pre}
who independently obtain convergence of an $N$-island model
using different methods.

Define $\mut_N(x):=\mu_N(x)-\mu_N(0)$ for $x\in I$.
For the $(N,\mu_N,\sigma_N^2)$-process to converge, we need
assumptions on $\mu_N,\sigma_N^2$ and on the initial distribution.
\begin{assumption}\label{a:A1_N_linear}
 Define $I:=[0,\infty)$.
 The functions 
 $\mu_N,\mu\colon I\to\R$
 are  locally Lipschitz continuous on $I$.
 The sequence $(\mu_N)_{N\in\N}$ converges pointwise to
 $\mu$ as $N\to\infty$.
 In addition, $N\mal\mu_N(0)\to\th\in[0,\infty)$ as $N\to\infty$
 and $0\leq N\mu_N(0)\leq 2\th$ for all $N\in\N$.
 The diffusion functions $\sigma_N^2$ and $\sigma^2$ are linear,
 that is, $\sigma_N^2(x)=\beta_N x$ and $\sigma^2 (x)=\beta x$ for
 some constants $\beta_N,\beta>0$ and all $N\in\N$. Furthermore $\beta_N$ converges to $\beta$
 as $N\to\infty$.
 Assumptions~\ref{a:A1} and~\ref{a:Hutzenthaler2009EJP} hold for $\mu$ and $\sigma^2$.
 Moreover $(\mu_N)_{N\in\N}$ is uniformly upward Lipschitz continuous,
 that is, $\mu_N(x)-\mu_N(y)\leq L_\mu \abs{x-y}$
 for all $x\geq y\in I$, $N\in\N$ and some constant $L_\mu\in[0,\infty)$.
\end  {assumption}
\noindent
Here is an example. If $\mu_N(x)=\tfrac{\th}{N}+C_1 x^{\kappa_1}-C_2 x^{\kappa_2}$
and $\sigma_N^2(x)=x$, then Assumption~\ref{a:A1_N_linear} is satisfied
if $1\leq\kappa_1<\kappa_2$ and $C_2(\kappa_1-1)>0$.
\begin{assumption}     \label{a:initial}
  The random variables $(X_0(i))_{i\in\N}$ and $(X_0^N(i))_{i\leq N}$
  are defined on the same probability space for each $N\in\N$.
  There exists a random permutation $\pi_1^N,\ldots,\pi_N^N$
  of $\{1,\ldots,N\}$ for each $N\in\N$ such that
  \begin{equation}
    \limN\E\eckB{\sum_{i=1}^N\abs{X_0(i)-X_0^N(\pi_i^N)}}=0.
  \end{equation}
  Furthermore the total mass of $X_0(\cdot)$ has finite
  expectation $\E\eckb{\sum_{i\in\N}X_0(i)}<\infty$.
\end{assumption}
\noindent
If $(x_i)_{i\in\N}\subset I$ is a summable sequence, then
Assumption~\ref{a:initial} is satisfied
for $X_0(i)=x_i$ and $X_0^N(i)=x_i$, $i\leq N\in\N$.

Next we introduce the topology for the weak convergence of the
$N$-island process.
What will be relevant here is not any specific numbering of the islands
but the statistics (or ``spectrum'') of their population sizes,
described by the sum of Dirac measures
at each time point, that is,
\begin{equation}   \label{eq:sum_dirac}
        \biggl(\sum_{i=1}^N \dl_{X_{t}^N(i)}\biggr)_{t\leq T}
\end{equation}
where $\dl_x$ is the Dirac measure on $x\in I$.
The state space of the measure-valued
process~\eqref{eq:sum_dirac} is the
set $\MCM_{\sigma}\rob{I}$ of $\sigma$-finite measures on $I$.
We equip the state space $\MCM_{\sigma}\rob{I}$ with the vague topology
on $I\setminus\{0\}$.
For weak convergence of $\MCM_{\sigma}\rob{I}$-valued processes,
we equip the space of \cadlag-functions from $[0,\infty)$
to $\MCM_\sigma\rob{I}$ with the Skorokhod topology (e.g.~\cite{EthierKurtz1986}).

Now we formulate the convergence of the
$(N,\mu_N,\sigma_N^2)$-process defined in~\eqref{eq:XN_mueN}.
\begin{theorem}    \label{thm:convergence}
  Suppose that $\rob{\mu_N}_{N\in\N}$ and $\rob{\sigma_N^2}_{N\in\N}$
  satisfy Assumption~\ref{a:A1_N_linear}
  and that the initial configurations $(X_0^N(i))_{i\leq N}$
  and $(X_0(i))_{i\in\N}$ satisfy Assumption~\ref{a:initial}.
  Then, for every $T\in[0,\infty)$, we have that
  \begin{equation}   \label{eq:convergence}
        \biggl(\sum_{i=1}^N \dl_{X_{t}^N(i)}\biggr)_{t\leq T}
    \varwlim
        \biggl(\sum_{(s,\eta)\in\VIM}\dl_{\eta_{t-s}}\biggr)_{t\leq T}
    \qqasN
  \end{equation}
  in distribution
  where $\VIM$ is the virgin island model with
  immigration rate
  $\th=\lim_{N\to\infty}N\mu_N(0)$ and initial configuration
  $(X_0(i))_{i\in\N}$.
\end{theorem}
\noindent
The proof is deferred to Section~\ref{sec:convergence_to_VIM}.

\begin{remark}  \label{r:convergence}
  For readability we rewrite the convergence in terms of test functions.
  The weak convergence~\eqref{eq:convergence} is equivalent to
  \begin{equation}  \label{eq:convergence_functional}
    \limN \E\eckbb{ F\biggl(\Bigl(\sum_{i=1}^N f\rub{X_t^N(i)}\Bigr)_{t\leq T}\biggr)}
    = \E\eckbb{ F\rubb{\Bigl(\sum_{(s,\eta)\in\VIM}
          f\rob{\eta_{t-s}}\Bigr)_{t\leq T}}}. 
  \end{equation}
  for every bounded continuous function
  $F$ on $\C\rob{[0,T],\R}$ and every continuous function
  $f\colon I\to\R$ with compact support
  in $(0,\abs{I})$.
  This equivalence follows from Theorem 2.2 of \cite{RoellyCoppoletta1986}
  if $\MCM_\sigma\rob{I}$ is equipped with the weak topology, the case of the vague topology
  follows analogously.
  Applications often require functions $f$ with non-compact support.
  The following condition might be useful in that case.
  Let $\Fb$ be a continuous function
  on $\C\rob{[0,T],\R}$ satisfying the Lipschitz condition
  \begin{equation}  \label{eq:Lip_F}
    \absb{\Fb\rob{\eta}-\Fb\rob{\etab}}
    \leq L_\Fb\sum_{j=1}^n \absb{\eta_{t_j}-\etab_{t_j}}
     \qquad\fa\eta, \etab\in\C\rob{[0,T],\R}
  \end{equation}
  for some $0\leq t_1\leq\cdots\leq t_n\leq T$.
  In addition let
  $\fb\colon I\to\R$ be a continuous function satisfying
  $\abs{\fb(x)}\leq L_\fb x$ for all $x\in I$.
  Then following the arguments in the proof of
  Lemma~\ref{l:vanishing_immigration_weak_process} below,
  one can show that~\eqref{eq:convergence_functional} holds
  with $F$ and $f$ replaced by $\Fb$ and $\fb$, respectively.
  \hfill\mbox{$\diamond$}
\end{remark}

The assumptions of Theorem~\ref{thm:convergence} are satisfied for branching
diffusions with local population regulation.
A prominent example is the $N$-island model
with logistic drift $\mu(x)=\gamma x(K-x)$ and with
$\sigma^2(x)=2\beta x$ for $x\in[0,\infty)$ and some constants $\gamma, K,\beta>0$.
More generally, Theorem~\ref{thm:convergence} can be applied
if $\mu(x)=\gamma x-c(x)$ and $\sigma^2(x)=2\beta x$ for $x\in[0,\infty)$
where $c(\cdot)$
is a concave function with $c^{'}(0)\in\R$.
We believe that Theorem~\ref{thm:convergence} also holds for
non-linear infinitesimal variances such as $\sigma^2(x)=x(1-x)$ in
case of the Wright-Fisher diffusion.
Our proof requires linearity only for one argument which is the step
from equation~\eqref{eq:phi_Delta} to equation~\eqref{eq:E_phi_Delta}.

In case of logistic branching, we obtain a noteworthy duality
of the total mass process
\begin{equation}  \label{eq:def:V}
  V_t:=\sum_{\ru{s,\chi}\in\VIM}\chi_{t-s},\quad t\geq0.
\end  {equation}
of the virgin island model with the mean
field model $(M_t)_{t\geq0}$ defined in~\eqref{eq:M}.
By Theorem 3 of~\cite{HutzenthalerWakolbinger2007}, systems of interacting Feller
branching diffusions with logistic drift satisfy
a duality relation which for the $(N,\gamma y(K-y),2\beta y)$-process
reads as
\begin{equation}     \label{eq:dual_XN}
  \E^{y\dl_1}\Big[e^{-\frac{\gamma}{\beta} x\sum_{i=1}^N{X_t^N(i)}}\Big]
  =\E^{x\uline{1}}\Big[e^{- \frac{\gamma}{\beta} X_t^N(1)y}\Big]
  \qquad \fa x,y,t\geq0
\end{equation}
where the notation $\E^{y\dl_1}$ refers to the initial
configuration $X_0^N=(y,0,\ldots,0)$
and $\E^{x\uline{1}}$ refers to $X_0^N=(x,\ldots,x)$.
This duality is established in~\cite{HutzenthalerWakolbinger2007} via a generator calculation,
in Swart (2006)\nocite{Swart2006} via dualities between Lloyd-Sudbury particle
models and in~\cite{AlkemperHutzenthaler2007} by following ancestral lineages of
forward and backward processes in a graphical representation.
Now let $N\to\infty$ in~\eqref{eq:dual_XN}.
Then the left-hand side converges to the Laplace transform of the
total mass process of the virgin island model (without immigration)
according to Theorem~\ref{thm:convergence}
and the right-hand side converges to the Laplace transform of the mean
field model~\eqref{eq:M} according to Proposition~\ref{p:MVL} below.
This proves the following corollary.
\begin{corollary}  \label{c:duality}
  Let $(V_t)_{t\geq0}$ be the total mass process of the virgin island model
  without immigration starting on only one island.
  Furthermore let $(M_t)_{t\geq0}$ be the solution of~\eqref{eq:M},
  both with coefficients $\mu(y)=\gamma y(K-y)$ and $\sigma^2(y)=2\beta y$
  for $y\in[0,\infty)$
  where $\gamma,K,\beta>0$.
  Then
  \begin{equation}
    \E^y \Big[e^{-\frac{\gamma}{\beta} x V_t}\Big]=\E^x \Big[e^{-\frac{\gamma}{\beta}  M_t  y}\Big]
    \qquad\fa x,y,t\geq0
  \end{equation}
  where $\E^y$ and $\E^x$ refer to $V_0=y$ and $M_0=x$, respectively.
\end{corollary}
Together with known results on the mean field model~\eqref{eq:M}, this corollary leads
to a computable expression for the extinction probability of the virgin island model.
\begin{corollary}
  Let $(V_t)_{t\geq0}$ be as in Corollary~\ref{c:duality}.
  Then $V_t$ converges  to a random
  variable $V_\infty$ in distribution as $t\to\infty$.
  If
  \begin{equation}  \label{eq:condition.logistic}
    \int_0^{\infty}
       \exp\robb{K\gamma x-\frac{\gamma\beta}{2}x^2}\mal\exp(-x)\,dx
    \leq1,
  \end{equation}
  then $\P[V_\infty=0]=1$. 
  If condition~\eqref{eq:condition.logistic} fails to hold, then
  \begin{equation}
    \P^y[V_\infty=0]=1-\P^y[V_\infty=\infty]
    =\int_0^\infty e^{-\frac{\gamma}{\beta}  yx}\,\Gamma_{\rho}(dx)\in(0,1)
  \end{equation}
  where $\P^y$ refers to $V_0=y\in(0,\infty)$.
  The parameter $\rho\in(0,\infty)$ is the unique solution of
  \begin{equation}
    \int_0^\infty y^{\frac{\rho}{\beta}}\rob{K-y}
               \exp\roB{\frac{\gamma K-1}{\beta}y-\frac{\gamma}{2\beta}y^2}\,dy=0
  \end{equation}
  and the probability distribution $\Gamma_\rho$ is defined by
  \begin{equation}  \label{eq:invariant.distribution}
    \Gamma_\rho(dx)=\frac{C_\rho}{\beta x}\exp\roB{\int_K^x\frac{(\rho-z)+\gamma z(K-z)}{\beta z}dz}\,dx
  \end{equation}
  on $(0,\infty)$ where $C_\rho$ is a normalizing constant.
\end{corollary}
\begin{proof}
  Theorem 2 of~\cite{Hutzenthaler2009EJP} shows convergence in distribution of $V_t$
  to $V_\infty$ as $t\to\infty$ and 
  $\P(V_\infty=0)=1$ if~\eqref{eq:condition.logistic} holds.
  If~\eqref{eq:condition.logistic} fails to hold, then
  Corollary~\ref{c:duality} together with convergence of $V_t$ implies convergence
  in distribution of $M_t$ to a variable $M_\infty$ as $t\to\infty$.
  The distribution of $M_\infty$ is necessarily an invariant distribution of
  the mean field model~\eqref{eq:M} and is nontrivial.
  Lemma 5.1 of~\cite{HutzenthalerWakolbinger2007} shows that there is exactly
  one nontrivial invariant distribution for~\eqref{eq:M} and this distribution
  is given by~\eqref{eq:invariant.distribution}.
\end{proof}

The second main result is a comparison of systems of locally regulated
diffusions with the virgin island model.
For its formulation, we introduce three stochastic orders which
are inspired by Cox et al.\ (1996)\nocite{CoxEtAl1996}.
Let $Z=(Z_t)_{t\geq0}$ and $\Zt=(\Zt_t)_{t\geq0}$ be two stochastic processes
with state space $I$. We say that $Z$ is dominated by $\Zt$ with respect
to a set $\SetF$ of test functions on path space if
\begin{equation}  \label{eq:stochastic.order}
  Z\leq_{\SetF} \Zt \;:\!\iff \E\sqb{ F(Z)}\leq \E\sqb{ F(\Zt)} \quad\fa F\in\SetF.
\end{equation}
The first order is 'the usual stochastic order' $\leq_{\operatorname{st}}$
in which $Z$ is dominated
by $\Zt$ if there is a coupling of $Z$ and $\Zt$ in which $Z_t$ is dominated
by $\Zt_t$ for all $t\geq0$ almost surely.
Assuming path continuity, an equivalent condition is as follows.
Denote the set of non-de\-crea\-sing test functions
of $n\in\N_{\geq0}$ arguments
by
\begin{equation}  \begin{split}  \label{eq:F+pm}
  \MCF_{+\pm}^{(n)}
  :=
  \MCF_{+\pm}^{(n)}(S)
  :=\Bigl\{f\colon S^n\to\R|
          \,\text{$f$ is non-decreasing, $f$ is bounded or $f\geq0$}\Bigr\}
\end{split}     \end{equation}
for a set $S\subseteq[0,\infty)$.
Furthermore let $\SetF_{+\pm}$ be the set of non-decreasing functions which
depend on finitely many time-space points
\begin{equation}  \begin{split}    \label{eq:SetF+pm}
  \SetF_{+\pm}
  &:=
    \SetF_{+\pm}(G,S)
  :=\Bigl\{F\colon\C\rub{[0,\infty)\times G,S}\to\R|
       \exists n\in\N_0\,\exists\,(t_1,i_1),\ldots,(t_n,i_n)\in[0,\infty)\times G\\
       &\qquad\qquad\qquad\qquad\qquad\qquad\ \;\exists f\in\MCF_{+\pm}^{(n)}(S) \text{ such that\ }
         F(\eta)=f\rob{\eta_{t_1}(i_1),\ldots,\eta_{t_n}(i_n)}\Bigr\}.
\end{split}     \end{equation}
If there is no space component, then we simply write
$\SetF_{+\pm}(S)$.
In this notation, $Z\leq_{\operatorname{st}}\Zt$ is equivalent to
$Z\leqFst\Zt$,
see Subsection 4.B.1 in Shaked and Shanthikumar (1994)\nocite{ShakedShanthikumar1994}.

We will use two more stochastic orders.
In the literature, the set of non-decreasing, convex functions is often used.
Here an adequate set is the collection of 
non-decreasing functions whose second order partial derivatives
are non-negative.
As we do not want to assume smoothness,
we slightly weaken the latter assumption.
We say for $1\leq i,j\leq n$ that
a function $f\colon\R^{n}\to\R$ is \emph{$(i,j)$-convex} if
\begin{equation}
  f(z+h_1 e_i+h_2 e_j)-f(z+h_1 e_i)-f(z+h_2 e_j)+f(z)\geq0
  \quad \fa z\in\R^n,h_1,h_2\geq0.
\end{equation}
Note that if $f$ is smooth, then this  is equivalent to $\delzdxidxj f\geq0$.
In addition note that $f$ is $(i,i)$-convex if and  only if $f$ is convex in 
the $i$-th component.
Moreover we say that $f$ is $(i,j)$-concave if $-f$ is $(i,j)$-convex.
A function is called directionally convex (e.g.~Shaked and Shanthikumar 1990\nocite{ShakedShanthikumar1990})
if it is $(i,j)$-convex for all $1\leq i,j\leq n$.
Such functions are also referred to as L-superadditive functions
(e.g.~R\"uschendorf 1983\nocite{Rueschendorf1983}).
Define the set of increasing, directionally convex functions as
\begin{equation}  \begin{split}  \label{eq:F++}
  \MCF_{++}^{(n)}
  :=\Bigl\{f\in\MCF_{+\pm}^{(n)}\colon
          \,f\text{ is $(i,j)$-convex for all $1\leq i,j\leq n$}\Bigr\}
\end{split}     \end{equation}
and
similarly $\MCF_{+-}$ with '$(i,j)$-convex' replaced by '$(i,j)$-concave'.
Furthermore define $\SetF_{++}$ and $\SetF_{+-}$ as in~\eqref{eq:SetF+pm}
with $\MCF_{+\pm}^{(n)}$ replaced by 
$\MCF_{++}^{(n)}$ and
$\MCF_{+-}^{(n)}$, respectively.
Now we have introduced three stochastic orders
$\leqFst$, $\leqFix$ and $\leqFiv$.
Note that $\MCF_{++}^{(n)}$ contains all mixed monomials and that
$\MCF_{+-}^{(n)}$ contains all functions $1-\exp\rob{-\sum_{i=1}^n\lambda_i x_i}$
with $\lambda_1,\ldots,\lambda_n\geq0$.

For the solution $(X_t)_{t\geq0}$ of~\eqref{eq:X} to be well-defined,
we additionally assume the migration matrix to be substochastic.
\begin{assumption}  \label{a:migration}
  The set $G$ is (at most) countable and the matrix $(m(j,i))_{j,i\in G}$ is non-negative
  and substochastic, i.e.,
    $m(j,k)\geq0$ and $\sum_{i\in G}m(j,i)\leq1$ for all $j,k\in G$.
\end{assumption}
\noindent
Note that
Assumption~\ref{a:A1} together with Assumption~\ref{a:migration} guarantees
existence and uniqueness of a strong solution of~\eqref{eq:X} with values
in $\{x\in I^G\colon \abs{x}<\infty\}$.
This follows from Proposition 2.1 and inequality (48) of~\cite{HutzenthalerWakolbinger2007}
by letting the weight function $\sigma_i\nearrow 1$
for $i\in G$ and using monotone convergence.
\begin{theorem}  \label{thm:comparison}
  Assume~\ref{a:A1}, \ref{a:Hutzenthaler2009EJP} and~\ref{a:migration}.
  If $\mu$ is concave and if $\sigma^2$ is superadditive, then
  we have that
  \begin{equation}  \label{eq:comparison}
    \biggl(\sum_{i\in G}{X_t(i)}\biggr)_{t\geq0}
    \leq_{\SetF_{+-}\ro{[0,\infty)}}
    \ruB{V_t}_{t\geq0}.
  \end{equation}
  If $\mu$ is concave and $\sigma^2$ is subadditive, then
  inequality~\eqref{eq:comparison} holds with
  $\SetF_{+-}$ replaced by $\SetF_{++}$.
  If $\mu$ is subadditive and $\sigma^2$ is additive, then
  inequality~\eqref{eq:comparison} holds with
  $\SetF_{+-}$ replaced by $\SetF_{+\pm}$.
\end{theorem}
\noindent
The proof is deferred to Section~\ref{sec:comparison_with_the_VIM}.
Comparisons of diffusions at fixed time points are well-known.
Cox et al.\ \nocite{CoxEtAl1996}(1996) establish a comparison
between the finite-dimensional distributions of two diffusions
where the test functions have product structure.
To the best of our knowledge,
Theorem~\ref{thm:comparison} is the first comparison result
with general test functions on
the finite-dimensional distributions.
The techniques we develop for this in  Subsection~\ref{ssec:Preservation of convexity}
might allow
to generalize the comparison results of
Cox et al.\ \nocite{CoxEtAl1996}(1996) on interacting diffusions
and the comparison results of
Bergenthum and R\"uschendorf
(2007)\nocite{BergenthumRueschendorf2007} on semimartingales.

The assumption of $\mu$ being subadditive is natural in the
following sense.
Let us assume that letting two $1$-island processes with initial
masses $x$ and $y$, respectively, evolve independently is better
in expectation
for the total mass than letting one
$1$-island process with initial mass $x+y$ evolve.
This assumption implies that
\begin{equation}
  \mu(x+y)
  =\limtO\frac{\E X_t^{x+y}-x-y}{t}
  \leq \limtO\frac{\E X_t^{x}+\E X_t^{y}-x-y}{t}
  =\mu(x)+\mu(y)
\end{equation}
for all $x,y,x+y\in I$ and thus
subadditivity of the infinitesimal mean $\mu$.
If $\sigma^2$ is not additive, then we need
the stronger assumption of $\mu$ being concave
for Lemma~\ref{l:preservation}.

From Theorem~\ref{thm:comparison} and a global extinction
result for the virgin island model,
we obtain a condition for global extinction of
systems of locally regulated diffusions.
According to Theorem 2 of~\cite{Hutzenthaler2009EJP}, the total mass
of the virgin island model converges in distribution to
zero as $t\to\infty$ if and only if condition~\eqref{eq:condition}
below is satisfied. Together with Theorem 2, this proves the
following corollary.
\begin{corollary}  \label{cor:global_extinction}
  Assume~\ref{a:A1}, \ref{a:migration} and~\ref{a:Hutzenthaler2009EJP}.
  Suppose that $\mu$ is subadditive and $\sigma^2$ is additive, or that
  $\mu$ is concave and
  $\sigma^2$ is superadditive.
  Then
  \begin{equation}  \label{eq:condition}
    \int_0^{\abs{I}}\frac{y}{\sigma^2(y)/2}
       \exp\robb{\int_0^y\frac{-x+\mu(x)}{\sigma^2(x)/2}\,dx}\,dy
    \leq1,
  \end{equation}
  implies global extinction of
  the solution $(X_t)_{t\geq0}$ of~\eqref{eq:X},
  that is, $\sum_{i\in G}{X_t(i)}\varwlim 0$ as $t\to\infty$ whenever
  $\sum_{i\in G}{X_0(i)}<\infty$ almost surely.
\end{corollary}
\begin{proof}
  The function $[0,\infty)\ni x\mapsto 1-e^{-\ld x}\in\mathcal{F}_{+\pm}\cap\mathcal{F}_{+-}$
  for every $\ld\in[0,\infty)$.
  If $\sigma^2$ is superadditive (or even additive),
  then Theorem 2 above and Theorem 2 of~\cite{Hutzenthaler2009EJP}
  imply that
  \begin{equation}
    \E\left[1-e^{-\ld\sum_{i\in G}{X_t(i)}}\right]
    \leq
    \E\left[1-e^{-\ld V_t}\right]
    \lrat 0
  \end{equation}
  for all $\ld \in[0,\infty)$. Convergence of the Laplace transform then implies weak convergence.
\end{proof}
In case of logistic branching ($\mu(y)=\gamma y(K-y)$, $\sigma^2(y)=2\beta y$),
condition~\eqref{eq:condition} simplifies to condition~\eqref{eq:condition.logistic}.

\section{Convergence to the virgin island model}
\label{sec:convergence_to_VIM}

\subsection{Outline}
\label{ssec:outline}

First we outline the intuition behind the proof.
The virgin island process is a tree of excursions whereas the $N$-island process
has no tree structure. It happens in the latter process that different 
emigrants colonize the same island. In addition, the $N$-island process is not
loop-free. An individual could migrate from island $1$ to island $3$ and then back
to island $1$. That these two effects vanish in the limit as the number
of islands tends to infinity will be established in two separate steps.

The first step
ensures that the limit of the $N$-island process as $N\to\infty$ is loop-free.
For this purpose, we decompose the $N$-island process according to the number
of migration steps. Throughout the paper, we say that an individual
has \emph{migration level $k\in\N_0$}
at time $t\in[0,\infty)$ if its ancestral lineage contains exactly $k$ migration steps.
For example, an individual starting on island $1$ at time $0$, moving to
island $3$ and then back to island $1$ has migration level $2$.
Let $N\in\N$.
We define a system
$\{(X_t^{N,k}(i))_{t\geq0}\colon k\in\N_0,i\leq N\}$
of diffusions such that
$X_t^{N,k}(i)$ consists of the mass on island $i$
at time $t$ with migration level $k$.
Recall $\mut_N(x)=\mu_N(x)-\mu_N(0)$ for all $x\in I$.
Define $X_0^{N,k}(i)=\1_{k=0}X_0^N(i)$ for all $i\in\{1,\ldots,N\}$
and $k\in\N_0$.
Let $\curlb{\rub{X_t^{N,k}(i),B_t^k(i)}_{t\geq0}\colon i\leq N, k\in\N_0}$
be a solution of
\begin{equation}  \begin{split}\label{eq:XN_k}
  dX_t^{N,k}(i)
    =&\left(\frac{1}{N}\sum_{j=1}^N X_t^{N,k-1}(j)-X_t^{N,k}(i)\right)dt
       +\1_{k=0}\mu_N(0)\,dt
     \\ &
     +\frac{X_t^{N,k}(i)}{\sum_{m\geq0}X_t^{N,m}(i)}
         \mut_N\rubb{\sum_{m\geq0}X_t^{N,m}(i)}dt\\
     & +\sqrt{\frac{X_t^{N,k}(i)}{\sum_{m\geq0}X_t^{N,m}(i)}
            \sigma_N^2\ruB{\sum_{m\geq0}X_t^{N,m}(i)}}\,dB_t^k(i)
\end  {split}     \end  {equation}
for all $t\in[0,\infty)$, $i=1,\ldots,N$ and $k\in\N_0$
where $X_t^{N,-1}:=0$ for $t\geq0$ and $i\leq N$
and where the family
$\curl{(B_t^k(i))_{t\geq0}\colon i,k\in\N_0}$
is a system of independent standard Brownian motions.
Here we implicitly used the continuous extension of
$\tfrac{x}{x+y}\mut_N(x+y)$ and of
$\tfrac{x}{x+y}\sigma^2_N(x+y)$ as functions of $(x,y)\in[0,\infty)^2\setminus\{(0,0)\}$
into the point $(0,0)$, where $N\in\N$.
Any weak solution of~\eqref{eq:XN_k} will be denoted
as a $(N,\mu_N,\sigma_N^2)$-process with migration levels.
See Lemma~\ref{l:existence_decomp}
for existence of a weak solution of~\eqref{eq:XN_k}.

Lemma~\ref{l:concentration} below indicates that the individuals with migration level $k$
at a fixed time are concentrated on essentially finitely many islands in the limit
$N\to\infty$.
A later migration event will not hit these
essentially
finitely many islands
because hitting a fixed island has
probability $\frac{1}{N}$.
Therefore we expect that all individuals on an island
have the same migration level.
Inserting this into~\eqref{eq:XN_k} suggests to consider
the solution
$\curlb{(Z_t^{N,k}(i))_{t\geq0}\colon i\leq N,k\in\N_0}$ of
\begin{equation}  \begin{split}   \label{eq:ZN}
  dZ_t^{N,k}(i)
    =&\rubb{\frac{1}{N}\sum_{j=1}^N Z_t^{N,k-1}(j)-Z_t^{N,k}(i)
         +\1_{k=0}\mu_N(0)+\mut_N\rub{Z_t^{N,k}(i)}}\,dt\\
     & +\sqrt{\sigma_N^2\rub{Z_t^{N,k}(i)}}\,dB_t^k(i),\quad
       Z_0^{N,k}(i)=X_0^{N,k}(i),\qquad i=1,...,N,
\end  {split}     \end  {equation}
where  $Z_t^{N,-1}:=0$ for all $t\geq0$ and $i\leq N$.
We will refer to this solution as the loop-free $(N,\mu_N,\sigma_N^2)$-process
or as loop-free $N$-island model.
Note that this is a $(\bar{G},\bar{m},\mu,\sigma^2)$-process with $\bar{G}:=\{1,2,\ldots,N\}\times\N_0$
and migration matrix $\bar{m}\big((i,k),(j,l)\big)=\tfrac{1}{N}\1_{l=k+1}$ for $(i,k),(j,l)\in\bar{G}$.
In particularly, we may and will choose
$\curlb{(Z_t^{N,k}(i))_{t\geq0}\colon i\leq N,k\in\N_0}$
to be the solution of~\eqref{eq:ZN} with respect to the Brownian motion of
a weak solution of~\eqref{eq:XN_k} for every $N\in\N$.
Consequently there exists a unique strong solution under Assumption~\ref{a:A1}.
Lemma~\ref{l:X_close_to_Z} below establishes
the assertion that the distance between the
$(N,\mu_N,\sigma_N^2)$-process with migration levels and
the loop-free $(N,\mu_N,\sigma_N^2)$-process
converges to zero in a suitable sense
as $N\to\infty$.
It turns out that some difficulties arise from the different
forms of the 
diffusion coefficients in the $(N,\mu_N,\sigma_N^2)$-process with migration
levels
and in the loop-free $(N,\mu_N,\sigma_N^2)$-process.
As we could not resolve these difficulties, we additionally assume for
Lemma~\ref{l:X_close_to_Z} that $\sigma_N^2$ is linear.
Then we have  that
$x\sigma_N^2(y)/y=\sigma_N^2(x)$ for $x=X_t^{N,k}(i)$ and
$y=\sum_{m=0}^\infty X_t^{N,m}(i)$
and the
diffusion coefficients in~\eqref{eq:XN_k} and in~\eqref{eq:ZN}
are similar.
Our proof of Lemma~\ref{l:X_close_to_Z} is a
moment estimate in the spirit of Yamada and Watanabe 
\nocite{YamadaWatanabe1971}(1971).

In Subsection~\ref{ssec:Poisson_limit} we show
that two emigrants colonize different islands in the limit $N\to\infty$.
Let us rephrase this more formally.
Recall that $\{(B_t^k(i))_{t\geq0}\colon i\leq N\}$ is independent
of $\{(B_t^{l}(i))_{t\geq0}\colon l<k,i\leq N\}$.
Thus, conditioned on $\{\sum_{j=1}^N Z_t^{N,k-1}(j)\colon t\geq0\}$,
the loop-free $(N,\mu_N,\sigma_N^2)$-process $(Z_t^{N,k}(i))_{t\geq0}$
on island $i$ with migration level $k\geq1$
evolves as the solution of
\begin{equation}  \label{eq:YN}
  dY_{t,s}^{N,\zeta}=\frac{\zeta_N(t)}{N}\,dt-Y_{t,s}^{N,\zeta}dt
                         +\mut_N\rub{Y_{t,s}^{N,\zeta}}\,dt
       +\sqrt{\sigma_N^2\rub{Y_{t,s}^{N,\zeta}}}dB_t,\ t\geq s,
\end{equation}
starting at time $s=0$ in $Y_{s,s}^{N,\zeta}=0$ driven
by the Brownian motion $(B_t)_{t\geq0}=(B_t^{k}(i))_{t\geq0}$
for each $i\leq N$
where
$\zeta_N(t):=\sum_{j=1}^N Z_t^{N,k-1}(j)$.
Note that $(Y_{t,s}^{N,\zeta}(i))_{t\geq s}$, $i\leq N$,
are independent and identically distributed.
Now let $\zeta_N\colon[0,\infty)\to N\cdot I:=\{N\cdot x\colon x\in I\}$
be a fixed path and
let $(\Yb_{t,s}^{N,\zeta}(i))_{t\geq s}$, $i\leq N$,
be independent solutions of~\eqref{eq:YN}.
We are interested in the total mass
$\big(\sum_{i=1}^N\Yb_{t,s}^{N,\zeta}(i)\big)_{t\geq s}$
as $N\to\infty$.
As $\Yb_{s,s}^{N,\zeta}(i)=0$ and as the immigration rate on island $i$
tends to zero, the
process $(\Yb_{t,s}^{N,\zeta}(i))_{t\geq s}$ converges to zero
as $N\to\infty$ for every $i\in\N$.
However, as mass of order $O(\frac{\zeta_N}{N})$ immigrates on a fixed island,
the probability that the excursion started by these immigrants reaches
a certain level $\dl>0$, say, is of order $O(\frac{\zeta_N}{N})$ as 
the convergence in~\eqref{eq:Q} indicates.
Now as there are $N$ independent trials, the Poisson limit theorem should
imply that
\begin{equation}      \label{eq:empirical_measure_converges}
   \biggl(\sum_{i=1}^{N}
    \dl_{\Yb_{t,s}^{N,\zeta}(i)}\biggr)_{s\leq t\leq T}\wlra
    \biggl(\int\dl_{\eta_{t-u}}\Pi(du,d\eta)\biggr)_{s\leq t\leq T}\qqasN
\end{equation}
where $\Pi$ is a Poisson point process with intensity measure
  $\limN\zeta_N(u) du\otimes Q(d\eta)$
if this limit exists.
We will prove~\eqref{eq:empirical_measure_converges}
in Lemma~\ref{l:vanishing_immigration_weak_process} by reversing time.

For convergence of the loop-free $(N,\mu_N,\sigma_N^2)$-process, we do
not need to assume linearity of the diffusion function. Here we may
replace Assumption~\ref{a:A1_N_linear} with the following weaker assumption.
\begin{assumption}\label{a:A1_N}
 The functions 
 $\mu_N,\mu\colon I\to\R$ and 
 $\sigma_N^2,\sigma^2\colon I\to[0,\infty)$
 are  locally Lipschitz continuous on $I$.
 The sequence $(\mu_N,\sigma_N^2)_{N\in\N}$ converges pointwise to
 $(\mu,\sigma^2)$ as $N\to\infty$.
 In addition, $N\mal\mu_N(0)\to\th\in[0,\infty)$ as $N\to\infty$
 and $N\mu_N(0)\leq 2\th$ for all $N\in\N$.
 The functions $\mu_N$ and $\sigma_N^2$ satisfy
 $\mu_N(0) \geq 0=\sigma_N^2(0)$ and
 if $\abs{I}<\infty$, then $\mu_N(\abs{I})\leq 0=\sigma_N^2(\abs{I})$. 
 Assumption~\ref{a:A1} and~\ref{a:Hutzenthaler2009EJP} hold for $\mu$ and $\sigma^2$.
 Moreover $(\mu_N)_{N\in\N}$ is uniformly upward Lipschitz continuous
 in zero, that is $\mu_N(x)-\mu_N(y)\leq L_\mu \abs{x-y}$
 for all $x\geq y\in I$, $N\in\N$ and some constant $L_\mu\in[0,\infty)$.
 The sequence $(\sigma_N^2)_{N\in\N}$ satisfies the
 uniform growth condition
 $\sigma_N^2(y)\leq L_{\sigma}(y+y^2)$
 for all $y\in I$, $N\in\N$ where $L_\sigma\in[0,\infty)$ is a finite constant
 and satisfies that $\liminf_{0<y\to 0}\inf_{N\in\N}\sigma_N^2(y)/\sigma^2(y)>0$.
\end  {assumption}
\noindent
Note that if $\sigma_N^2$ is linear, then
Assumption~\ref{a:A1_N} implies Assumption~\ref{a:A1_N_linear}.

some steps of our proof are based on second-moment estimates
and require the following assumption of uniformly finite second
moments of the initial distribution.
This assumption is then relaxed in further steps.

\begin{assumption}  \label{a:second_moments}
  The initial distribution satisfies that
  \begin{equation}
    \sup_{N\in\N}\E\eckB{\Bigl(\sum_{i=1}^N X_0^N(i)\Bigr)^2}<\infty.
  \end{equation}
\end{assumption}

\subsection{Preliminaries}
\label{ssec:preliminaries}
%
In this subsection we establish preliminary results such as moment estimates
and existence of the processes.
The quick reader might want to skip this subsection.
We begin with weak existence of the $N$-island process with migration levels.
%
\begin{lemma}  \label{l:existence_decomp}
  Assume~\ref{a:A1_N}.
  The $(N,\mu_N,\sigma_N^2)$-process with migration levels exists in the weak
  sense, that is,
  equation~\eqref{eq:XN_k} has a weak solution for every $N\in\N$.
\end{lemma}
\begin{proof}
  As the proof is fairly standard, we
  only give an outline.
  Approximate~\eqref{eq:XN_k} with stochastic differential
  equations for which weak solutions exist.
  For example, approximate $\mu_N$ and $\sigma_N^2$ locally uniformly with bounded continuous
  functions $\mu_{N,n}$ and $\sigma_{N,n}^2$, respectively.
  Consider the solution $(X_t^{N,k,n})_{t\geq0}$ of~\eqref{eq:XN_k} 
  with
  $\mu_N$ and $\sigma_N^2$ replaced by
  $\mu_{N,n}$ and $\sigma_{N,n}^2$, respectively,
  and which only depends on the migration levels $k\leq n$.
  Then this solution has a weak solution according to 
  Theorem V.23.5 and Theorem V.20.1 of~\cite{RogersWilliams2000b}
  as the coefficients are bounded and continuous and the stochastic differential
  equation is fi\-nite-di\-men\-sio\-nal.
  Show that the formal generator hereof converges to the formal
  generator associated with \eqref{eq:XN_k}. 
  In addition establish tightness of $(X^{N,k,n})_{n\in\N}$
  using moment estimates as in Lemma~\ref{l:second_moment_estimate_X_Z}
  for fixed $N\in\N$ but uniformly in $n\in\N$.
  Then apply the tightness criterion of Aldous (1978)\nocite{Aldous1978}.
  Then there exists a converging subsequence and its limit
  solves the martingale problem
  associated with
  \eqref{eq:XN_k}, see Lemma 4.5.1 in~\cite{EthierKurtz1986}.
  From this solution of the martingale problem, construct 
  a weak solution of
  \eqref{eq:XN_k} as in Theorem V.20.1 of~\cite{RogersWilliams2000b}.
\end{proof}
Next we prove that the $N$-island model with migration levels
is indeed a decomposition of the $N$-island model~\eqref{eq:XN_mueN}.
%
%
%
%
%
\begin{lemma}      \label{l:decomposition}
  Assume~\ref{a:A1_N}.
  Fix $N\in\N$ and
  let $\curlb{\rub{X_t^{N,k}(i),B_t^k(i)}_{t\geq0}\colon k\in\N_0,i\leq N}$
  be a solution of~\eqref{eq:XN_k}.
  Then $\curlb{(\Xt_t^{N}(i))_{t\geq0}\colon i\leq N}$
  defined through 
  \begin{equation}   \label{eq:all_levels}
    \Xt_t^N(i):=\sum_{k\geq 0}X_{t}^{N,k}(i),\quad t\geq0, i\leq N,
  \end{equation}
  is the unique solution of the $N$-island model~\eqref{eq:XN_mueN}
  corresponding to standard Brownian motions
  defined through
  \begin{equation}   \label{eq:corresponding_BM}
    d B_t(i)=\1_{\sum_{m\geq0}X_t^{N,m}(i)>0}
       \sum_{k\geq0}
       \sqrt{\frac{X_t^{N,k}(i)}{\sum_{m\geq0}X_t^{N,m}(i)}}dB_t^k(i)
       +\1_{\sum_{m\geq0}X_t^{N,m}(i)=0}dB_t^0(i)
    \quad t\geq0
  \end{equation}
  for every $i\leq N$.
\end{lemma}
\begin{proof}
  The process $(B_t(i))_{t\geq0}$ is a continuous martingale with
  quadratic variation process
  $[B(i),B(j)]_{t}=\dl_{ij}t$.
  Therefore L{\'e}vy's characterization (e.g.\ Theorem IV.33.1 in~\cite{RogersWilliams2000b})
  implies that
  \eqref{eq:corresponding_BM} defines a standard Brownian motion.
  Moreover it follows from summing~\eqref{eq:XN_k}
  over $k\in\N_0$ that $\rub{\Xt_t^N(i)}_{t\geq0}$
  solves~\eqref{eq:XN_mueN}.
  Pathwise uniqueness of~\eqref{eq:XN_mueN} has been established
  in Proposition 2.1 of~\cite{HutzenthalerWakolbinger2007}.
\end{proof}
  In the following lemmas,
  let the process $\curl{\ru{X_t^{N,k}(i),B_t^k(i)}_{t\geq0}\colon k\in\N_0,i\leq N }$
  be a solution of~\eqref{eq:XN_k} and let
  the process
  $\curl{\ru{Z_t^{N,k}(i),B_t^k(i)}_{t\geq0}\colon k\in\N_0,i\leq N }$
  be the solution of~\eqref{eq:ZN}.
  Define a stopping time $\taut_K^N\in[0,\infty)$ through
  \begin{equation}   \label{eq:tautKN}
    \taut_K^N:=\inf\curlB{t\geq0\colon
      \sum_{i=1}^N\sum_{m\geq0}\rub{X_t^{N,m}(i)+Z_t^{N,m}(i)}
      \geq K}
  \end{equation}
  for every $K\in[0,\infty)$ and every $N\in\N$.
%
\begin{lemma}  \label{l:first_moment_estimate_X_Z}
  Assume~\ref{a:A1_N}.
  Then we have that
  \begin{equation}  \label{eq:first_moment_estimate_X_Z}
    \sup_{t\leq T}\E\eckbb{\sum_{i=1}^N \sum_{k=0}^\infty
                 X_t^{N,k}(i)\Big|(X_0^{N,l})_{l\in\N_0}=(x^N\1_{l=0})_{l\in\N_0}}
    \leq \roB{2\th T+\sum_{i=1}^N x_i^N}e^{L_\mu T}
  \end{equation}
  for every configuration $x^N\in I^N$, every $T\in[0,\infty)$ and every $N\in\N$.
  The analogous assertion holds for the
  loop-free $(N,\mu_N,\sigma_N^2)$-processes, $N\in\N$.
\end{lemma}
\begin{proof}
  Fix $N\in\N$ and $T<\infty$.
  By Assumption~\ref{a:A1_N} we have that $\mu_N(x)\leq L_{\mu} x+\tfrac{2\th}{N}$
  for all $x\in I$.
  According to Lemma~\ref{l:decomposition},
  $\curlb{(\Xt_t^{N}(i))_{t\geq0}\colon i\leq N}$
  defined through~\eqref{eq:all_levels}
  is a solution of~\eqref{eq:XN_mueN}.
  Sum~\eqref{eq:XN_mueN} over $i\leq N$,
  stop at time $\taut_K^N$ and take expectations to obtain that
  \begin{equation}    \label{eq:first_moment_XN}
    \E\eckbb{\sum_{i=1}^N \Xt_{t\wedge\taut_K^N}^N(i)}
    \leq\sum_{i=1}^N x_i^N
     +\int_0^t L_\mu  \E\eckbb{\sum_{i=1}^N \Xt_{s\wedge\taut_K^N}^N(i)}+2\th\,ds
  \end{equation}
  for every $t\leq T$ and $K\in[0,\infty)$.
  Note that the right-hand side is finite.
  Now Gronwall's inequality implies that
  \begin{equation}
    \E\eckbb{\sum_{i=1}^N \Xt_{t\wedge\taut_K^N}^N(i)}
    \leq
    \roB{\sum_{i=1}^N x_i^N+2\th T} e^{L_\mu  T}
  \end{equation}
  for all $t\leq T$ and $K\in[0,\infty)$.
  Letting $K\to\infty$, path continuity and Fatou's lemma yield that
  \begin{equation}  \begin{split}
    \sup_{t\leq T}\E\eckbb{\sum_{i=1}^N \Xt_{t}^N(i)}
    &=
    \sup_{t\leq T}\E\eckbb{\sum_{i=1}^N \liminf_{K\to\infty}
                         \Xt_{t\wedge\taut_K^N}^N(i)}\\
    &\leq
    \sup_{t\leq T}\liminf_{K\to\infty}\E\eckbb{\sum_{i=1}^N 
                         \Xt_{t\wedge\taut_K^N}^N(i)}
    \leq
    \roB{\sum_{i=1}^N x_i^N+2\th T} e^{L_\mu  T}.
  \end{split}     \end{equation}
  This proves inequality~\eqref{eq:first_moment_estimate_X_Z}.
  The inequality for the
  loop-free $N$-island process follows similarly.
\end{proof}
%
\begin{lemma}   \label{l:essentially_finitely_many_generations}
  Assume~\ref{a:initial} and \ref{a:A1_N}.
  Then we have that
  \begin{equation}  \label{eq:essentially_finitely_many_generations}
    \sum_{k\geq0}\sup_{N\in\N}\sup_{t\leq T}\sum_{i=1}^N 
       \ruB{\E X_t^{N,k}(i)+\E Z_t^{N,k}(i)}
       <\infty
  \end{equation}
  for every $T\in[0,\infty)$.
\end{lemma}
\begin{proof}
  We prove inequality~\eqref{eq:essentially_finitely_many_generations}
  for the solution of~\eqref{eq:XN_k}.
  The estimate for the solution of~\eqref{eq:ZN} is analogous.

  Recall $\mut_N(x)=\mu_N(x)-\mu_N(0)$ for all $x\in I$.
  Apply It{\^o}'s formula to $\sum_{i=1}^N X_{t}^{N,k}(i)$,
  take expectations, estimate $\mut_N(x)\leq L_\mu x$ for all $x\in I$
  and take suprema to obtain that
  \begin{equation}  \begin{split}
    \sup_{N\in\N}\sup_{t\leq T}\E\sum_{i=1}^N& X_{t}^{N,k}(i)
    \leq \1_{k=0}\sup_{N\in\N}\E\sum_{i=1}^N X_0^N(i)\\
    &+\int_0^T 2\theta\1_{k=0}+\sup_{N\in\N}\sup_{s\leq t}\E\sum_{i=1}^N
              \rubb{ X_{s}^{N,k-1}(i)
              +L_\mu X_{s}^{N,k}(i)}\,dt
  \end{split}     \end{equation}
  for all $T\geq0$ and all $k\in\N_0$.
  Note that the right-hand side is finite due to
  Lemma~\ref{l:first_moment_estimate_X_Z} and Assumption~\ref{a:initial}.
  Summing over $k\leq K\in\N$ and applying Gronwall's inequality
  implies that
  \begin{equation}
    \sum_{k=0}^K\sup_{N\in\N}\sup_{t\leq T}\E\sum_{i=1}^N X_{t}^{N,k}(i)
    \leq \roB{\sup_{M\in\N}\E\sum_{i=1}^M X_0^M(i)+2\theta T}
    \cdot e^{(1+L_\mu)T}
  \end{equation}
  for every $K\in\N$. Letting $K\to\infty$
  proves~\eqref{eq:essentially_finitely_many_generations}.
\end{proof}
%
\begin{lemma}  \label{l:second_moment_estimate_X_Z}
  Assume~\ref{a:A1_N}.
  Then we have that
  \begin{equation}  \begin{split} \label{eq:second_moment_estimate_X_Z}
  \lefteqn{
    \E\eckbb{\biggl(\sup_{t\leq T}
        \sum_{i=1}^N\sum_{m\geq0}X_t^{N,m}(i)\biggr)^2
        \Big|(X_0^{N,l})_{l\in\N_0}=(x^N\1_{l=0})_{l\in\N_0}
        }
   }\\
    &\leq 4\biggl[\Bigl(\sum_{i=1}^N x_i^N \Bigr)^2+2\theta T+1\biggr]
               \Big(1+T(4\theta +L_\sigma)e^{L_\mu T}\Big)e^{(2L_\mu+L_\sigma)T}
  \end{split}     \end{equation}
  for every configuration $x^N\in I^N$, every $T\in[0,\infty)$ and every $N\in\N$.
  The analogous assertion holds for the
  loop-free $(N,\mu_N,\sigma_N^2)$-process.
\end{lemma}
\begin{proof}
  Fix $N\in\N$, $T\in[0,\infty)$ and a configuration $x^N\in I^N$.
  According to Lemma~\ref{l:decomposition},
  $(\sum_{m=0}^\infty X_t^{N,m})_{t\geq0}$ is an
  $N$-island model.
  Recall from Assumption~\ref{a:A1_N} that
  $\mu_N(x)\leq L_\mu x+\tfrac{2\th}{N}=:\mub_N(x)$ for all $x\in I$.
  Thus Lemma 3.3 of~\cite{HutzenthalerWakolbinger2007} implies that
  the $(N,\mu_N,\sigma_N^2)$-process is dominated by the
  $(N,\mub_N,\sigma_N^2)$-process $(\Xb_t^N)_{t\geq0}$ starting in $\Xb_0^N=x^N$.
  Using It{\^o}'s formula and $\sigma_N^2(x)\leq L_\sigma(x+x^2)$ for
  all $x\in I$, we get that
  \begin{equation}  \begin{split}    \label{eq:se_mo_XN}
    \lefteqn{\E\eckbb{\Bigl(\sum_{i=1}^N
                                    {\Xb_{t\wedge\taut_K^N}^{N}(i)}
                      \Bigr)^2
                     }
            }\\
   &= \Bigl(\sum_{i=1}^N x_i^N \Bigr)^2
    +\E\eckbb{\int_0^{t\wedge\taut_K^N} 2L_\mu \Bigl(\sum_{i=1}^N
                                    {\Xb_{s}^{N}(i)}
                      \Bigr)^2
     +4\th \sum_{i=1}^N \Xb_{s}^N(i)
     +\sum_{i=1}^N \sigma_N^2\rob{\Xb_{s}^N(i)}\,ds}
   \\
   &\leq \Bigl(\sum_{i=1}^N x_i^N \Bigr)^2
    +\int_0^t (2 L_\mu+L_\sigma) \E\eckbb{\Bigl(\sum_{i=1}^N
                                                 {\Xb_{s\wedge\taut_K^N}^{N}(i)}
                                      \Bigr)^2
                                     }
      +(4\th+L_\sigma)\Big(2\theta T+\sum_{i=1}^N x_i^N\Big)e^{L_\mu T}\,ds
   \\
   &\leq \Bigl(\sum_{i=1}^N x_i^N \Bigr)^2\big(1+T(4\th+L_\sigma)e^{L_\mu T}\big)
    +\!\!\int_0^t (2 L_\mu+L_\sigma) \E\eckB{\Bigl(\sum_{i=1}^N
                                                 {\Xb_{s\wedge\taut_K^N}^{N}(i)}
                                      \Bigr)^2
                                     }ds
      +\!T\big(4\th+L_\sigma\big)\big(2\theta T+1\big)e^{L_\mu T}
  \end{split}     \end{equation}
  for every $t\leq T$ and every $K\in\N$.
  We used Lemma~\ref{l:first_moment_estimate_X_Z}
  for the last but one inequality
  and the estimate $a\leq 1+a^2$ for $a\in\R$ for the last inequality.
  Note that the right-hand side is finite.
  Applying Doob's $L^2$ submartingale inequality
  (e.g.\ Theorem II.70.2 in~\cite{RogersWilliams2000a})
  to the submartingale $\sum_{i=1}^N \Xb_t^N(i)$,
  using Fatou's lemma
  and
  applying Gronwall's inequality to~\eqref{eq:se_mo_XN}, we conclude that
  \begin{equation}  \begin{split}
    &\E\eckbb{\sup_{t\leq T}
        \Bigl(\sum_{i=1}^N\sum_{m\geq0}X_t^{N,m}(i)\Bigr)^2}
    \leq 4\sup_{t\leq T}
    \E\eckB{ \Bigl(\sum_{i=1}^N\Xb_t^{N}(i)\Bigr)^2}
    \leq 4\sup_{t\leq T}
    \E\eckB{\liminf_{K\to\infty} \Bigl(\sum_{i=1}^N\Xb_{t\wedge\taut_K^N}^{N}(i)\Bigr)^2}
    \\
    &
    \leq 4\sup_{t\leq T}\liminf_{K\to\infty}
    \E\eckB{ \Bigl(\sum_{i=1}^N\Xb_{t\wedge\taut_K^N}^{N}(i)\Bigr)^2}
    \leq 4\biggl[\Bigl(\sum_{i=1}^N x_i^N \Bigr)^2+2\theta T+1\biggr]
               \Big(1+T(4\theta +L_\sigma)e^{L_\mu T}\Big)e^{(2L_\mu+L_\sigma)T}.
  \end{split}     \end{equation}
  The proof in the case of the loop-free $N$-island model is analogous.
\end{proof}
Recall $\taut_K^N$ from~\eqref{eq:tautKN}.
Next we show that stopping at the time $\taut_K^N$
has no impact within a finite time interval in the limit $K\to\infty$.
%
\begin{lemma}  \label{l:tau_theta}
  Assume~\ref{a:A1_N}
  and~\ref{a:second_moments}.
  Then any solution of~\eqref{eq:XN_k} satisfies that
  \begin{equation}
    \limsupK\sup_{N\in\N}
    \E\eckbb{\sup_{t\leq T}\sum_{i=1}^N\sum_{m\geq0}X_t^{N,m}(i)\1_{\taut_K^N\leq T}}
    =0
  \end{equation}
  for every $T\in[0,\infty)$.
  The analogous assertion holds for
  the loop-free $(N,\mu_N,\sigma_N^2)$-process.
\end{lemma}
\begin{proof}
  Rewriting
  $\{\taut_K^N\leq T\}
  =\{\sup_{t\leq T}\sum_{i=1}^N\sum_{m\geq0}(X_t^{N,m}(i)+Z_t^{N,m}(i))\geq K\}$,
  the assertion follows from the Markov inequality and from the second-moment estimate of
  Lemma~\ref{l:second_moment_estimate_X_Z}.
\end{proof}
%
\begin{lemma}  \label{l:average.square.Nisland}
  Assume~\ref{a:A1_N}.
  Then the $N$-island process $(X_t^N)_{t\geq0}$
  solving the SDE~\eqref{eq:XN_mueN}
  satisfies that
  \begin{equation}  \begin{split}
  \lefteqn{
    \E\left[\tfrac{1}{N}\sum_{i=1}^N\sup_{t\in[0,T]}\left(X_t^N(i)\right)^2\Big|X_0^N=x^N\right]
  }\\
    &\leq
    2\left(\tfrac{1}{N}\sum_{i=1}^N\big(x_i^N\big)^2
       +\tfrac{24T\theta^2}{N^2}
       +8L_\sigma T \Big(\tfrac{2\theta T}{N}+\tfrac{1}{N}\sum_{i=1}^N x_i^N\Big)e^{L_\mu T}
    \right)\exp\big( 40(1+T)(1+L_\mu+L_\sigma)^2T\big)
  \end{split}     \end{equation}
  for every configuration $x^N\in I^N$, every $T\in[0,\infty)$ and every $N\in\N$.
\end{lemma}
\begin{proof}
  Fix $N\in\N$ and $x^N\in I^N$ throughout the proof.
  Assumption~\ref{a:A1_N} implies that
  $\mu_N(x)\leq L_\mu x+\tfrac{2\th}{N}=:\mub_N(x)$ for all $x\in I$.
  Thus Lemma 3.3 of~\cite{HutzenthalerWakolbinger2007} implies that
  the $(N,\mu_N,\sigma_N^2)$-process $(X_t^N)_{t\geq0}$ is dominated by the
  $(N,\mub_N,\sigma_N^2)$-process $(\Xb_t^N)_{t\geq0}$ starting in $\Xb_0^N=x^N$.
  Applying Doob's $L^2$ submartingale inequality
  (e.g.\ Theorem II.70.2 in~\cite{RogersWilliams2000a}),
  Jensen's inequality and $\sigma_N^2(x)\leq L_\sigma(x+x^2)$ for all $x\in I$,
  we get that
  \begin{equation*}  \begin{split}
  \lefteqn{
    \E\left[\sup_{t\in[0,T]}\left(\bar{X}_t^N(i)-x_i^N\right)^2\right]
  }\\
    &=
    \E\left[\sup_{t\in[0,T]}\bigg(\int_0^t \tfrac{1}{N}\sum_{j=1}^N \bar{X}_s^N(j)-\bar{X}_s^N(i)
                                       +L_\mu \bar{X}_s^N(i)+\tfrac{2\theta}{N}ds
                                       +\int_0^t \sigma_N(\bar{X}_s^N(i))dB_s(i)\bigg)^2\right]
    \\
    &\leq
    2\E\bigg[\bigg(\int_0^T \tfrac{1}{N}\sum_{j=1}^N \bar{X}_s^N(j)
                     +(L_\mu+1) \bar{X}_s^N(i)+\tfrac{2\theta}{N}\,ds\bigg)^2\bigg]
   +2\E\left[\sup_{t\in[0,T]}\left(\int_0^t \sigma_N\big(\Xb_s^N(i)\big)dB_s(i)\right)^2\right]
    \\
    &\leq
    6T\int_0^T\E\bigg[\Big( \tfrac{1}{N}\sum_{j=1}^N \bar{X}_s^N(j)\Big)^2
                     +(L_\mu+1)^2 \Big(\bar{X}_s^N(i)\Big)^2+\tfrac{4\theta^2}{N^2}\bigg]ds
   +8\int_0^T \E\left[\sigma_N^2\big(\Xb_s^N(i)\big)\right]ds
    \\
    &\leq
    \int_0^T\E\bigg[6T\tfrac{1}{N}\sum_{j=1}^N \Big( \bar{X}_s^N(j)\Big)^2
       +\big(6T(L_\mu+1)^2+8L_\sigma\big)\Big( \bar{X}_s^N(i)\Big)^2
       +\tfrac{24T\theta^2}{N^2}
       +8L_\sigma\bar{X}_s^N(i)\bigg]ds
  \end{split}     \end{equation*}
  for all $i\in\{1,2,\ldots,N\}$
  and, using Lemma~\ref{l:first_moment_estimate_X_Z},
  \begin{equation*}  \begin{split}
  \lefteqn{
    \E\left[\tfrac{1}{N}\sum_{i=1}^N\sup_{t\in[0,T]}\left(\bar{X}_t^N(i)\right)^2\right]
    -2\tfrac{1}{N}\sum_{i=1}^N\big(x_i^N\big)^2
    \leq
    2\E\left[\tfrac{1}{N}\sum_{i=1}^N\sup_{t\in[0,T]}\left(\bar{X}_t^N(i)-x_i^N\right)^2\right]
  }\\
    &\leq
    40(1+T)(1+L_\mu+L_\sigma)^2
    \int_0^T\E\bigg[\tfrac{1}{N}\sum_{i=1}^N \sup_{r\in[0,s]} \Big( \bar{X}_r^N(i)\Big)^2\bigg]ds
       +\tfrac{48T\theta^2}{N^2}
       +16L_\sigma T \tfrac{1}{N}\Big(2\theta T+\sum_{i=1}^N x_i^N\Big)e^{L_\mu T}
  \end{split}     \end{equation*}
  for all $T\in[0,\infty)$.
  The right-hand side is finite due to Lemma~\ref{l:second_moment_estimate_X_Z}.
  Therefore, Gronwall's lemma implies that
  \begin{equation}  \begin{split}
  \lefteqn{
    \E\left[\tfrac{1}{N}\sum_{i=1}^N\sup_{t\in[0,T]}\left(X_t^N(i)\right)^2\right]
    \leq
    \E\left[\tfrac{1}{N}\sum_{i=1}^N\sup_{t\in[0,T]}\left(\bar{X}_t^N(i)\right)^2\right]
  }\\
    &\leq
    \left(2\tfrac{1}{N}\sum_{i=1}^N\big(x_i^N\big)^2
       +\tfrac{48T\theta^2}{N^2}
       +16L_\sigma T \tfrac{1}{N}\Big(2\theta T+\sum_{i=1}^N x_i^N\Big)e^{L_\mu T}
    \right)\exp\big( 40(1+T)(1+L_\mu+L_\sigma)^2T\big)
  \end{split}     \end{equation}
  for all $T\in[0,\infty)$ and this finishes the proof.
\end{proof}

Next we prove some preliminary results for the solution $(Y_{t,s}^{N,\zeta})_{t\geq0}$
of~\eqref{eq:YN}.
%
\begin{lemma}  \label{l:second_moment_estimate_Y}
  Assume~\ref{a:A1_N}.
  Let $\zeta_N\colon[0,\infty)\to N\mal I=\{N\cdot x\colon x\in I\}$
  be a locally square Lebesgue integrable function for every $N\in\N$.
  Then we have that
  \begin{equation}   \label{eq:second_moment_estimate}
      \E^x\eckB{\rub{\sup_{s\leq t\leq T}Y_{t,s}^{N,\zeta}}^2}
      \leq C_T\eckB{x+x^2+\int_s^T\! \frac{\zeta_N(r)}{N}
                  +\Bigl(\frac{\zeta_N(r)}{N}\Bigr)^2\,dr}
  \end{equation}
  for all $x\in I$, $0\leq s\leq T$, $N\in\N$ and some constant $C_T<\infty$
  which does not depend on $x$, $N$ or on $\zeta_N$.
\end  {lemma}
\begin{proof}
  The proof is similar to the proof of Lemma~\ref{l:second_moment_estimate_X_Z}, so
  we omit it.
\end{proof}
%
\begin{lemma}  \label{l:Y_2_supsum}
  Assume~\ref{a:A1_N}.
  Let $\zeta_N\colon[0,\infty)\to N\mal I$
  be a locally square Lebesgue integrable function for every $N\in\N$.
  Furthermore let $\rob{Y_{t,s}^{N,\zeta}(i)}_{t\geq s}$, $i\leq N$,
  be independent solutions of~\eqref{eq:YN} starting in 
  $Y_{s,s}^{N,\zeta}(i)=0$, $i\leq N$,
  for every $N\in\N$.
  Then we have that
  \begin{equation}   \label{eq:Y_2_supsum}
      \E\eckbb{\sup_{s\leq t\leq T}
          \Bigl(\sum_{i=1}^N Y_{t,s}^{N,\zeta}(i)\Bigr)^2
                }
      \leq C_T\eckB{\int_s^T\! \zeta_N(r)
                  +\frac{\bigl(\zeta_N(r)\bigr)^2}{N}\,dr}
  \end{equation}
  for all $0\leq s\leq T$, $N\in\N$ and some constant $C_T<\infty$
  which does not depend on $N$ or $\zeta_N$.
\end  {lemma}
\begin{proof}
  The proof is similar to the proof of Lemma~\ref{l:second_moment_estimate_X_Z}, so
  we omit it.
\end  {proof}
%
\begin{lemma}  \label{l:first_moment_estimate}
  Assume~\ref{a:A1_N} and fix $T\in[0,\infty)$.
  Let $(Y_{t,s}^{N,\zeta})_{t\geq s}$
  and $(Y_{t,s}^{N,\zetat})_{t\geq s}$ be two solutions of~\eqref{eq:YN} with respect
  to the same Brownian motion such that
  $Y_{s,s}^{N,\zeta}=x$ and $Y_{s,s}^{N,\zetat}=y$.
  If $\zeta_N,\zetat_N\colon[0,T]\to N\mal I$ are square Lebesgue integrable,
  then
  we have that
  \begin{equation}   \label{eq:moment_estimate}
    \E\big[(Y_{t,s}^{N,\zeta}-Y_{t,s}^{N,\zetat})^+\big]
    \leq e^{L_\mu (t-s)}\rubb{\frac{1}{N}\int_s^t(\zeta_N(r)-\zetat_N(r))^+\,dr+(x-y)^+}
  \end{equation}
  for all $N\in\N$, for all $0\leq s\leq t\leq T$ and  all $x,y\in I$
  where $z^+=\max(z,0)$ for all $z\in\R$.
\end{lemma}
\begin{proof}
  As in Theorem 1 of Yamada and Watanabe (1971)~\cite{YamadaWatanabe1971}, an approximation
  of $x\to x^+$ with $\C^2$-functions (see also the proof of Lemma~\ref{l:X_close_to_Z} for this approximation)
  results in
  \begin{equation}
    d\big(Y_{t,s}^{N,\zeta}-Y_{t,s}^{N,\zetat}\big)^+
    =\1_{Y_{t,s}^{N,\zeta}-Y_{t,s}^{N,\zetat}\geq 0}
       d\rub{Y_{t,s}^{N,\zeta}-Y_{t,s}^{N,\zetat}}.
  \end{equation}
  Taking expectations,
  the upward Lipschitz continuity
  of $\mut_N$ implies that
  \begin{equation}  \begin{split}
    \lefteqn{\E\big[\big(Y_{t,s}^{N,\zeta}-Y_{t,s}^{N,\zetat}\big)^+\big]}\\
    &\leq(x-y)^++\frac{1}{N}\int_s^t(\zeta_N(r)-\zetat_N(r)\big)^+\,dr
      +L_\mu \int_s^t \E\big[\big(Y_{r,s}^{N,\zeta}-Y_{r,s}^{N,\zetat}\big)^+\big]\,dr
  \end{split}     \end{equation}
  for all $t\geq s$.
  The right-hand side is finite due to Lemma~\ref{l:second_moment_estimate_Y}.
  Therefore, Gronwall's inequality implies~\eqref{eq:moment_estimate}.
\end{proof}
Now we study the solution $(Y_{t,s}^{N,c})_{t\geq s}$ of~\eqref{eq:YN}
in which $\zeta_N\equiv c$ is a constant $c\in[0,\infty)$ for every $N\in\N$
with $N\geq c/|I|$.
Let a point $\alpha\in(0,\abs{I})$ be fixed.
Recall the scale function $S$ from~\eqref{eq:S}.
Define the scale function $S_N$ of $(Y_{t,s}^{N,c})_{t\geq s}$ 
through
\begin{equation} \label{eq:SN}
  s_N(z):=\exp\ruB{-\int_{\alpha}^z\frac{c/N}{\sigma_N^2(x)/2}\,dx
                   -\int_0^z\frac{-x+\mut_N(x)}{\sigma_N^2(x)/2}\,dx},\ 
  S_N(y):=\int_0^y s_N(z)\,dz
\end  {equation}
for all $y,z\in I$.
We point out that using two reference points ($\al$ and $0$) in the
definition of $S_N$ is unusual but this definition differs from the
standard definition with a single reference point $\al$ only by
a constant factor.
The next two lemmas involve the speed measures
\begin{equation}  \label{eq:m}
  m(dy):=\frac{2}{\sigma^2(y)s(y)}\,dy,\quad
  m_N(dy):=\frac{2}{\sigma_N^2(y)s_N(y)}\,dy
\end{equation}
as measures on $(0,\abs{I})$ for every $N\in\N$.
Note that under Assumption~\ref{a:A1_N}, $m_N(\cdot)$
converges as a measure on $(0,\abs{I})$ vaguely to $m(\cdot)$ as $N\to\infty$.

%
\begin{lemma}   \label{l:replace_initial_point}
  Assume~\ref{a:A1_N} and~\ref{a:Hutzenthaler2009EJP}.
  If $\zeta_N\colon[0,\infty)\to N\mal I$ is locally square Lebesgue integrable for every $N\in\N$,
  then
  \begin{equation}
    \limdl\sup_{N\in\N}\int_0^\dl\absB{\E^x \big[f_N(Y_{t,s}^{N,\zeta})\big]
                                  -\E^0 \big[f_N(Y_{t,s}^{N,\zeta})\big]}\,m_N(dx)
    =0
  \end{equation}
  for all $s\leq t<\infty$ and 
  all functions $f_N\colon I\to\R$, $N\in\N$,
  with $\sup_{x\neq y\in I}\sup_{N\in\N}\tfrac{|f_N(x)-f_N(y)|}{|x-y|}<\infty$.
\end{lemma}
\begin{proof}
  Fix $s\leq t<\infty$ and
  define $C:=\sup_{x\neq y\in I}\sup_{N\in\N}\tfrac{|f_N(x)-f_N(y)|}{|x-y|}\in[0,\infty)$.
  Let $(Y_{u,s}^{N,\zeta,x})_{u\geq s}$ and
  $(Y_{u,s}^{N,\zeta,0})_{u\geq s}$ be solutions of~\eqref{eq:YN} with respect to the
  same Brownian motion satisfying $Y_{s,s}^{N,\zeta,x}=x$ and $Y_{s,s}^{N,\zeta,0}=0$.
  According to Assumptions~\ref{a:A1_N} and~\ref{a:Hutzenthaler2009EJP},
  there exist real numbers $\eps,\dl_0\in(0,|I|\wedge\alpha)$ such that
  $\sigma_N^2(y)\geq \eps\sigma^2(y)$ for all $y\in[0,\dl_0]$ and all $N\in\N$
  and such that $\int_0^{\dl_0}y/\sigma^2(y)\,dy<\infty$.
  The first-moment estimate of Lemma~\ref{l:first_moment_estimate} provides us with
  the inequality
  \begin{equation}  \begin{split}  \label{eq:replace_initial_point}
    \lefteqn{\sup_{N\in\N}\int_0^\dl\absB{\E^x\big[ f_N(Y_{t,s}^{N,\zeta})\big]
                                   -\E^0 \big[f_N(Y_{t,s}^{N,\zeta})\big]}\,m_N(dx)}\\
    &\leq
    \sup_{N\in\N}\int_0^\dl\E\eckB{C\absb{ Y_{t,s}^{N,\zeta,x} -Y_{t,s}^{N,\zeta,0}}}\,m_N(dx)
    \leq \sup_{N\in\N}\int_0^\dl C\cdot e^{L_\mu t}\cdot x\, m_N(dx)
    \\&
    \leq C e^{L_\mu t} \int_0^\dl x\frac{2}
                                        {\eps\sigma^2(x)}\,dx
           \cdot\exp\Big(\int_0^{\dl_0}\frac{2L_\mu z} {\eps\sigma^2(z)}\,dz\Big)
  \end{split}     \end{equation}
  for every $\dl\in(0,\dl_0)$.
  The right-hand side
  of~\eqref{eq:replace_initial_point}
  converges to zero as $\dl\to0$ by the dominated convergence
  theorem.
\end{proof}
%
\begin{lemma}  \label{l:order_of_m_N}
  Assume~\ref{a:A1_N}.
  If $c>0$, then
  \begin{equation}
    \frac{c}{N}m_N((0,\dl))\lra 1\qqasN
  \end{equation}
  for every $\dl\in(0,|I|\wedge\alpha)$ such that $\int_0^{\dl}y/\sigma^2(y)\,dy<\infty$ and
  $\inf_{y\in(0,\dl)}\inf_{N\in\N}\sigma_N^2(y)/\sigma^2(y)>0$.
\end{lemma}
\begin{proof}
  Fix
  $\dl\in(0,|I|\wedge\alpha)$ such that $\int_0^{\dl}y/\sigma^2(y)\,dy<\infty$ and
  $\inf_{y\in(0,\dl)}\inf_{N\in\N}\sigma_N^2(y)/\sigma^2(y)>0$
  and fix $c\in(0,\infty)$.
  Integration by parts yields that
  \begin{equation}  \begin{split}   \label{eq:conv_int_by_parts}
    \frac{c}{N} m_N\rob{(0,\dl)}
    &=\int_0^\dl\frac{c/N}{\sigma_N^2(y)/2}
         \exp\ruB{\int_{\alpha}^y\frac{c/N}{\sigma_N^2(x)/2}\,dx}
         \exp\ruB{\int_0^y\frac{-x+\mut_N(x)}{\sigma_N^2(x)/2}\,dx}\,dy\\
    &=\eckbb{\exp\ruB{\int_{\alpha}^y\frac{c/N}{\sigma_N^2(x)/2}\,dx
            +\int_0^y\frac{-x+\mut_N(x)}{\sigma_N^2(x)/2}\,dx}}_0^\dl\\
      &\quad-\int_0^\dl\frac{-y+\mut_N(y)}{\sigma_N^2(y)/2}
       \exp\ruB{\int_{\alpha}^y\frac{c/N}{\sigma_N^2(x)/2}\,dx
       +\int_0^y\frac{-x+\mut_N(x)}{\sigma_N^2(x)/2}\,dx}\,dy
  \end{split}     \end{equation}
  for every $N\in\N$.
  As $\sigma_N^2$ is Lipschitz continuous in $[0,\alpha]$ and $\sigma_N(0)=0$,
  $\int_{\alpha}^0\tfrac{1}{\sigma_N^2(x)/2}dx=-\infty$ for every $N\in\N$.
  Letting $N\to\infty$ in~\eqref{eq:conv_int_by_parts}
  and applying the dominated convergence
  theorem shows that
  \begin{equation}  \begin{split}
    \lefteqn{\limN\frac{c}{N}\, m_N\rob{(0,\dl)}}\\
    &=\exp\ruB{\int_0^\dl\frac{-x+\mu(x)}{\sigma^2(x)/2}\,dx}
      -\int_0^\dl\frac{-y+\mu(y)}{\sigma^2(y)/2}
       \exp\ruB{\int_0^y\frac{-x+\mu(x)}{\sigma^2(x)/2}\,dx}\,dy\\
    &=\exp\ruB{\int_0^\dl\frac{-x+\mu(x)}{\sigma^2(x)/2}\,dx}
      -\eckB{\exp\ruB{\int_0^y\frac{-x+\mu(x)}{\sigma^2(x)/2}\,dx}}_0^\dl
  \end{split}     \end{equation}
  which is equal to one.
\end{proof}
We recall the following lemma from~\cite{Hutzenthaler2009EJP}, see
Lemma 9.8 there.
\begin{lemma}   \label{l:finite_excursion_area}
  Assume~\ref{a:A1} and~\ref{a:Hutzenthaler2009EJP}.
  Let $Q$ be the excursion measure defined through~\eqref{eq:Q}.
  Then
  \begin{equation}
    \int \left(\int_0^\infty \chi_t\,dt\right) Q(d\chi)
    =\int_0^{\abs{I}}\frac{y}{\sigma^2(y)/2}
       \exp\robb{\int_0^y\frac{-x+\mu(x)}{\sigma^2(x)/2}\,dx}\,dy
    <\infty.
  \end{equation}
\end{lemma}
The last result of this subsection is a variation of the second moment
estimate of Lemma~\ref{l:second_moment_estimate_X_Z}.
Define a stopping time $\tau_K^N\in[0,\infty]$ through
\begin{equation}   \label{eq:tauKN}
  \tau_K^N:=\inf\Big(\curlB{t\geq0\colon
    \sum_{i=1}^N X_t^{N}(i) \geq K}\cup\{\infty\}\Big)
\end{equation}
for every $K\in[0,\infty)$ and every $N\in\N$.
%
\begin{lemma}  \label{l:second_moment_estimate_X}
  Assume~\ref{a:A1_N}, \ref{a:Hutzenthaler2009EJP} and~\ref{a:second_moments}.
  Then we have that
  \begin{equation} \label{eq:second_moment_estimate_X}
    \sup_{N\in\N}\sum_{i=1}^N \E\eckbb{\biggl(\sup_{t\in[0,T]}
        X_{t\wedge\tau_K^N}^{N}(i)\biggr)^2}<\infty
  \end{equation}
  for all $T\in[0,\infty)$ and all $K\in\N$.
\end{lemma}
\begin{proof}
  Fix $T\in[0,\infty)$ and $K\in\N$.
  Lemma 3.3 in~\cite{HutzenthalerWakolbinger2007} shows that, on the event $\{\tau_K^N\geq t\}$,
  $X_t^N(i)$
  is bounded from above by
  $Y_{t,0}^{N,K+N\mu_N(0)}$
  for all $t\in[0,T]$ almost surely for every $N\in\N$.
  By Assumption~\ref{a:A1_N} we have that
  $N\mu_N(0)\leq 2\theta$ for all $N\in\N$.
  Together with the second-moment estimate of
  Lemma~\ref{l:second_moment_estimate_Y},
  this implies that
  \begin{equation}  \begin{split} \label{nivea}
    \sum_{i=1}^N &\,\E\eckbb{\biggl(\sup_{t\in[0,T]}
        X_{t\wedge\tau_K^N}^{N}(i)\biggr)^2}
    \leq
    \sum_{i=1}^N \E\eckbb{\E^{X_0^N(i)}\eckbb{\biggl(\sup_{t\in[0,T]}
        Y_{t,0}^{N,K+2\theta}\biggr)^2}}\\
    &\leq C_T\sum_{i=1}^N
    \roB{\E \eckb{X_0^N(i)}+\E\eckB{\roB{X_0^N(i)}^2}
         +T\frac{K+2\theta}{N}+T\frac{(K+2\theta)^2}{N^2}}\\
    &\leq C_T\robb{\sup_{M\in\N}\E\Big[\sum_{i=1}^M X_0^M(i)\Big]
       +\sup_{M\in\N}\E\eckbb{\Big|\sum_{i=1}^M X_0^M(i)\Big|^2}+T (K+2\theta)+T(K+2\theta)^2}
  \end{split}     \end{equation}
  for every $N\in\N$ and some constant $C_T<\infty$.
  The right-hand side is finite due to
  Assumption~\ref{a:second_moments}.
\end{proof}
\subsection{Poisson limit of independent diffusions with vanishing immigration}
\label{ssec:Poisson_limit}
In this subsection, we prove~\eqref{eq:empirical_measure_converges} which is the central
step in the proof of Theorem~\ref{thm:convergence}.
Our proof is based on reversing time in the stationary process.
For the time reversal,
we consider the following stationary situation.
Excursions from zero of
the solution process $(Y_t)_{t\geq0}$ of the
SDE~\eqref{eq:Y} start at times given by the points
of an homogeneous Poisson point process on $\R$ with rate $1$.
This process of immigrating excursions is invariant for the dynamics
of $(Y_t)_{t\geq0}$ restricted to non-extinction, see~\eqref{eq:relation_m_Q}.
Now the time reversal of an excursion is again governed by the excursion measure,
see Lemma~\ref{l:time_reversal}.
As a consequence, reversing time in the process of immigrating excursions
does not change the distribution.

Let us retell this argument more formally. Consider a Poisson point process $\Pi$
on $(-\infty,\infty)\times U$ with intensity measure $ds\otimes Q$.
Then
$\sum_{(s,\eta)\in\Pi}\dl_{(\eta_{t-s})_{t\geq0}}$ is the process
of immigrating excursions. Note that at a fixed time $t$,
$\sum_{(s,\eta)\in\Pi}\dl_{\eta_{t-s}}$
is a Poisson point process
on $(0,\infty)$ with intensity measure
\begin{equation}  \label{eq:relation_m_Q}
  \int_{-\infty}^\infty Q(\eta_{t-s}\in\dup y)\,ds=
  \int_{-\infty}^0 Q(\eta_{-s}\in\dup y)\,ds=m(dy)
\end{equation}
where $m$ is the speed measure defined in~\eqref{eq:m}.
Here we used that $\eta_t=0$ for $t\in(-\infty,0]$ for all $\eta\in U$.
The relation~\eqref{eq:relation_m_Q} between the speed measure
and the excursion measure has been
established in Lemma 9.8 of~\cite{Hutzenthaler2009EJP}
by exploiting a well-known explicit formula for $\E^y\int_0^\infty f(Y_s)\,ds$.
It is also well-known (e.g.~(15.5.34) in~\cite{KarlinTaylor1981b})
that the speed measure $m$ is an
invariant measure for the sub-Markov semigroup
$\E^{\cdot}f(Y_t)\1_{Y_t>0}$.
This can also be seen from~\eqref{eq:relation_m_Q} by noting
that
$Q(\eta_s\in\dup y)$ is an entrance measure for this sub-Markov semigroup.
Thus the process of immigrating excursions is indeed invariant for the dynamics
of $(Y_t)_{t\geq0}$ restricted to non-extinction.

Now we show that reversing time in the process of immigrating excursions
does not change the distribution of the process.
The process $(Y_t)_{t\geq0}$ restricted to non-extinction is time-reversible
when started in the invariant measure $m$, that is,
\begin{equation} \label{eq:time_reversal_Y}
   \int_I \E^x F\rub{\ru{Y_t}_{t\leq T}}\1_{Y_T>0} m(dx)
  =\int_I \E^x F\rub{\ru{Y_{T-t}}_{t\leq T}}\1_{Y_T>0} m(dx)
\end{equation}
for every $T\in[0,\infty)$ and every non-negative measurable function
on $\C\rob{[0,T]}$,
see Section 13 of Chapter 15 in \cite{KarlinTaylor1981b}.
The next lemma shows that
if the speed measure $m$ is replaced by the left-hand side of~\eqref{eq:relation_m_Q},
then~\eqref{eq:time_reversal_Y} can be extended to allow for extinction.
First we state the Markov property of the excursion measure.
Definition~\eqref{eq:Q} of $Q$ as rescaled law of $(Y_t)_{t\geq0}$
together with the Markov property of $(Y_t)_{t\geq0}$ implies that
\begin{equation}  \begin{split}  \label{eq:Q_and_MP}
   \int F\ruB{\rub{\eta_t}_{0\leq t\leq T}}
                  \Ft\ruB{\rub{\eta_{T+t}}_{t\geq 0}}Q(d\eta)
  =\int F\ruB{\rub{\eta_t}_{0\leq t\leq T}}
       \E^{\eta_T}\Ft\ruB{\rub{Y_t}_{t\geq 0}}Q(d\eta)
\end{split}     \end{equation}
for all measurable functions $F,\Ft\colon\C\rob{[0,\infty),[0,\infty)}\to[0,\infty)$
satisfying $F(\uline{0})=0=\Ft(\uline{0})$
and every $T\in[0,\infty)$.
Here and below, $\uline{0}$ denotes the function which is $\equiv0$.
%
\begin{lemma}  \label{l:time_reversal}
  Assume~\ref{a:A1} and~\ref{a:Hutzenthaler2009EJP}.
  Then
  \begin{equation}  \label{eq:time_reversal_Q}
    \int \int_{-\infty}^\infty F\rub{\ru{\eta_{t-s}}_{t\in\R}}\dup s\, Q(\dup\eta)
   =\int \int_{-\infty}^\infty F\rub{\ru{\eta_{T-t-s}}_{t\in\R}}\dup s\, Q(\dup\eta)
  \end{equation}
  for all $T\in\R$ and all measurable functions
  $F\colon\C\rub{[0,\infty)}\to[0,\infty)$.
\end{lemma}
\begin{proof}
  It suffices (see e.g.~Theorem 14.12 in~\cite{Klenke2008})
  to establish~\eqref{eq:time_reversal_Q} for
  $F_n(\eta):=\prod_{i=1}^n f_i\ru{\eta_{t_i}}$ where $t_1<\ldots <t_n\in\R$
  and $f_1,\ldots,f_n\in\C_b\rub{[0,\infty),[0,\infty)}$.
  If $F_n(\uline{0})>0$, then both sides of~\eqref{eq:time_reversal_Q}
  are infinite.
  For the rest of the proof, we assume $F_n(\uline{0})=0$,
  that is, $f_i(0)=0$ for at least one $i\in\{1,\ldots,n\}$.
  We may even assume
  $f_i\in\C_c\rub{(0,\infty),[0,\infty)}$ for at least one $i\in\{1,\ldots,n\}$.
  Otherwise approximate $f_i$ monotonically from below with test functions
  which have compact support.
  In addition, we may without loss of generality assume
  $t_1=0=T$.
  Otherwise use a time translation.

  If $F_n$ vanishes on $\{\eta\colon\eta_0=0\text{ or }\eta_{t_n}=0\}$,
  then~\eqref{eq:time_reversal_Q} is essentially~\eqref{eq:time_reversal_Y}.
  To see this, consider
  \begin{equation}  \begin{split} \label{eq:have_2b_positive}
    &\int_{-\infty}^\infty\int \1_{\eta_{-s}>0} F_n\rub{\ru{\eta_{t-s}}_{t\in\R}}
       \1_{\eta_{t_n-s}>0}\, Q(\dup\eta)\dup s\\
    &\underset{\eqref{eq:Q_and_MP}}{=}
         \int_{-\infty}^\infty\int \1_{\eta_{-s}>0}\E^{\eta_{-s}}\eckB{F_n\rub{\ru{Y_{t}}_{t\in[0,t_n]}}
        \1_{Y_{t_n}>0}}\, Q(\dup\eta)\dup s\\
    &\underset{\eqref{eq:relation_m_Q}}{=}
       \int_I \E^{x}\eckB{F_n\rub{\ru{Y_{t}}_{t\in[0,t_n]}}
        \1_{Y_{t_n}>0}}m(\dup x).
  \end{split}     \end{equation}
  Applying~\eqref{eq:time_reversal_Y} with $T=t_n$, reversing 
  the calculation in~\eqref{eq:have_2b_positive} and substituting $s-t_n\mapsto s$
  shows that
  \begin{equation}  \begin{split} \label{eq:positive_paths}
   & \int_{-\infty}^\infty\int \1_{\eta_{-s}>0} F_n\rub{\ru{\eta_{t-s}}_{t\in[0,t_n]}}
       \1_{\eta_{t_n-s}>0}\, Q(\dup\eta)\dup s\\
   &= \int_{-\infty}^\infty \int
       \1_{\eta_{-s}>0}
       F_n\rub{\ru{\eta_{t_n-t-s}}_{t\in[0,t_n]}}
       \1_{\eta_{t_n-s}>0}
       \, Q(\dup\eta)\dup s.
       \\
   &= \int_{-\infty}^\infty \int \1_{\eta_{-t_n-s}>0}F_n\rub{\ru{\eta_{-t-s}}_{t\in[0,t_n]}}
       \1_{\eta_{-s}>0}\, Q(\dup\eta)\dup s.
  \end{split}     \end{equation}
  We prove~\eqref{eq:time_reversal_Q} with $F$ replaced by $F_n$ by induction
  on $n\in\N$.
  The base case $n=1$ follows from a time translation.
  The induction step $n-1\to n$ follows directly from~\eqref{eq:positive_paths}
  if $f_1(0)=0=f_n(0)$. If $f_n(0)>0$, then
  \begin{equation}  \begin{split}   \label{eq:last_point_zero}
    \lefteqn{\int \int_{-\infty}^\infty F_n\rub{\ru{\eta_{t-s}}_{t\in[0,t_n]}}
        \1_{\eta_{-s}>0}\1_{\eta_{t_n-s}=0}\,ds\, Q(d\eta)}\\
    &=f_n(0)\rubb{\int \int_{-\infty}^\infty F_{n-1}\rub{\ru{\eta_{t-s}}_{t\in[0,t_n]}}
        \1_{\eta_{-s}>0}\rub{1-\1_{\eta_{t_n-s}>0}}\,ds\, Q(d\eta)}\\
    &=f_n(0)\rubb{\int \int_{-\infty}^\infty F_{n-1}\rub{\ru{\eta_{-t-s}}_{t\in[0,t_n]}}
        \1_{\eta_{-s}>0}\rub{1-\1_{\eta_{-t_n-s}>0}}\,ds\, Q(d\eta)}\\
    &=\int \int_{-\infty}^\infty F_n\rub{\ru{\eta_{-t-s}}_{t\in[0,t_n]}}
        \1_{\eta_{-s}>0}\1_{\eta_{-t_n-s}=0}\,ds\, Q(d\eta).
  \end{split}     \end{equation}
  For the second step we used linearity, applied the induction hypothesis
  and equation~\eqref{eq:positive_paths} and again used linearity.
  Adding~\eqref{eq:positive_paths} and~\eqref{eq:last_point_zero}
  proves the induction step in case of $f_1(0)=0$.
  The remaining case $f_1(0)>0$ follows from a similar calculation
  as in~\eqref{eq:last_point_zero}.
  This completes the proof of Lemma~\ref{l:time_reversal}.
\end{proof}
%
\begin{lemma} \label{l:speed_excursion_measure}
  Assume~\ref{a:A1} and~\ref{a:Hutzenthaler2009EJP}.
  Then the solution process $(Y_t)_{t\geq0}$ of the SDE~\eqref{eq:Y} satisfies that
  \begin{equation}
    \int_I \E^x F\rub{\ru{Y_t}_{t\in[0,T]}}m(dx)
    =\int \int_{-\infty}^T F\rub{\ru{\eta_{T-t-s}}_{t\in[0,T]}}\dup s\, Q(\dup\eta)
  \end{equation}
  for all measurable functions
  $F\colon\C\rub{[0,T],I}\to[0,\infty)$ satisfying $F(\uline{0})=0$ and all $T\in[0,\infty)$.
\end{lemma}
\begin{proof}
  Express the speed measure in terms of the excursion measure as in~\eqref{eq:relation_m_Q}
  \begin{equation}  \begin{split}
    \int_I \E^x F\rub{\ru{Y_t}_{t\in[0,T]}}m(dx)
    &=\int_{-\infty}^0 \int \1_{\eta_{-s}>0}\E^{\eta_{-s}}F\rub{\ru{Y_t}_{t\in[0,T]}}Q(d\eta)\,ds\\
    &= \int \int_{-\infty}^\infty\1_{\eta_{-s}>0}
      F\rub{\ru{\eta_{t-s}}_{t\in[0,T]}}\dup s\, Q(\dup\eta)\\
    &=\int \int_{-\infty}^T F\rub{\ru{\eta_{T-t-s}}_{t\in[0,T]}}Q(\dup\eta)\, \dup s.
  \end{split}     \end{equation}
  The last two steps are the Markov property~\eqref{eq:Q_and_MP}
  and Lemma~\ref{l:time_reversal}, respectively.
\end{proof}
With Lemma~\ref{l:speed_excursion_measure} in hand, we now
reverse time to prove a first version of the
Poisson approximation~\eqref{eq:empirical_measure_converges}.
%
\begin{lemma}   \label{l:vanishing_immigration}
  Assume~\ref{a:A1_N} and~\ref{a:Hutzenthaler2009EJP}.
  Let $c,s\in[0,\infty)$ be real numbers.
  Let $(Y_{t,s}^{N,c})_{t\geq s}$ be the solution of~\eqref{eq:YN}
  with $\zeta_N(\cdot)\equiv c$ and $Y_{s,s}^{N,c}=0$ for every $N\in\N$ with $N\geq c/|I|$.
  Then
  \begin{equation}  \label{eq:vanishing_immigration}
    \limN N\E^0\big[ f_N(Y_{t,s}^{N,c})\big]
       = c\int_s^t\int f_0(\chi_{t-r})Q(d\chi)\,dr
  \end{equation}
  for all $t\in[s,\infty)$ and
  all functions $f_N\colon I\to\R$, $N\in\N_0$,
  with $\sup_{x\neq y\in I}\sup_{N\in\N}\tfrac{|f_N(x)-f_N(y)|}{|x-y|}<\infty$,
  with $\lim_{N\to\infty}N|f_N(0)|=0$
  and
  with $\lim_{N\to\infty}f_N(y)=f_0(y)$ for all $y\in I$.
\end{lemma}
\begin{proof}
  If $c=0$, then both sides of~\eqref{eq:vanishing_immigration} are equal to zero.
  So for the rest of the proof we assume $c>0$.
  Let $\eps,\dl_0\in(0,|I|\wedge\alpha)$ be such that
  $\sigma_N^2(y)\geq \eps\sigma^2(y)$ for all $y\in[0,\dl_0]$
  and such that $\int_0^{\dl_0}y/\sigma^2(y)\,dy<\infty$.
  Lemma~\ref{l:order_of_m_N} provides us with $N/m_N\rub{(0,\dl)}\to c$
  as $N\to\infty$
  for all $\dl\in(0,\dl_0)$.
  Thus
  \begin{equation}  \begin{split}  \label{eq:sancho}
    \limN N\E^0 \big[f_N(Y_{t,s}^{N,c})\big]
    &=\limN\frac{N}{m_N\rub{(0,\dl)}}
      \limN\int_0^\dl \E^0 \big[f_N(Y_{t,s}^{N,c})\big]m_N(dx)\\
   &=:c\limN\int_0^\dl \E^x \big[f_N(Y_{t,s}^{N,c})\big]m_N(dx)+C(\dl)
  \end{split}     \end{equation}
  for all $\dl\in(0,\dl_0)$ and all $t\in[s,\infty)$
  where $C\colon[0,\dl_0]\to\R$ is a suitable function.
  The term $C(\dl)$ converges to zero as $\dl\to0$ according
  to Lemma~\ref{l:replace_initial_point}.

  The speed measure $m_N$ is an invariant (non-probability) measure
  for $(Y_{t,s}^{N,c})_{t\geq0}$, see e.g.~(15.5.34) in~\cite{KarlinTaylor1981b}.
  Thus we may reverse time.
  As we let $N\to\infty$, we will exploit that $(Y_{t,s}^{N,c})_{t\geq0}$ converges
  weakly to
  the solution process $(Y_t)_{t\geq s}$ of the SDE~\eqref{eq:Y}.
  In addition, $m_N(dy)$ converges vaguely to
  $m(dy)$ as $N\to\infty$
  due to the dominated convergence theorem and Assumptions~\ref{a:A1_N}
  and~\ref{a:Hutzenthaler2009EJP}
  as the densities converge.
  These observations imply that
  \begin{equation}  \begin{split}  \label{eq:pancho}
  \lefteqn{
    \limN\int_0^\dl \E^x\big[ f_N\ru{Y_{t,s}^{N,c}}\big]m_N(dx)
    =\limN\int_I f_N(y) \P^y\eckb{Y_{t,s}^{N,c}\leq\dl} m_N(dy)
  }
    \\ &
    =\int_I f_0(y) \P^y\eckb{ Y_{t-s}\leq\dl} m(dy)
    =\int\int_{-\infty}^t f_0\rub{\chi_{t-r}}\1_{\chi_{s-r}\leq\dl}dr\,Q(d\chi)
    \lradlO\int\int_{s}^t f_0\rub{\chi_{t-r}}dr\,Q(d\chi).
  \end{split}     \end{equation}
  The last but one step is Lemma~\ref{l:speed_excursion_measure} and the
  last step follows from the dominated convergence theorem together with
  Lemma~\ref{l:finite_excursion_area}.
  Putting~\eqref{eq:sancho} and~\eqref{eq:pancho} together completes the proof
  of Lemma~\ref{l:vanishing_immigration}.
\end{proof}

Next we use induction to generalize Lemma~\ref{l:vanishing_immigration} to
test functions which depend on finitely many time coordinates.
For this
let $\MCE_{s,T}$ be the following set of bounded functions on $\C\rob{[s,T],I}$
for $0\leq s\leq T<\infty$
which
depend on finitely many coordinates and which are globally Lipschitz continuous
in every coordinate
\begin{equation}  \begin{split}  \label{eq:MCET}
  \MCE_{s,T}&:=\biggl\{\C\rob{[s,T],I}\ni\eta\mapsto \prod_{i=1}^n f_i(\eta_{t_i})\in\R\colon n\in\N,
  s\leq t_1<\ldots< t_n\leq T,
  \\
  &\qquad\qquad\qquad\qquad\qquad f_1,...,f_n\in\C\rob{I,\R}\text{ are bounded and globally Lipschitz continuous}\}
\end{split}     \end{equation}
for every $T\in[s,\infty)$ and every $s\in[0,\infty)$.
  Due to the Lipschitz continuity and boundedness 
  of $f_1,\ldots,f_n$,
  there exists a constant $L_F\in(0,\infty)$ such that
  \begin{equation}  \label{eq:F_is_globally_Lipschitz}
    \abs{F\ru{\eta}-F(\etab)}\leq L_F\sum_{j=1}^n\abs{\eta_{t_j}-\etab_{t_j}}
    \qquad \fa \eta,\etab\in \C\rob{[s,T],I}
  \end{equation}
  for all $F\in\MCE_{s,T}$ and all $0\leq s\leq T<\infty$.
  Note that the set $\MCE_{s,T}$ is closed under multiplication and separates points
  in $\C([s,T],I)$
  for all $0\leq s\leq T<\infty$.
  Thus the linear span of $\MCE_{s,T}$ is an algebra which separates points
  in $\C([s,T],I)$
  for all $0\leq s\leq T<\infty$.
  According to Theorem 3.4.5 in~\cite{EthierKurtz1986} the linear span of $\MCE_{s,T}$
  is distribution determining for measures on
  $\C([s,T],I)$
  and so is $\MCE_{s,T}$
  for all $0\leq s\leq T<\infty$.
%
\begin{lemma}   \label{l:vanishing_immigration_general}
  Assume~\ref{a:A1_N} and~\ref{a:Hutzenthaler2009EJP}.
  Let $0\leq s\leq T<\infty$.
  Suppose that $\zeta\colon[s,\infty)\to[0,\infty)$ 
  and that $\zetah_N\colon[s,\infty)\to N\mal I$
  are square Lebesgue integrable
  and that $\int_s^T |\zeta(r)-\zetah_N(r)|dr\to 0$ as $N\to\infty$.
  Let $(Y_{t,s}^{N,\zetah})_{t\in[s,\infty)}$ satisfy~\eqref{eq:YN}.
  Then
  \begin{equation}  \label{eq:vanishing_immigration_general}
     \limN N\E^0 F\rub{(Y_{t,s}^{N,\zetah})_{t\in[s,T]}}=
      \int_s^\infty \zeta(r)\int F\rub{(\chi_{t-r})_{t\in[s,T]}}Q(d\chi)\,dr\in\R
  \end{equation}
  for
  all functions $F\in\MCE_{s,T}$.
\end{lemma}
\begin{proof}
  Let $L_F\in(0,\infty)$ be such that $F$ satisfies~\eqref{eq:F_is_globally_Lipschitz}
  and let $F$ be bounded by $C_F\in(0,\infty)$.
  Moreover,  let
  $(Y_t)_{t\in[s,\infty)}$ be the pathwise unique stochastic process
  such that
  $(Y_t)_{t\in[s,\infty)}$ 
  and $(Y_{t,s}^{N,\zetah})_{t\in[s,\infty)}$ are solutions of~\eqref{eq:YN}
  with respect to the same Brownian motion.
  Lemma~\ref{l:first_moment_estimate} implies that
  \begin{equation}  \begin{split}
    \lefteqn{
    \absB{\limN N\E^0  F\roB{\rob{Y_{t,s}^{N,\zetah}}_{t\in[s,T]}}
        - \limN N\E^0  F\roB{\rob{Y_{t,s}^{N,\zeta}}_{t\in[s,T]}}
         }
    }\\
    &\leq\limN N L_F\sum_{j=1}^n\E^0 
      \Big|Y_{t_j,s}^{N,\zetah}-Y_{t_j,s}^{N,\zeta}
      \Big|\\
    &\leq L_T\cdot n e^{L_\mu T}\limN\int_0^T\absb{\zetah_N(r)-\zeta(r)} dr =0.
  \end{split}     \end{equation}
  Therefore, it suffices to prove~\eqref{eq:vanishing_immigration_general}
  with $\zetah_N$ replaced by $\zeta$.
  A similar argument shows that we may assume $\zeta$ to be bounded;
  otherwise replace $\zeta$ by $\min(\zeta,M)$ and let $M\to\infty$.

  We begin with the case of $\zeta$ being a simple function.
  W.l.o.g.\ we consider $\zeta(\cdot)=\sum_{i=1}^n c_i\1_{[t_{i-1},t_i)}(\cdot)$
  where $c_1,...,c_n\geq0$ and $t_0=s$
  as we may let $F$ depend trivially on further time points.
  The proof of~\eqref{eq:vanishing_immigration_general}
  is by induction on $n$.
  The case $n=1$ has been settled in 
  Lemma~\ref{l:vanishing_immigration}.
  For the induction step we split up the left-hand side
  of~\eqref{eq:vanishing_immigration_general} into two terms according
  to whether the process at time $t_1$ is essentially zero or not.
  In order to formalize the notion ``essentially zero'',
  choose a function
  $\phi_\dl\in\C^2\rub{I,[0,1]}$ such that $\phi_\dl(x)=1$
  for $x\geq2\dl$ and $\phi_\dl(x)=0$ for $x\leq\dl$ for every $\dl\in(0,|I|)$.
  Furthermore
  define $\Fb_2(\eta):=\prod_{i=2}^n f_i(\eta_{t_i})$
  for all $\eta\in\C([t_1,T],I)$.

  First we consider the case that the process is away from $0$ at time $t_1$.
  The following equation~\eqref{eq:no_immig} shows that we may discard immigration
  after time $t_1$. 
  For this, note that the moment estimate of Lemma \ref{l:first_moment_estimate} implies that
  \begin{equation}  \begin{split}
     &\absbb{
      \absB{\E^{y}\Fb_2\ruB{\rub{Y_{t,t_1}^{N,\zeta}}_{t\in[t_1,T]}}
         -\E^{y}\Fb_2\ruB{\rub{Y_{t-t_1}}_{t\in[t_1,T]}}}
      -
      \absB{\E^{z}\Fb_2\ruB{\rub{Y_{t,t_1}^{N,\zeta}}_{t\in[t_1,T]}}
         -\E^{z}\Fb_2\ruB{\rub{Y_{t-t_1}}_{t\in[t_1,T]}}}
     }
     \\
     &\leq
      \absB{
         \E^{y}\Fb_2\ruB{\rub{Y_{t,t_1}^{N,\zeta}}_{t\in[t_1,T]}}
         -
         \E^{z}\Fb_2\ruB{\rub{Y_{t-t_1}^{N,\zeta}}_{t\in[t_1,T]}}
      }
      +
      \absB{
         \E^{y}\Fb_2\ruB{\rub{Y_{t-t_1}}_{t\in[t_1,T]}}
         -
         \E^{z}\Fb_2\ruB{\rub{Y_{t-t_1}}_{t\in[t_1,T]}}
      }
    \\
    &\leq
    2L_F\sum_{j=2}^n e^{L_\mu(t_j-s)}|y-z|
  \end{split}     \end{equation}
  for all $y,z\in I$ and all $N\in\N$.
  Consequently, the sequence of functions 
  \begin{equation}
    I\ni y\mapsto \phi_\dl(y)f_1(y)
      \absB{\E^{y}\Fb_2\ruB{\rub{Y_{t,t_1}^{N,\zeta}}_{t\in[t_1,T]}}
         -\E^{y}\Fb_2\ruB{\rub{Y_{t-t_1}}_{t\in[t_1,T]}}}
  \in\R, \quad N\in\N,
  \end{equation}
  is uniformly globally Lipschitz continuous
  and satisfies that
  \begin{equation}  \label{eq:which.fN}
      \lim_{N\to\infty}\phi_\dl(y)f_1(y)
      \absB{\E^{y}\Fb_2\ruB{\rub{Y_{t,t_1}^{N,\zeta}}_{t\in[t_1,T]}}
         -\E^{y}\Fb_2\ruB{\rub{Y_{t-t_1}}_{t\in[t_1,T]}}}
      =0
  \end{equation}
  for all $y\in I$
  and
  all $\dl\in(0,|I|)$.
  Lemma~\ref{l:vanishing_immigration} thus implies that
  \begin{equation}  \begin{split}  \label{eq:no_immig}
    \lim_{N\to\infty}N\E^0\eckbb{\phi_\dl\rub{Y_{t_1,s}^{N,c_1}}
     f_1\rub{Y_{t_1,s}^{N,c_1}}\cdot
      \absB{\E^{Y_{t_1,s}^{N,c_1}}\Fb_2\ruB{\rub{Y_{t,t_1}^{N,\zeta}}_{t\in[t_1,T]}}
         -\E^{Y_{t_1,s}^{N,c_1}}\Fb_2\ruB{\rub{Y_{t-t_1}}_{t\in[t_1,T]}}}}
    =0
  \end{split}     \end{equation}
  for all $\dl\in(0,|I|)$.
  If the process is essentially zero at time $t_1$
  (that is $1-\phi_\dl\rub{Y_{t_1,s}^{N,\zeta}}\approx 1$),
  then we may restart the process at time $t_1$ in the state $0$.
  The Lipschitz continuity of $f_1,\ldots,f_n$ together with the moment estimate
  of Lemma~\ref{l:first_moment_estimate} provides us with
  \begin{equation}  \begin{split}  \label{eq:smaller_dl}
    \limsupN &N\E^0\eckbb{\rub{1-\phi_\dl\rub{Y_{t_1,s}^{N,c_1}}}
     \absB{
       f_1\rub{Y_{t_1,s}^{N,c_1}}\E^{Y_{t_1,s}^{N,c_1}}\Fb_2\ruB{\rub{Y_{t,t_1}^{N,\zeta}}_{t\in[t_1,T]}}
       -f_1(0)\E^{0}\Fb_2\ruB{\rub{Y_{t,t_1}^{N,\zeta}}_{t\in[t_1,T]}}}}\\
    &\leq \limsupN N
       \E^0\eckB{\rub{1-\phi_\dl\rub{Y_{t_1,s}^{N,c_1}}}n L_F e^{L_\mu t_n} Y_{t_1,s}^{N,c_1}}\\
    &=n L_F e^{L_\mu t_n} c_1 \int_s^{t_1}\int\eckB{\ruB{1-\phi_\dl(\chi_{t_1-r})}\chi_{t_1-r}}Q(d\chi)\,dr
  \end{split}     \end{equation}
  for all $\dl\in(0,|I|)$.
  The last step is Lemma~\ref{l:vanishing_immigration}.
  By the dominated convergence
  theorem together with Lemma~\ref{l:finite_excursion_area},
  the right-hand side of~\eqref{eq:smaller_dl}
  converges to zero as $\dl\to0$.
  Therefore, using~\eqref{eq:no_immig}, \eqref{eq:smaller_dl},
  the Markov property,
  Lemma~\ref{l:vanishing_immigration},
  $\lim_{N\to\infty}\E^0 \big[\phi_\dl(Y_{t_1,s}^{N,c_1})\big]=0$ for all $\dl>0$
  and applying the induction hypothesis
  \begin{equation}  \begin{split}
    \lefteqn{\limN N\E^0\eckB{ F\ruB{\rub{Y_{t,s}^{N,\zeta}}_{t\in[s,T]}}}}\\
    &=
    \limN N\E^0\eckB{\ruB{
      \phi_\dl\rub{Y_{t_1,s}^{N,c_1}}
      +1-\phi_\dl\rub{Y_{t_1,s}^{N,c_1}}}
      F\ruB{\rub{Y_{t,s}^{N,\zeta}}_{t\in[s,T]}}}\\
    &=\limdlO\limN N\E^0\eckB{\phi_\dl\rub{Y_{t_1,s}^{N,c_1}}
      f_1\rub{Y_{t_1,s}^{N,c_1}}
     \E^{Y_{t_1,s}^{N,c_1}}\Fb_2\ruB{\rub{Y_{t-t_1}}_{t\in[t_1,T]}}}\\
    &\qquad+\limdlO\limN\E^0\eckB{1-\phi_\dl\rub{Y_{t_1,s}^{N,c_1}}}
       f_1(0)\limN N
       \E^0 \eckB{\Fb_2\ruB{\rub{Y_{t,t_1}^{N,\zeta}}_{t\in[t_1,T]}}}\\
    &=\limdlO\int_s^{t_1}c_1\int\phi_\dl\rub{\chi_{t_1-r}}f_1\rub{\chi_{t_1-r}}
      \E^{\chi_{t_1-r}}\eckB{\Fb_2\ruB{\rub{Y_{t-t_1}}_{t\in[t_1,T]}}}Q(d\chi)\,dr\\
    &\qquad+f_1(0)\int_{t_1}^{\infty}\zeta(r)\int \Fb_2\ruB{\ru{\chi_{t-r}}_{t\in[t_1,T]}}Q(d\chi)\,dr\\
    &=\int_s^\infty \zeta(r)\int F\ruB{(\chi_{t-r})_{t\in[s,T]}}Q(d\chi)\,dr.
  \end{split}     \end{equation}
  The last step follows from
  the pointwise convergence $\phi_\dl(x)\to 1$ as $\dl\to0$ for every $x\in(0,|I|)$
  together with the dominated convergence theorem
  and from the Markov property~\eqref{eq:Q_and_MP}.
  Finally let $\zeta$ be integrable. Approximate $\zeta$ with simple
  functions $(\zeta_n)_{n\in\N}$.
  Applying Lemma~\ref{l:first_moment_estimate},
  it is straight forward to show that
  equation~\eqref{eq:vanishing_immigration_general} with $\zeta$ replaced by $\zeta_n$
  converges to equation~\eqref{eq:vanishing_immigration_general} as $n\to\infty$.
\end{proof}

%
\begin{lemma}  \label{l:vanishing_immigration_weak_process}
  Assume~\ref{a:A1_N} and~\ref{a:Hutzenthaler2009EJP}.
  Fix $s\in[0,\infty)$.
  Suppose that
  $\zeta\colon[s,\infty)\to[0,\infty)$ and
  $\zetah_N\colon[s,\infty)\to N\cdot I$ are locally square 
  integrable functions for $N\in\N$
  and that
  $\zetah_N\to\zeta$ as $N\to\infty$ in $L^1_{loc}$.
  Let $(Y_{t,s}^{N,\zetah}(i))_{t\in[s,\infty)}$, $i\leq N$, be independent
  solutions of~\eqref{eq:YN}
  satisfying  $Y_{s,s}^{N,\zetah}(i)=0$ for every $N\in\N$.
  Then we have that
  \begin{equation} \label{eq:vanishing_immigration_weak_process}
     \Bigl(\sum_{i=1}^{N}
      \dl_{Y_{t,s}^{N,\zetah}(i)}\Bigr)_{s\leq t\leq T}\wlra
      \Bigl(\int\dl_{\eta_{t-u}}\Pi(du,d\eta)\Bigr)_{s\leq t\leq T}\qqasN
  \end{equation}
  where $\Pi$ is a Poisson point process on
  $[s,\infty)\times U$
  with intensity measure
     $\zeta(u) du\otimes Q(d\eta)$.
  Moreover let $\Fb$ be a continuous function
  from $\C\rob{[s,\infty),\R}$ to $\R$ satisfying the Lipschitz condition
  \begin{equation}  \label{eq:F_Lipschitz}
    \absb{\Fb\rob{\eta}-\Fb\rob{\etab}}
    \leq L_\Fb\sum_{j=1}^n \absb{\eta_{t_j}-\etab_{t_j}}
     \qquad\fa\eta, \etab\in\C\rob{[s,\infty),\R}
  \end{equation}
  for some $s\leq t_1\leq\cdots\leq t_n\leq T$, some $n\in\N$
  and some $L_\Fb\in(0,\infty)$.
  In addition let
  $\fb\colon I\to\R$ be a continuous function satisfying
  $\abs{\fb(x)}\leq L_\fb x$ for all $x\in I$.
  Then we have that
  \begin{equation}  \begin{split}  \label{eq:convY}
    \limN\E\eckbb{\Fb\robb{\Bigl(\sum_{i=1}^N
           \fb\rob{Y_{t,s}^{N,\zetah}(i)}\Bigr)_{t\in[s,\infty)}}}
    =\E\eckbb{\Fb\robb{\Bigl(\int \fb\rob{\eta_{t-u}}\Pi(du,d\eta)\Bigr)_{t\in[s,\infty)}}}.
  \end{split}     \end{equation}
\end{lemma}
\begin{proof}
  Let $\MCF\subseteq\C_c^2\rob{I\cap(0,\infty),\R}$ be a countable dense subset of
  $\C_c\rob{I\cap(0,\infty),\R}$.
  Tightness of
  \begin{equation}  \label{eq:tightnessTOshow}
     \curlbb{\Bigl(\sum_{i=1}^{N}
      \dl_{Y_{t,s}^{N,\zetah}(i)}\Bigr)_{t\in[s,\infty)}\colon N\in\N}
  \end{equation}
  follows from tightness of
  \begin{equation}   \label{eq:tightnessfromhere}
     \curlbb{\Bigl(\sum_{i=1}^{N}
      f\rob{Y_{t,s}^{N,\zetah}(i)}\Bigr)_{t\in[s,\infty)}\colon N\in\N}
      \quad\fa f\in\MCF.
  \end{equation}
  This type of argument has been established in
  Theorem 2.1 of Roelly-Coppoletta \nocite{RoellyCoppoletta1986}(1986)
  for the weak topology and $\C_0$.
  Following the proof hereof, one can show the analogous argument 
  for the vague topology and $\C_c$.
  Fix $f\in\MCF$ and define $S_{t,s}^N:=\sum_{i=1}^N f\rob{Y_{t,s}^{N,\zetah}(i)}$ for all
  $t\in[s,\infty)$ and $N\in\N$.
  Note that $f$ is globally Lipschitz continuous.
  For $K\in\N$ and fixed $t\in[s,\infty)$, global Lipschitz continuity of $f$
  implies that
  \begin{equation}
    \P\eckB{\absb{S_{t,s}^N}\geq K}
    \leq \frac{L_f}{K}
         \E\eckbb{\sum_{i=1}^N Y_{t,s}^{N,\zetah}(i)}
    \leq \frac{L_f}{K}
         \sup_{M\in\N}M\E^0\eckbb{Y_{t,s}^{M,\zetah}(1)}
  \end{equation}
  for some constant $L_f<\infty$ and for all $N\in\N$.
  The right-hand side is finite according to Lemma~\ref{l:first_moment_estimate}
  and converges to zero
  as $K\to\infty$.
  This proves tightness of $S_{t,s}^N$, $N\in\N$, for every $t\in[s,\infty)$.
  For the second part of the Aldous criterion, fix $T\in[s,\infty)$
  and let $\tau_N$, $N\in\N$, be stopping times which values in $[s,T]$.
  In addition define 
  \begin{equation}
    \Gen_t^N f(x):=\roB{\frac{\zetah_N(t)}{N}-x+\mut_N(x)}f^{'}(x)
     +\frac{1}{2}\sigma_N^2(x)f^{''}(x)\qquad \fa x\in I,
     t\in[s,\infty),N\in\N.
  \end{equation}
  The functions $\mut_N$, $\sigma_N$ and $\sigma_N^2$ are uniformly
  globally Lipschitz continuous on the support of $f$
  according to Assumption~\ref{a:A1_N}.
  Therefore there exists a constant
  $C_f\in[1,\infty)$ such that
  $\abs{\Gen_t^N f(x)}\leq C_f\rob{ \tfrac{\zetah_N(t)}{N}+x}$
  and $\rob{f^{'}\sigma_N}^2(x)\leq C_f^2 x$ for all $x\in I$, $t\in[s,\infty)$
  and all $N\in\N$.
  For fixed $\eta>0$ and $\dlb\in[0,1]$, we use It\^o's formula
  to obtain that
  \begin{equation}  \begin{split} \label{eq:esti_tightiY}
    &\eta^2\P\eckB{\absb{S_{\tau_N+\dl,s}^N-S_{\tau_N,s}^N}>\eta}
    \leq\E\eckB{\Bigl(S_{\tau_N+\dl,s}^N-S_{\tau_N,s}^N\Bigr)^2}\\
    &=\E\eckbb{\Bigl(\sum_{i=1}^N\int_{\tau_N}^{\tau_N+\dl}
            \Gen_{u}^N f\rob{Y_{u,s}^{N,\zetah}(i)}\,du
            +\sum_{i=1}^N\int_{\tau_N}^{\tau_N+\dl}
                    \rob{f^{'}\mal\sigma_N}\rob{Y_{u,s}^{N,\zetah}(i)}dB_u(i)\Bigr)^2}\\
    &\leq    3C_f^2\E\eckbb{\Bigl(\int_{\tau_N}^{\tau_N+\dl}\zetah_N(u)du\Bigr)^2}
            +3C_f^2\E\eckbb{\Bigl(\sum_{i=1}^N
              \int_0^\dl Y_{\tau_N+u,s}^{N,\zetah}(i) du \Bigr)^2}\\
    &\qquad\quad\qquad\qquad +3\sum_{i=1}^N\E\eckB{\int_{0}^{\dl}
                    \rob{f^{'}\sigma_N}^2\rob{Y_{\tau_N+u,s}^{N,\zetah}(i)}du}\\
    &\leq    3C_f^2\dl\E\eckbb{\int_{\tau_N}^{\tau_N+\dl}
            \bigl(\zetah_N(u)\bigr)^2du}
            +3C_f^2\dl\E\eckbb{\int_0^\dl \Bigl(\sum_{i=1}^N Y_{\tau_N+u,s}^{N,\zetah}(i)
                            \Bigr)^2 du}\\
    &\qquad\quad\qquad\qquad +3C_f^2\E\eckB{\int_{0}^{\dl}
                   \sum_{i=1}^N Y_{\tau_N+u,s}^{N,\zetah}(i) du}\\
    &\leq \dlb\mal 3C_f^2\sup_{M\in\N}\int_s^{T+1}\roB{\zetah_M(u)}^2 du
    +\dlb\mal 6C_f^2\sup_{M\in\N}
      \E^0\eckbb{\sup_{s\leq u\leq T+1}\Bigl(\sum_{i=1}^M Y_{u,s}^{M,\zetah}(i)\Bigr)^2}+\dlb 3C_f^2
  \end{split}     \end{equation}
  for all $\dl\leq \dlb$ and  for all $N\in\N$.
  The right-hand side of~\eqref{eq:esti_tightiY}
  is finite by Lemma~\ref{l:Y_2_supsum}.
  Letting $\dlb\to0$, the left-hand side of~\eqref{eq:esti_tightiY}
  converges to zero uniformly in $N\in\N$ and $\dl\leq\dlb$.
  This proves tightness of~\eqref{eq:tightnessTOshow} according to
  the Aldous criterion, see Aldous (1978)\nocite{Aldous1978}.

  Next we prove convergence of fi\-nite-di\-men\-sio\-nal distributions.
  Let $n\in\N$, $f\in\MCF$ with $f\geq 0$, $s\leq t_1\leq\cdots\leq t_n$ and $\ld_1,\cdots,\ld_n\in[0,\infty)$
  be arbitrary.
  Using independence we obtain that
  \begin{equation}  \begin{split}  \label{eq:fdd_total_mass}
    \E\eckB{ e^{-\sum_{j=1}^{n} \ld_j S_{t_j,s}^N}}
    &= \E \eckB{e^{-\sum_{i=1}^N\sum_{j=1}^{n} \ld_j f\rob{Y_{t_j,s}^{N,\zetah}(i)}}}
    =\prod_{i=1}^N\eckB{ \E e^{-\sum_{j=1}^{n} \ld_j f\rob{Y_{t_j,s}^{N,\zetah}(i)}}}\\
    &=\biggl(1-\frac{N \E^0\eckB{ 1-e^{-\sum_{j=1}^{n}
           \ld_j f\rob{Y_{t_j,s}^{N,\zetah}(1)}}}}{N}\biggr)^N
  \end{split}     \end{equation}
  for all $N\in\N$.
  Note that 
  $
    1-e^{-\sum_{j=1}^n\ld_j f(x_j)}
    =\sum_{j=1}^n \rob{1- e^{-\ld_j f(x_j)}}e^{-\sum_{i=1}^{j-1}\ld_i f(x_i)}
  $
  for all $x_1,\ldots,x_n\in[0,\infty)$.
  Applying Lemma~\ref{l:vanishing_immigration_general}
  to each summand of this telescope sum we get that
  \begin{equation}  \begin{split}
    \limN\E\eckB{ e^{-\sum_{j=1}^{n} \ld_j S_{t_j,s}^N}}
    &
    = \exp\robb{-\lim_{N\to\infty}
           N \E^0\eckB{ 1-e^{-\sum_{j=1}^{n} \ld_j f\rob{Y_{t_j,s}^{N,\zetah}(1)}}}
           }
    \\ &
    = \exp\robb{-\int\int_s^\infty\zeta(r)
           \roB{1-e^{-\sum_{j=1}^{n} \ld_j f\rob{\eta_{t_j-r}}}}dr Q(d\eta)}
    \\ &
    = \E\eckbb{ \exp\robb{-\sum_{j=1}^{n}\ld_j  \int f\roB{\eta_{t_j-u}}\Pi(du,d\eta)}}.
  \end{split}     \end{equation}
  This proves convergence of fi\-nite-di\-men\-sio\-nal distributions.

  It remains to
  prove~\eqref{eq:convY}, which includes non-bounded test functions.
  Let $\Fb$ and $\fb$ as in~\eqref{eq:convY}.
  By the previous step and by 
  the Skorokhod representation of weak convergence
  (e.g.\ Theorem II.86.1 in~\cite{RogersWilliams2000a})
  there exists a version of 
  $
    \curlb{\bigl(\sum_{i=1}^{N}
      \dl_{Y_{t,s}^{N,\zetah}(i)}
            \bigr)_{t\in[s,\infty)}
         \colon N\in\N
        }
  $
  which converges almost surely as $N\to\infty$.
  For every $K\in\N$, let $h_K\colon[0,\infty)\to[0,1]$
  be a continuous function
  which satisfies $h_K(x)=1$ for every $x\in[\tfrac{1}{K},K]$ and $h_K(x)=0$ for every
  $x\not\in[\frac{1}{2K},2K]$.
  Then $\fb\mal h_K\in\C_c\rob{I\cap(0,\infty),\R}$ for all $K\in\N$. Thus
  \begin{equation} \label{eq:conv_almostsurely}
    \curlbb{M\wedge \Fb\robb{\biggl(\sum_{i=1}^{N}
            \rob{\fb\mal h_K}\rob{  Y_{t,s}^{N,\zetah}(i) }
          \biggr)_{t\in[s,\infty)}
                  }
         \colon N\in\N
        }
  \end{equation}
  converges almost surely as $N\to\infty$ for all $M,K\in\N$.
  Next we observe that
  \begin{equation}  \begin{split}
    &
    \varlimsup_{N\to\infty}
    \biggl|
       \E\eckbb{\Fb\robb{\Bigl(\sum_{i=1}^N
          \fb\rob{Y_{t,s}^{N,\zetah}(i)}\Bigr)_{t\in[s,\infty)}}}
       -
       \E\eckbb{M\wedge\Fb\robb{\Bigl(\sum_{i=1}^N
          \fb\rob{Y_{t,s}^{N,\zetah}(i)}\Bigr)_{t\in[s,\infty)}}}
    \biggr|
    \\
    &=
    \varlimsup_{N\to\infty}
    \biggl|
       \E\eckbb{
          \Fb\robb{\Bigl(\sum_{i=1}^N
          \fb\rob{Y_{t,s}^{N,\zetah}(i)}\Bigr)_{t\in[s,\infty)}}
          \1_{\Fb\big((\sum_{i=1}^N
          \fb(Y_{t,s}^{N,\zetah}(i)))_{t\in[s,\infty)}\big)>M}
          }
    \biggr|
    \\
    &\leq
    \sup_{N\in\N}\frac{1}{M}\E\eckbb{\biggl(\Fb\robb{\Bigl(\sum_{i=1}^N
          \fb\rob{Y_{t,s}^{N,\zetah}(i)}\Bigr)_{t\in[s,\infty)}}\biggr)^2}
    \\ &
    \leq
    \frac{1}{M}2\Fb\rob{\uline{0}}+
    \frac{1}{M}\sup_{N\in\N}2L_\Fb^2 L_\fb^2\E\eckbb{\biggl(\sum_{j=1}^n \sum_{i=1}^N
       Y_{t_j,s}^{N,\zetah}(i)\biggr)^2}\\
    &\leq 
    \frac{1}{M}2\Fb\rob{\uline{0}}+
    \frac{1}{M}2L_\Fb^2 L_\fb^2 n^2\sup_{N\in\N}\E\eckbb{\sup_{t\in[s,t_n]}
          \Big(\sum_{i=1}^N Y_{t,s}^{N,\zetah}(i)\Big)^2}
  \end{split}     \end{equation}
  for all $M\in\N$.
  The right-hand side is finite due to Lemma~\ref{l:Y_2_supsum} for every $M\in\N$
  and converges to zero as $M\to\infty$.
  The Lipschitz condition~\eqref{eq:F_Lipschitz} implies that
  \begin{equation}  \begin{split} \label{eq:truncate.Fb}
    \varlimsup_{N\to\infty}\biggl|&\E\eckbb{M\wedge\Fb\robb{\Bigl(\sum_{i=1}^N
             \fb\rob{Y_{t,s}^{N,\zetah}(i)}\Bigr)_{t\in[s,\infty)}}}
           -\E\eckbb{M\wedge\Fb\robb{\Bigl(\sum_{i=1}^N \rob{\fb\mal h_K}
                    \rob{Y_{t,s}^{N,\zetah}(i)}\Bigr)_{t\in[s,\infty)}}}
          \biggr|\\
    &\leq\varlimsup_{N\to\infty}\sum_{j=1}^n L_\Fb L_\fb\E\eckbb{\sum_{i=1}^N 
           Y_{t_j,s}^{N,\zetah}(i)\roB{1-h_K\rob{Y_{t_j,s}^{N,\zetah}(i)}}}\\
    &\leq\varlimsup_{N\to\infty}\sum_{j=1}^n L_\Fb L_\fb\E\eckbb{\sum_{i=1}^N  Y_{t_j,s}^{N,\zetah}(i)
        \roB{\1_{\sup_{t\in[s,t_n]}Y_{t,s}^{N,\zetah}(i)\geq K}
            +\1_{Y_{t_j,s}^{N,\zetah}(i)\leq\frac{1}{K}}
            }                            }
    \\ &
    \leq\frac{n L_\Fb L_\fb}{K}\sup_{N\in\N}\E\eckbb{\sup_{t\in[s,t_n]}\Big(\sum_{i=1}^NY_{t,s}^{N,\zetah}(i)\Big)^2}
        +L_\Fb L_\fb\sum_{j=1}^n\varlimsup_{N\to\infty} N \E^0\eckB{ Y_{t_j,s}^{N,\zetah}(1)\wedge\tfrac{1}{K}}
    \\ &
    \leq\frac{n L_\Fb L_\fb}{K}\sup_{N\in\N}\E^0\eckbb{\sup_{t\in[s,t_n]}\Big(\sum_{i=1}^NY_{t,s}^{N,\zetah}(i)\Big)^2}
        +L_\Fb L_\fb\sum_{j=1}^n
        \int_s^{t_j} \zeta(r)\int \chi_{t_j-r}\wedge\tfrac{1}{K}\,\,Q(d\chi)\,dr
  \end{split}     \end{equation}
  for all $M,K\in\N$.
  For the last step, we applied Lemma~\ref{l:vanishing_immigration_general}.
  The right-hand side of~\eqref{eq:truncate.Fb} is finite according to
  Lemmas~\ref{l:Y_2_supsum} and~\ref{l:finite_excursion_area}
  and converges to zero as $K\to\infty$ for every $M\in\N$ according to the dominated convergence theorem.
  Thus
  \begin{equation}  \begin{split}
    &\limN\E\eckbb{\Fb\robb{\Bigl(\sum_{i=1}^N \fb\rob{Y_{t,s}^{N,\zetah}(i)}\Bigr)_{t\in[s,\infty)}}}
    =
    \lim_{M\to\infty}\limN\E\eckbb{M\wedge\Fb\robb{\Bigl(\sum_{i=1}^N \fb\rob{Y_{t,s}^{N,\zetah}(i)}\Bigr)_{t\in[s,\infty)}}}
    \\
    &=\lim_{M\to\infty}\limK\limN\E\eckbb{M\wedge\Fb\robb{\Bigl(\sum_{i=1}^N \rob{\fb\mal h_K}
                   \rob{Y_{t,s}^{N,\zetah}(i)}\Bigr)_{t\in[s,\infty)}}}\\
    &=\lim_{M\to\infty}\limK\E\eckbb{M\wedge\Fb\robb{\Bigl(
         \int \rob{\fb\mal h_K}\rob{\eta_{t-u}}\Pi(du,d\eta)\Bigr)_{t\in[s,\infty)}}}
    \\&
    =\E\eckbb{\Fb\robb{\Bigl(\int \fb\rob{\eta_{t-u}}\Pi(du,d\eta)\Bigr)_{t\in[s,\infty)}}}.
  \end{split}     \end{equation}
  The last equality follows from the dominated convergence theorem
  together with Assumption~\ref{a:Hutzenthaler2009EJP}.
  This proves~\eqref{eq:convY}.
\end{proof}
Remark: Instead of referring to Theorem 2.1 of Roelly-Coppoletta (1986)\nocite{RoellyCoppoletta1986}
one could alternatively apply Theorem 3.1 of Jakubowski (1986)\nocite{Jakubowski1986}.
\subsection{Convergence of the loop-free process}
\label{ssec:convergence_of_the_loop_free_process}
Recall the loop-free $N$-island process from~\eqref{eq:ZN}.
The following lemma shows that the loop-free $N$-island process
converges to the virgin island model.
\begin{lemma} \label{l:convergence_of_the_loop_free_process}
  Assume~\ref{a:A1_N}, \ref{a:Hutzenthaler2009EJP}, \ref{a:initial}
  and \ref{a:second_moments}.
  Then we have that
  \begin{equation}   \label{eq:convergence_of_the_loop_free_process}
    \biggl(\sum_{i=1}^N\sum_{k=0}^\infty
      \dl_{Z_{t}^{N,k}(i)}\biggr)_{t\leq T}\wlra
      \biggl(\sum_{(s,\eta)\in\VIM}\dl_{\eta_{t-s}}\biggr)_{t\leq T}\qqasN
  \end{equation}
  for every $T\in[0,\infty)$.
\end{lemma}
\begin{proof}
  Fix $T\in[0,\infty)$.
  Let $\MCF\subseteq\C_c^2\rob{(0,\infty),\R}$ be a countable dense subset of
  $\C_c\rob{(0,\infty),\R}$.
  Tightness of the left-hand side of~\eqref{eq:convergence_of_the_loop_free_process}
  in $N\in\N$
  follows from tightness of
  \begin{equation}  \label{eq:tightnesstoprove}
     \curlbb{\Bigl(\sum_{i=1}^{N}\sum_{k=0}^\infty
      f\rob{Z_{t}^{N,k}(i)}\Bigr)_{t\leq T}\colon N\in\N}
      \quad\fa f\in\MCF.
  \end{equation}
  This type of argument has been established in
  Theorem 2.1 of Roelly-Coppoletta \nocite{RoellyCoppoletta1986}(1986)
  for the weak topology and $\C_0$.
  Following the proof hereof, one can show the analogous argument 
  for the vague topology and $\C_c$.
  Fix $f\in\MCF$ and define
  $S_t^N:=\sum_{i=1}^N\sum_{k=0}^\infty f\rob{Z_{t}^{N,k}(i)}$ for all
  $t\leq T$ and $N\in\N$.
  For fixed $t\in[0,\infty)$, global Lipschitz continuity of $f$
  implies that
  \begin{equation}
    \P\eckB{\absb{S_t^N}\geq K}
    \leq \frac{L_f}{K}
         \E\eckbb{\sum_{i=1}^N\sum_{k=0}^\infty Z_{t}^{N,k}(i)}
    \leq \frac{L_f}{K}
         \sup_{M\in\N}\E\eckbb{\sum_{i=1}^M\sum_{k=0}^\infty Z_{t}^{M,k}(i)}
  \end{equation}
  for some constant $L_f<\infty$ and for all $K,N\in\N$.
  The right-hand side is finite according to Lemma~\ref{l:first_moment_estimate_X_Z}.
  This proves tightness of~\eqref{eq:tightnesstoprove} for every fixed time point.
  For the second part of the Aldous criterion, 
  let $\tau_N$, $N\in\N$, be stopping times which are uniformly bounded by $T$.
  In addition define
  $H(\uline{x}):=\sum_{i=1}^N\sum_{k=0}^\infty f(x_i^k)$ and
  \begin{equation}  \begin{split}
    \Gen^N H(\uline{x})&:=\sum_{i=1}^N\sum_{k=0}^\infty
      \roB{\frac{\1_{k>0}}{N}\sum_{j=1}^N x_j^{k-1}-x_i^{k}+\1_{k=0}\mu_N(0)+\mut_N(x_i^k)}
       f^{'}(x_i^k)
       +\sum_{i=1}^N\sum_{k=0}^\infty\frac{1}{2}\sigma_N^2(x_i^k)f^{''}(x_i^k)
  \end{split}     \end{equation}
  for all $\uline{x}=(x_i^k)_{i\leq N,k\in\N_0}\in I^{N\times\N_0}$.
  Assumption~\ref{a:A1_N} implies that the functions $\mu_N$, $\sigma_N$ and $\sigma_N^2$ are 
  uniformly globally Lipschitz continuous
  on the support of $f$.
  Moreover $N\mu_N(0)$ is bounded by $2\th$ uniformly in $N\in\N$.
  Therefore there exists a constant $C_H\in[1,\infty)$ such that
  $\abs{\Gen^N H(\uline{x})}\leq C_H\rob{2\th+\sum_{i=1}^N\sum_{k=0}^\infty x_i^k}$
  for all
  $\uline{x}=(x_i^k)_{i\leq N,k\in\N_0}$ and all $N\in\N$ and
  such that
  $\rob{(f^{'}\sigma_N)(y)}\leq C_H y$
  for all $y\in I$ and $N\in\N$.
  For fixed $\eta>0$ and $\dlb\in[0,1]$, we use It\^o's formula
  to obtain that
  \begin{equation}  \begin{split} \label{eq:esti_tightiZ}
    &\eta^2\P\eckB{\absb{S_{\tau_N+\dl}^N-S_{\tau_N}^N}>\eta}
    \leq\E\eckB{\Bigl(S_{\tau_N+\dl}^N-S_{\tau_N}^N\Bigr)^2}\\
    &=\E\eckbb{\Bigl(\int_{\tau_N}^{\tau_N+\dl}
            \Gen^N H\rob{Z_{u}^{N,\cdot}(\cdot)}\,du
            +\sum_{i=1}^N\sum_{k=0}^\infty\int_{\tau_N}^{\tau_N+\dl}
                    \rob{f^{'}\mal\sigma_N}\rob{Z_{u}^{N,k}(i)}dB_u^k(i)\Bigr)^2}\\
    &\leq    2C_H^2\E\eckbb{\Bigl(\int_0^{\dl} 2\th+\sum_{i=1}^N\sum_{k=0}^\infty
             Z_{\tau_N+u}^{N,k}(i)\,du\Bigr)^2}
    +2\sum_{i=1}^N\sum_{k=0}^\infty\E\eckB{\int_{0}^{\dl}
                    \rob{f^{'}\sigma_N}^2\rob{Z_{\tau_N+u}^{N,k}(i)}du}\\
    &\leq \dlb\mal 4C_H^2\sup_{M\in\N}\E\eckbb{\sup_{t\leq T+1}
          \Bigl(2\th+\sum_{i=1}^M\sum_{k=0}^\infty Z_{t}^{M,k}(i)\Bigr)^2}
  \end{split}     \end{equation}
  for all $\dl\leq \dlb$ and  for all $N\in\N$.
  The second step follows from $(a+b)^2\leq 2a^2+2b^2$ for $a,b\in\R$
  and from It{\^o}'s isometry.
  In the last step we used that
  $\rob{\int_0^\dl h(u)du}^2\leq\dl\int_0^\dl (h(u))^2du$
  for every integrable function $h$.
  The right-hand side of~\eqref{eq:esti_tightiZ}
  is finite according to Lemma~\ref{l:second_moment_estimate_X_Z}.
  Letting $\dlb\to0$, the left-hand side of~\eqref{eq:esti_tightiZ}
  converges to zero uniformly in $N\in\N$ and $\dl\leq\dlb$.
  Tightness of~\eqref{eq:tightnesstoprove} now follows from
  the Aldous criterion, see Aldous (1978)\nocite{Aldous1978}.

  It remains to identify the limit 
  of the left-hand side of~\eqref{eq:convergence_of_the_loop_free_process}
  by proving that
  \begin{equation}  \label{eq:loop_free_identify_limit}
    \limN\E\exp\ruB{-\sum_{j=1}^n\ld_j\sum_{k=0}^m\sum_{i=1}^N f\rob{Z_{t_j}^{N,k}(i)}
                   }
    =\E\exp\ruB{-\sum_{j=1}^n\ld_j\sum_{k=0}^m\sum_{(s,\eta)\in\VIM^{(k)}
                                  }
             f\rub{\eta_{t_j-s}}
              }
  \end{equation}
  for all $\ld_1,\ldots,\ld_n\geq0$, for all $0\leq t_1\leq\cdots\leq t_n$,
  for all $m,n\in\N_0$
  and for all non-negative $f\in\C_c^{2}\left((0,\infty)\right)$.
  Lemma~\ref{l:essentially_finitely_many_generations} justifies to restrict
  the summation over $k$ to finitely many summands.
  We prove~\eqref{eq:loop_free_identify_limit}
  by induction on $m\in\N_0$ using the Poisson
  limit~\eqref{eq:vanishing_immigration_weak_process}
  for independent one-dimensional diffusions.
  Define $F(\eta):=\sum_{j=1}^n \ld_j f(\eta_{t_j})$ for 
  every $\eta\in\C([0,\infty),I)$.
  Note that $F$ satisfies the Lipschitz condition~\eqref{eq:Lip_F}
  for some constant $L_F\in(0,\infty)$.
  Let $\rob{Y_{t,0}^{N,\zeta}(i)}_{t\geq0}$
  and $\rob{\Yt_{t,0}^{N,\zeta}(i)}_{t\geq0}$
  be solutions of~\eqref{eq:YN}
  with respect to the Brownian motion
  $(B_t^0(i))_{t\geq0}$
  such that
  $Y_{0,0}^{N,\zeta}(i)=X_0(i)$,
  $\Yt_{0,0}^{N,\zeta}(i)=0$ and $\zeta_N(\cdot)\equiv N\mu_N(0)$
  for every $N\in\N$ and $i\in\N$.
  In addition let $\rob{Y_t(i)}_{t\geq0}$, $i\in\N$, be independent solutions
  of~\eqref{eq:Y} with $Y_0(\cdot)\equiv X_0(\cdot)$.
  Note that $Z_{\cdot}^{N,0}(i)$ is the solution of~\eqref{eq:YN}
  with respect to the Brownian motion $(B_t^0(i))_{t\geq0}$
  with $Z_0^{N,0}(i)=X_0^N(i)$ for every $i\leq N$ and every $N\in\N$.
  Recall the random permutation $\pi^N$ from Assumption~\ref{a:initial}.
  The first-moment estimate of 
  Lemma~\ref{l:first_moment_estimate} implies that
  \begin{equation*}  \begin{split}    \label{eq:ZminusY}
    &\absbb{\E\eckbb{\exp\robb{-\sum_{i=1}^N F\roB{\rob{Z_t^{N,0}(i)}_{t\geq0}}}}
    - \E\eckbb{\exp\robb{-\sum_{i=1}^N F\roB{\rob{Y_{t,0}^{N,\zeta}(i)}_{t\geq0}}}}}\\
    &\leq L_F \sum_{j=1}^n \sum_{i=1}^N\E\eckbb{\absB{Z_{t_j,0}^{N,\zeta}(\pi_i^N)-Y_{t_j,0}^{N,\zeta}(i)}}
    \leq L_F e^{L_\mu t_n}\sum_{j=1}^n \E\eckbb{\sum_{i=1}^N\absB{X_0^N(\pi_i^N)-X_0(i)}}
  \end{split}     \end{equation*}
  for all $N\in\N$.
  of Lemma~\ref{l:first_moment_estimate}.
  Letting $N\to\infty$ the right-hand side converges to zero
  according to Assumption~\ref{a:initial}.
  The process 
  $\rob{Y_{\cdot,0}^{N,\zeta}(i)}_{i\leq N}$
  in turn
  is close to $\rob{\Yt_{\cdot,0}^{N,\zeta}(i)}_{i\leq N}$ except for islands with
  a significant amount of mass at time zero.
  Formalizing this we use Lemma~\ref{l:first_moment_estimate}
  to obtain that
  \begin{equation*}  \begin{split}    \label{eq:YminusYt}
    &\limK\limN\absbb{\E\eckbb{\exp\robb{-\sum_{i=K+1}^N F\roB{\rob{Y_{t,0}^{N,\zeta}(i)}_{t\geq0}}}}
    - \E\eckbb{\exp\robb{-\sum_{i=1}^N F\roB{\rob{\Yt_{t,0}^{N,\zeta}(i)}_{t\geq0}}}}}\\
    &\leq \limK L_F n e^{L_\mu t_n}\sum_{i=K+1}^\infty \E\eckb{X_0(i)}
      +\limK\limN\E\eckbb{1-\exp\robb{-\sum_{i=1}^K F\roB{\rob{\Yt_{t,0}^{N,\zeta}(i)}_{t\geq0}}}}.
  \end{split}     \end{equation*}
  The first summand on the right-hand side is zero according to Assumption~\ref{a:initial}.
  Note that $\mu_N\to\mu$ and $\sigma_N\to\sigma$
  as $N\to\infty$ by Assumption~\ref{a:A1_N}.
  Thus
  $\rob{\Yt_{t,0}^{N,\zeta}(i)}_{t\geq 0}$ converges in distribution to the zero
  function as $N\to\infty$ for every fixed $i\in\N$.
  Consequently the second summand on the right-hand side is zero.
  Moreover $\rob{Y_{t,0}^{N,\zeta}(i)}_{t\geq 0}$ converges in distribution to
  $\rob{Y_{t}(i)}_{t\geq 0}$
  as $N\to\infty$ for every fixed $i\in\N$.
  These observations imply that
  \begin{equation*}  \begin{split}
    &\limN \E\eckbb{\exp\robb{-\sum_{i=1}^N F\roB{\rob{Z_t^{N,0}(i)}_{t\geq0}}}}\\
    &=\limK\limN \E\eckbb{\exp\robb{-\sum_{i=1}^K F\roB{\rob{Y_{t,0}^{N,\zeta}(i)}_{t\geq0}}}}
         \E\eckbb{\exp\robb{-\sum_{i=K+1}^N F\roB{\rob{Y_{t,0}^{N,\zeta}(i)}_{t\geq0}}}}\\
    &=\limK\E\eckbb{\exp\robb{-\sum_{i=1}^K F\roB{\rob{Y_t(i)}_{t\geq0}}}}
      \limN\E\eckbb{\exp\robb{-\sum_{i=1}^N F\roB{\rob{\Yt_{t,0}^{N,\zeta}(i)}_{t\geq0}}}}\\
    &=\E\eckbb{\exp\robb{-\sum_{i=1}^\infty F\roB{\rob{Y_t(i)}_{t\geq0}}}}
       \E\eckbb{\exp\robb{-\!\!\sum_{(s,\eta)\in\Pi^{\th}}
                F\rub{(\eta_{t-s})_{t\geq0}}
                        }
               }\\
    &=
    \E\eckbb{\exp\ruB{-\!\!\sum_{(s,\eta)\in\VIM^{(0)}
                      }
             F\rub{(\eta_{t-s})_{t\geq0}}
               }
            }.
  \end{split}     \end{equation*}
  The last but one step follows from
  Lemma~\ref{l:vanishing_immigration_weak_process}
  with $\zetah_N(\cdot)\equiv N\mu_N(0)$
  and $\zeta(\cdot)\equiv\th$.
  This proves~\eqref{eq:loop_free_identify_limit} in the base case $m=0$.

  For the induction step $m\to m+1$,
  we observe that a version of $(Z_t^{N,m+1}(i))_{t\geq0}$ conditioned on
  $\zeta_N(r):={\sum_{j=1}^N Z_r^{N,m}(j)}$, $r\geq0$, is given by
  the one-dimensional diffusion
  $(Y_{t,0}^{N,\zeta})_{t\geq0}$ with vanishing immigration.
  Thus we may realize $(Z_t^{N,m+1}(i))_{t\geq0}$ by choosing 
  a suitable version of $(Z_t^{N,m}(j))_{t\geq0}$, $j=1,\dots,N$,
  and by independently sampling
  a version of $(Y_{t,0}^{N,\zeta})_{t\geq0}$ whose driving Brownian
  motion is independent of
  $\{(Z_t^{N,m}(j))_{t\geq0}\colon j=1,\ldots,N\}$.
  Tightness of $\big\{\big(\sum_{i=1}^N\sum_{k=0}^{\tilde{m}}\delta_{Z_t^{N,k}(i)}\big)_{t\leq T}\big\}$
  together with the induction hypothesis implies that
  \begin{equation}   \label{eq:convergence_of_the_loop_free_process_finite_summation}
    \biggl(\sum_{i=1}^N\sum_{k=0}^\mt
      \dl_{Z_{t}^{N,k}(i)}\biggr)_{t\leq T}\wlra
      \biggl(\sum_{k=0}^\mt\sum_{(s,\eta)\in\VIM^{(k)}}
                       \dl_{\eta_{t-s}}\biggr)_{t\leq T}\qqasN
  \end{equation}
  for every $\mt\leq m$.
  Thanks to the Skorokhod
  representation of weak convergence
  (e.g.\ Theorem II.86.1 in~\cite{RogersWilliams2000a}),
  we may assume that the convergence
  in~\eqref{eq:convergence_of_the_loop_free_process_finite_summation}
  holds almost surely.
  As a consequence we obtain that
  \begin{equation} 
    \biggl(\sum_{i=1}^N
      \dl_{Z_{t}^{N,m}(i)}\biggr)_{t\leq T}\lra
      \biggl(\sum_{(s,\eta)\in\VIM^{(m)}}
                       \dl_{\eta_{t-s}}\biggr)_{t\leq T}\qqasN
  \end{equation}
  holds almost surely.
  Using arguments from the proof of~\eqref{eq:convY}, one can 
  deduce from this that
  \begin{equation}  \label{eq:lasteq}
    \biggl(\sum_{i=1}^N
      Z_{t}^{N,m}(i)\biggr)_{t\leq T}\wlra
      \biggl(V_t^{(m)}\biggr)_{t\leq T}\qqasN
  \end{equation}
  holds almost surely 
  where
  the total mass of the $n$-th generation of the virgin island model is defined as
\begin{equation}  \label{eq:def:Vn}
  V_t^{(n)}:=\sum_{\ru{s,\eta}\in\VIM^{(n)}}\eta_{t-s},\quad t\geq0,
\end  {equation}
  for every $n\in\N_0$.
  Together with continuity
  of $\rob{\eta_{s}}_{s\leq T}\mapsto \int_0^t\abs{\eta_s-V_s^{(m)}}\,ds$, \eqref{eq:lasteq} implies that
  \begin{equation}
    \int_0^t\absb{\sum_{i=1}^N Z_s^{N,m}(i)-V_s^{(m)}}\,ds
    \lraN 0\qquad\text{almost surely}.
  \end{equation}
  Now the main step of the proof is 
  Lemma~\ref{l:vanishing_immigration_weak_process}
  with $\zetah_N(t):=\sum_{i=1}^N Z_t^{N,m}(i)$
  and $\zeta(t):=V_t^{(m)}$ for all $t\geq0$.
  Lemma~\ref{l:vanishing_immigration_weak_process}
  implies that
  \begin{equation}  \begin{split}
    \E&\eckB{\exp\ruB{-\sum_{i=1}^N F\roB{\rub{Z_t^{N,m+1}}_{t\leq T}}}|
         \rub{Z_\cdot^{M,m}}_{M\in\N}}\\
    &=\E\eckB{\exp\ruB{-\sum_{i=1}^N F\roB{\rub{Y_{t,0}^{N,\zeta
                                               }
                                               }_{t\leq T}
                                          }
                      }
             \Big|\rub{Z_\cdot^{M,m}}_{M\in\N}
            }\\
    &\lra
    \E\eckB{\exp\ruB{-\int F\rub{(\eta_{t-s})_{t\leq T}}\Pi^{(m+1)}(ds,d\eta)}\Big|
         \VIM^{(m)}}\qqasN
  \end{split}     \end{equation}
  almost surely.
  Here $\Pi^{(m+1)}$ conditioned on $\VIM^{(m)}$ is a Poisson point process
  on $[0,\infty)\times U$ with intensity measure
  \begin{equation}  \begin{split}
    \E&\eckbb{  \Pi^{(m+1)}(ds,d\eta)|\VIM^{(m)}
            } 
    = V_s^{(m)}\,ds\otimes Q(d\eta)\\
    &=\!\!\!\sum_{(r,\psi)\in\VIM^{(m)}}\!\!\!\psi_{s-r}ds\otimes Q(d\eta)
    =\!\!\!\sum_{(r,\psi)\in\VIM^{(m)}}\!\!\!\E\Pi^{(r,\psi)}(ds,d\eta).
  \end{split}     \end{equation}
  Due to this decomposition of $\Pi^{(m+1)}$, we may realize $\Pi^{(m+1)}$
  conditioned on $\VIM^{(m)}$ as the independent superposition of
  $\{\Pi^{(m,s,\psi)}\colon(s,\psi)\in\VIM^{(m)}\}$.
  In other words, $\Pi^{(m+1)}$ is equal in distribution to the $(m+1)$st-generation
  of the virgin island model.
  Therefore we get that
  \begin{equation*}  \begin{split}
    \lefteqn{\limN\E\eckbb{\exp\ruB{-\sum_{k=0}^{m+1}\sum_{i=1}^N
              F\roB{\rub{Z_t^{N,k}(i)}_{t\leq T}}
            }      }      }\\
    &= \E\eckbb{\limN\exp\ruB{-\sum_{k=0}^{m}\sum_{i=1}^N F\rub{Z_{\cdot}^{N,k}(i)}
                             }
                  \E\eckB{\exp\ruB{-\sum_{i=1}^N F\rub{Z_{\cdot}^{N,m+1}(i)}
                                  }\Big|\rub{Z_\cdot^{M,m}}_{M\in\N}
                        }        }\\
    &=\E\biggl[\exp\rubb{-\sum_{k=0}^m\sum_{(s,\eta)\in\VIM^{(k)}
                                  }\!\!\!
             F\rub{(\eta_{t-s})_{t\leq T}}
               }
      \E\eckB{\exp\ruB{-\!\!\!\sum_{(s,\eta)\in\VIM^{(m+1)}}\!\!\!
             F\rub{(\eta_{t-s})_{t\leq T}}}\Big| \VIM^{(m)}
              }\biggr]\\
    &=\E\eckbb{\exp\ruB{-\sum_{k=0}^{m+1}\sum_{(s,\eta)\in\VIM^{(k)}
                                  }
             F\rub{(\eta_{t-s})_{t\leq T}}
              }        }
  \end{split}     \end{equation*}
  which proves~\eqref{eq:loop_free_identify_limit}
  and completes the proof
  of Lemma~\ref{l:convergence_of_the_loop_free_process}.
\end{proof}

\subsection{Reducing to the loop-free $(N,\mu_N,\sigma_N^2)$-process}
\label{ssec:reducing}
Next we show that the $(N,\mu_N,\sigma_N^2)$-process
with migration levels
and the loop-free $(N,\mu_N,\sigma_N^2)$-process are identical in the limit $N\to\infty$.
Our proof formalizes the following intuition. The individuals of a certain 
migration level
are concentrated on essentially finitely many islands. That these
finitely many islands are populated by migrants of a different
migration level has a probability of order $\tfrac{1}{N}$.
As a consequence, all individuals on one fixed
island have the same migration level in the limit $N\to\infty$.
This intuition is subject of Lemma~\ref{l:one_generation_per_island}.

First we show that a generation cannot be dispersed uniformly over all islands.
To obtain this interpretation from the following lemma,
assume $X_t^{N,k}(i)\approx\tfrac{1}{N}$ for all $i\leq N$ and some time $t\geq0$.
Then the cutting operation in~\eqref{eq:l:concentration}
has no effect for $N$ large enough.
However it is clear that
the total mass of all individuals with migration level $k$ does
not tend to zero as $N\to\infty$.
Thus, as a consequence of the following lemma, $X_t^{N,k}(i)\approx\tfrac{1}{N}$ cannot be true.
%
\begin{lemma}       \label{l:concentration}
  Assume~\ref{a:A1_N}, \ref{a:initial}
  and \ref{a:second_moments}.
  Then any solution of~\eqref{eq:XN_k} satisfies
  \begin{equation} \label{eq:l:concentration}
    \sum_{k\geq0}\sup_{N\in\N}\E\eckbb{\sum_{i=1}^N
       \ruB{X_{t}^{N,k}(i)\wedge\dl}}
    \lradl 0
  \end{equation}
  for all $t\in[0,\infty)$.
  The assertion is also true if $X_t^{N,k}(i)$ is
  replaced by $Z_t^{N,k}(i)$.
\end{lemma}
\begin{proof}
  Fix $T\in[0,\infty)$ and $t\in[0,T]$.
  According to Lemma~\ref{l:essentially_finitely_many_generations},
  for each $\eps>0$,
  there exists a $k_0\in\N$
  such that
  \begin{equation}
    \sum_{k\geq k_0}\sup_{N\in\N}\sup_{s\leq T}
      \E\eckbb{\sum_{i=1}^N X_s^{N,k}(i)}\leq\eps.
  \end{equation}
  Thus it suffices to prove convergence of every summand
  in~\eqref{eq:l:concentration}.
  In addition,
  if we forget the migration levels
  in the $(N,\mu_N,\sigma_N^2)$-process with migration levels, then we obtain the
  $N$-island model.
  More formally, Lemma~\ref{l:decomposition} shows that
  \begin{equation}
    \Xb_t^N(i):=\sum_{m\geq0}X_t^{N,m}(i)\quad\fa i\leq N,t\geq0
  \end{equation}
  defines an $N$-island model.
  Recall $\taut_{K}^N$ from~\eqref{eq:tautKN} and fix $K\in\N$.
  According to Lemma~\ref{l:tau_theta} it suffices to prove~\eqref{eq:l:concentration}
  with expectation being restricted to the event
  $\{\taut_K^N> t\}$ for every $N\in\N$.
  Now Lemma 3.3 in~\cite{HutzenthalerWakolbinger2007} shows that, on the event $\{\taut_K^N>t\}$,
  $\Xb_t^N(i)$
  is stochastically bounded above by
  $Y_{t,0}^{N,K+N\mu_N(0)}$.
  Hence
  we get for all $N,K\in\N$, $k\in\N_0$ and $\dlt>0$ that
  \begin{equation}  \begin{split}  \label{eq:tau_not_yet}
    \lefteqn{\E\eckbb{\sum_{i=1}^N\ruB{X_t^{N,k}(i)\wedge\dl}\1_{\taut_K^N> t}}}\\
    &\leq 
      \sum_{i=1}^N\E\eckbb{\1_{\Xb_0^N(i)\leq\dlt}
        \ruB{X_t^{N,k}(i)\wedge\dl}\1_{\taut_K^N> t}}
      +\dl\E\eckbb{\sum_{i=1}^N\1_{\Xb_0^N(i)>\dlt}}\\
    &\leq 
      \sum_{i=1}^N\E^{\Xb_0^N(i)\wedge\dlt} \eckB{Y_{t,0}^{N,K+N\mu_N(0)}\wedge\dl}
      +\frac{\dl}{\dlt}\E\eckB{\sum_{i=1}^N\Xb_0^N(i)}\\
    &\leq 
      N\E^{0} \eckB{Y_{t,0}^{N,K+N\mu_N(0)}\wedge\dl}
      +C_t\sum_{i=1}^N\E\eckB{\ruB{\Xb_0^N(i)\wedge\dlt}}
      +\frac{\dl}{\dlt}\sup_{M\in\N}\E\eckB{\sum_{i=1}^M\Xb_0^M(i)}
  \end{split}     \end{equation}
  for some constant $C_t<\infty$.
  The last step follows from
  Lemma~\ref{l:first_moment_estimate}.
  Next we let $N\to\infty$ and $\dl\to0$
  in~\eqref{eq:tau_not_yet}.
  Applying Lemma~\ref{l:vanishing_immigration_weak_process}
  with $\zetah_N:=K+N\mu_N(0)$
  and using Assumption~\ref{a:initial},
  we see that the limit of~\eqref{eq:tau_not_yet}
  as $N\to\infty$ and as $\dl\to0$ is bounded above by
  \begin{equation}    \label{eq:tau_not_yet_upper_bound}
      \limdl(K+\th)\int\int_0^t\rub{\chi_{t-r}\wedge\dl}\,dr\, Q(d\chi)
      +C_t\E\eckB{\sum_{i=1}^\infty\ruB{\Xb_0(i)\wedge\dlt}}
  \end{equation}
  for every $K\in\N$ and $\dlt>0$.
  The first summand in~\eqref{eq:tau_not_yet_upper_bound} is zero by the
  dominated convergence theorem and Lemma~\ref{l:finite_excursion_area}.
  The last summand converges to zero as $\dlt\to0$
  by the dominated convergence theorem and Assumption~\ref{a:initial}.
  This completes the proof of~\eqref{eq:l:concentration}.
  The proof of~\eqref{eq:l:concentration} 
  with $X_t^{N,k}(i)$ replaced by $Z_t^{N,k}(i)$ is similar
\end{proof}

%
Now we prove that all individuals on a fixed island have the same
migration level in the limit $N\to\infty$.
More precisely, we show that $X_t^{N,k}(i)$
and $\sum_{m\neq k}X_t^{N,m}(i)$ cannot be big at the same time for any
$i\leq N$ or any $k\in\N_0$.
\begin{lemma}       \label{l:one_generation_per_island}
  Assume~\ref{a:A1_N}, \ref{a:initial} and~\ref{a:second_moments}.
  Then any solution of~\eqref{eq:XN_k} satisfies
  \begin{equation} \label{eq:l:one_generation_per_island}
    \limN\sup_{s,t\leq T}\E\eckbb{\sum_{i=1}^N\sum_{k\geq0}
        X_s^{N,k}(i)\sum_{m\neq k}X_t^{N,m}(i)}
    =0
  \end{equation}
  for every $T\in[0,\infty)$. The assertion is also true if $X_t^{N,k}(i)$ is
  replaced by $Z_t^{N,k}(i)$.
\end{lemma}
\begin{proof}
  We begin with the case $s=t$.
  Fix $K\in\N$ and
  recall $\taut_K^N$ from~\eqref{eq:tautKN}.
  By the uniform local Lipschitz continuity of $\mu_N$,
  there exists a finite constant $L_K$
  such that $\mut_N(x)\leq L_K x$ for all $x\leq K$ and all $N\in\N$.
  Applying It{\^o}'s formula, we see that
    \begin{equation}  \begin{split}
      d&\sum_{i=1}^N\sum_{k\geq0} X_{t\wedge\taut_K^N}^{N,k}(i)
          \sum_{m\neq k}X_{t\wedge\taut_K^N}^{N,m}(i)
      \leq 2\sum_{i=1}^N\sum_{k\geq 0} X_{t\wedge\taut_K^N}^{N,k}(i)\mu_N(0)\,dt\\
      &\quad +\sum_{i=1}^N\sum_{k\geq0} X_{t\wedge\taut_K^N}^{N,k}(i)
          \sum_{m\neq k}
            \rubb{\frac{1}{N}\sum_{j=1}^N X_{t\wedge\taut_K^N}^{N,m-1}(j)\,dt
                +L_K X_{t\wedge\taut_K^N}^{N,m}(i)\,dt+d\,M_t^{N,m}(i)}\\
      &\quad +\sum_{i=1}^N\sum_{k\geq0} 
          \sum_{m\neq k}X_{t\wedge\taut_K^N}^{N,m}(i)
            \rubb{\frac{1}{N}\sum_{j=1}^N X_{t\wedge\taut_K^N}^{N,k-1}(j)\,dt
                +L_K X_{t\wedge\taut_K^N}^{N,k}(i)\,dt+d\,\Mt_t^{N,k}(i)}
    \end{split}     \end{equation}
  where $(M_t^{N,m}(i))_{t\geq0}$
  and
  $(\Mt_t^{N,k}(i))_{t\geq0}$ are suitable martingales
  for each $i\leq N$ and $k,m\in\N_0$.
  Now take expectations to obtain that
  \begin{equation}  \begin{split}
    \lefteqn{\E\eckbb{\sum_{i=1}^N \sum_{k\geq0}X_{t\wedge\taut_K^N}^{N,k}(i)
          \sum_{m\neq k}X_{t\wedge\taut_K^N}^{N,m}(i)}
          \leq2\mu_N(0)\int_0^t
               \E\eckB{\sum_{i=1}^N\sum_{k\geq0} X_{s\wedge\taut_K^N}^{N,k}(i)
                     }
                       \,ds
           }\\
    &\qquad\qquad\qquad\quad+\int_0^t\E\eckbb{\sum_{i=1}^N\sum_{k\geq0} X_{s\wedge\taut_K^N}^{N,k}(i)
          \frac{1}{N}
              \sum_{j=1}^N\sum_{m\geq 0} X_{s\wedge\taut_K^N}^{N,m-1}(j)}\,ds\\
    &\qquad\qquad\qquad\quad+\int_0^t\E\eckbb{
         \sum_{k\geq0} 
         \frac{1}{N}\sum_{j=1}^N X_{s\wedge\taut_K^N}^{N,k-1}(j)\sum_{i=1}^N\sum_{m\geq 0}
          X_{s\wedge\taut_K^N}^{N,m}(i)}\,ds\\
    &\qquad\qquad\qquad\quad+2L_K\int_0^t\E\eckbb{\sum_{i=1}^N\sum_{k\geq0}
          X_{s\wedge\taut_K^N}^{N,k}(i)
          \sum_{m\neq k} X_{s\wedge\taut_K^N}^{N,m}(i)}\,ds
  \end{split}     \end{equation}
  for every $N\in\N$ and $t\leq T$.
  Note that the right-hand side is finite.
  Using Gronwall's inequality, $\mu_N(0)\leq 2\th/N$ and
  $\sum_{i=1}^N \sum_{k\geq0}X_{t\wedge\taut_K^N}^{N,k}(i)\leq K$, we conclude that
  \begin{equation}  \begin{split}
    \sup_{t\leq T}\E\eckbb{\sum_{i=1}^N\sum_{k\geq0} X_{t\wedge\taut_K^N}^{N,k}(i)
          \sum_{m\neq k}X_{t\wedge\taut_K^N}^{N,m}(i)}
          \leq \frac{1}{N}T(4\th K+2K^2)\mal e^{2L_K T}
  \end{split}     \end{equation}
  for every $N\in\N$.
  Letting $N\to\infty$ proves \eqref{eq:l:one_generation_per_island}
  in the case $s=t$.
  For the case $s<t$, apply It{\^o}'s formula in
  \begin{equation}
    X_{s\wedge\taut_K^N}^{N,k}(i)\sum_{m\neq k}
      \int_{s\wedge\taut_K^N}^{t\wedge\taut_K^N}d X_r^{N,m}(i)
  \end{equation}
  and use similar estimates as above. The case $s>t$ is analogous to the case $s<t$.
%
\end{proof}

Knowing that there is asymptotically at most one generation on every island,
we are now in a position to prove that
$\rub{X_t^{N,k}(i)}_{t\geq0}$
and $\rub{Z_t^{N,k}(i)}_{t\geq0}$ are close to each other.
%
\begin{lemma}   \label{l:X_close_to_Z}
  Assume~\ref{a:A1_N_linear}, \ref{a:initial}
  and \ref{a:second_moments}.
  For each $N\in\N$, let 
  \begin{equation}
    \curlb{\rub{X_t^{N,k}(i),\Bb_t^k(i)}_{t\geq0}\colon k\in\N_0,i\leq N}
  \end{equation}
  be a solution of~\eqref{eq:XN_k} and let
  \begin{equation}
    \curlb{\ru{Z_t^{N,k}(i)}_{t\geq0}\colon k\in\N_0,i\leq N}
  \end{equation}
  be the unique solution of~\eqref{eq:ZN} with
  $\curlb{(  B_t^k(i))_{t\geq0}}$
  replaced by 
  $\curlb{(\Bb_t^k(i))_{t\geq0}}$.
  Then
  \begin{equation}  \label{eq:X_close_to_Z}
    \E\eckbb{\sum_{i=1}^N\sum_{k=0}^\infty
       \absb{X_t^{N,k}(i)-Z_t^{N,k}(i)}}\lra0
       \qqasN
  \end{equation}
  for every $t\in[0,\infty)$.
\end{lemma}
\begin{proof}
  In the first step, assume that $\sigma_N^2$ and $\mu_N$ are uniformly
  globally Lipschitz continuous
  and bounded.
  The general case will later be handled by a stopping argument.

  Define
  $x^+:=\max(x,0)$ and $x^-:=(-x)^+$ for all $x\in\R$.
  We first prove~\eqref{eq:X_close_to_Z} with $\abs{x}$ replaced by
  $x^+$ and $x^-$, respectively, separately.
  As $\R\ni x\mapsto x^+$ is not differentiable, we will apply It{\^o}'s formula
  to $\phi_n$ being defined as follows.
  Let $1=a_0>a_1>\cdots>a_n>\cdots>0$ satisfy
  \begin{equation}
    \int_{a_1}^1 \frac{1}{u}\,du=1,\;\int_{a_2}^{a_1}\frac{1}{u}\,du=2,
    ...,\int_{a_n}^{a_{n-1}}\frac{1}{u}\,du=n,...
  \end  {equation}
  Note that $a_n\to0$ as $n\to\infty$.
  For every $n\in\N$, there exists a continuous function
  $\psi_n\colon\R\to[0,\infty)$ with support in $(a_n,a_{n-1})$ such that
  \begin{equation}
    0\leq \psi_n(u)\leq\frac{2}{nu}\text{ for all $u>0$ and }
    \int_{a_n}^{a_{n-1}}\psi_n(x)\,dx=1.
  \end  {equation}
  Note that $\psi_n(u)\leq\tfrac{2}{na_n}$ for all $u\geq0$ and $n\in\N$.
  With this, we define the function
  \begin{equation}
    \phi_n(x):=\1_{x>0}\int_0^{{x}}dy\int_0^y\psi_n(u)\,du,\quad x\in\R,
  \end  {equation}
  for every $n\in\N$.
  These functions
  satisfy $\phi_n\in \C^2(\R)$, $\abs{\phi_n^{'}(x)}\leq1$,
  $\phi_n^{''}(x)=\1_{x>0}\psi_n({x})$, $\phi_n(x)\leq{x}^+$ and $\phi_k(x)
  {\ra}{x}^+$ as $k\to\infty$ for every $x\in\R$ and $n\in\N$.

  Denote the difference process by
  $\Delta_t^{k}(i):=X_{t}^{N,k}(i)-Z_{t}^{N,k}(i)$ for all $i\in G$, $t\geq0$,
  $k\in\N_0$ and $N\in\N$.
  The dependence on $N$ is suppressed for the sake of a more compact
  notation.
  By definition of $\phi_n$,
  $x^{+}\leq\phi_{n }(x)+a_{n -1}\wedge x^{+}$ for all $x\in\R$
  and $n\in\N$.
  Thus
  \begin{equation}  \begin{split}     \label{eq:conv78}
    \E&\eckbb{\sum_{i=1}^N\sum_{k=0}^{\infty}
        \rub{\Delta_t^{k}(i)}^{+}}
    \leq
    \sum_{k=K+1}^\infty \sup_{M\in\N}\sum_{i=1}^M\E\eckbb{X_t^{M,k}(i)}\\
    &\qquad+\E\eckbb{\sum_{i=1}^N \sum_{k=0}^{K}\phi_{N^2}
        \rub{\Delta_t^{k}(i)}}
    + \sum_{k=0}^{\infty}\sup_{M\geq 1}\E\eckbb{\sum_{i=1}^{M}
         \ruB{a_{N^2-1}\wedge X_{t}^{M,k}(i)}}
  \end{split}     \end{equation}
  for all $K\in\N$ and $N\in\N$.
  The first summand on the right-hand side converges to zero
  as $K\to\infty$ uniformly in $N\in\N$ according to
  Lemma~\ref{l:essentially_finitely_many_generations}.
  The last summand on the right-hand side converges to zero
  as $N\to\infty$ according to Lemma~\ref{l:concentration}.
  Consequently 
  \begin{equation}  \label{eq:an_estimate_as}
    \limsupN\,\E\eckbb{\sum_{i=1}^{N}\sum_{k=0}^{\infty}\ruB{\Delta_{t}^k(i)}^{+}}
    \leq\limsupK\,\limsupN\,\E\eckbb{\sum_{i=1}^{N}
         \sum_{k=0}^{K}\phi_{N^2}\ruB{\Delta_{t}^k(i)}}.
  \end{equation}
  To estimate the right-hand side, we apply
  It{\^o}'s formula to obtain that
  \begin{equation}  \begin{split}    \label{eq:phi_Delta}
    \lefteqn{d\sum_{i=1}^N \phi_{N^2}\rub{\Delta_t^k(i)}
        =
        }\\
    &\ \phantom{+} \sum_{i=1}^N\phi_{N^2}^{'}
        \left(\sqrt{
           \frac{X_t^{N,k}(i)
           \mal \sigma_N^2\ruB{\sum_{m=0}^\infty X_t^{N,m}(i)}}
                   {\sum_{m=0}^\infty X_t^{N,m}(i)}
                   }
             -\sqrt{
                   \sigma_N^2\rub{Z_t^{N,k}(i)}
                   }
        \right)\,d\Bb_t^k(i)
        \\
    &\ +\sum_{i=1}^N\phi_{N^2}^{'}\rub{\Delta_t^k(i)}
       \rubb{
            \frac{1}{N}\sum_{j=1}^N \Delta_t^{k-1}(j)-\Delta_t^{k}(i)
            }\,dt\\
    &\ +\sum_{i=1}^N\phi_{N^2}^{'}\rub{\Delta_t^k(i)}
       \rubb{
           \frac{X_t^{N,k}(i)
           \mal \mut_N\ruB{\sum_{m=0}^\infty X_t^{N,m}(i)}}{\sum_{m=0}^\infty X_t^{N,m}(i)}
           -\mut_N\ruB{Z_t^{N,k}(i)}
            }\,dt\\
    &\ 
      +\sum_{i=1}^N\frac{\phi_{N^2}^{''}\rub{\Delta_t^k(i)}}{2}
        \left(\sqrt{
           \frac{X_t^{N,k}(i)
           \mal \sigma_N^2\ruB{\sum_{m=0}^\infty X_t^{N,m}(i)}}{\sum_{m=0}^\infty X_t^{N,m}(i)}
                   }
             -\sqrt{
                   \sigma_N^2\rub{Z_t^{N,k}(i)}
                   }
        \right)^2\,dt
  \end{split}     \end{equation}
  for every $k\in\N_0$.
  Now we simplify the last summand on the right-hand side using
  the assumption of $\sigma_N^2$ being linear.
  This linearity implies that
  $x\cdot\sigma_N^2(y)/y=\sigma_N^2(x)$
  for $x=X_t^{N,k}(i)$ and $y=\sum_{m=0}^\infty X_t^{N,m}(i)$.
  In addition
  take expectations, estimate $\abs{\phi_{N^2}^{'}(x)}\leq 1$
  and use $-\phi_{N^2}^{'}(x)\cdot x\leq0$ for $x\in\R$.
  Thus we have that
  \begin{equation}  \begin{split}    \label{eq:E_phi_Delta}
    \lefteqn{\ddt \E\eckB{\sum_{i=1}^N \phi_{N^2}\rub{\Delta_t^k(i)}}
            }\\
    &\leq \E\eckB{\sum_{j=1}^N\absb{ \Delta_t^{k-1}(j) }}
         +\E\eckB{\sum_{i=1}^N\absb{\mut_N\rob{X_t^{N,k}(i)}
                                 -\mut_N\rob{Z_t^{N,k}(i)}
                                   }
               }\\
    &\ +\E\eckbb{\sum_{i=1}^N\absbb{
           \frac{X_t^{N,k}(i)
                   \mal \mut_N\ruB{\sum_{m=0}^\infty X_t^{N,m}(i)}}
                {\sum_{m=0}^\infty X_t^{N,m}(i)}
           -\mut_N\ruB{X_t^{N,k}(i)}
                                  }
               }\\
    &\ +\E\eckbb{\sum_{i=1}^N
           \frac{\phi_{N^2}^{''}\rob{\Delta_t^{k}(i)}}{2}
        \left(\sqrt{
                   \sigma_N^2\rub{X_t^{N,k}(i)}
                   }
             -\sqrt{
                   \sigma_N^2\rub{Z_t^{N,k}(i)}
                   }
        \right)^2
               }
  \end{split}     \end{equation}
  for every $t\geq0$, $k\in\N_0$ and $N\in\N$.
  The last but one summand on the right-hand side
  is estimated as follows.
  Let $\dl>0$.
  In case of $X_t^{N,k}(i)\geq\dl$, apply the inequality
  $\abs{\tfrac{a}{a+b}\mut_N(a+b)-\mut_N(a)}\leq 2L_{\mu} b$ with
  $a=X_t^{N,k}(i)$ and $b=\sum_{m\neq k}^\infty X_t^{N,m}(i)$.
  In case of $X_t^{N,k}<\dl$, use $\abs{\mut_N(x)}\leq L_{\mu} x$ for all $x\geq0$.
  For the last summand on the right-hand side, use
  $(\sqrt{x}-\sqrt{y})^2\leq\abs{x-y}$ and
  $\abs{\sigma_N^2(x)-\sigma_N^2(y)}\leq L_{\sigma}\abs{x-y}$ for
  $x,y\in[0,\infty)$.
  Moreover recall $\phi_{N^2}^{''}(x)\leq\tfrac{2}{N^2x^{+}}$
  for all $x>0$.
  Define $L:=L_\mu\vee L_\sigma$.
  Hence, we get that
  \begin{equation}  \begin{split}  \label{eq:Delta_after_expectations}
    \lefteqn{\ddt \E\eckB{\sum_{i=1}^N \phi_{N^2}\rub{\Delta_t^k(i)}}
            }\\
    &\leq \E\eckB{\sum_{j=1}^N \absb{\Delta_t^{k-1}(j)}}
        +L\E\eckB{\sum_{i=1}^N\absb{
                                     X_t^{N,k}(i)-Z_t^{N,k}(i)
                                   }
                }\\
    &\ +2L\E\eckbb{\sum_{i=1}^N \1_{X_t^{N,k}(i)\geq\dl}
         \sum_{m\neq k}X_t^{N,m}(i)
                  }\\
    & +2 L\E\eckbb{ \sum_{i=1}^{N}\ruB{X_t^{N,k}(i)\wedge\dl}
                  }
      +\E\eckbb{\sum_{i=1}^N
         \frac{L\rub{\Delta_t^k(i)}^{+}}{N^2\rub{\Delta_t^k(i)}^{+}}
               }
  \end{split}     \end{equation}
  for every $t\geq0$, $k\in\N_0$ and $N\in\N$.
  Summing over $k\in\{0,\ldots,K\}$ leads to
  \begin{equation}  \begin{split}  \label{eq:summed_over_k}
    \lefteqn{\ddt \E\eckbb{\sum_{k=0}^K\sum_{i=1}^N
              \phi_{N^2}\rub{\Delta_t^k(i)}}
            }\\
    &\leq (1+L)\E\eckbb{\sum_{k=0}^\infty\sum_{i=1}^N\absB{\Delta_t^{k}(i)}}
       +\frac{2L}{\dl}\E\eckbb{\sum_{k=0}^\infty\sum_{i=1}^N X_t^{N,k}(i)
         \sum_{m\neq k}X_t^{N,m}(i)
                  }\\
    & +2 L\E\eckbb{\sum_{k=0}^\infty \sum_{i=1}^{N}\ruB{X_t^{N,k}(i)\wedge\dl}
                  }
      +(K+1)\frac{L}{N}
  \end{split}     \end{equation}
  for every $t\geq0$ and $K,N\in\N$.
  The last three summands on the right-hand
  side of~\eqref{eq:summed_over_k}  converge to zero
  uniformly in $t\in[0,T]$
  if we first let $N\to\infty$ and then $\dl\to0$.
  This follows from
  Lemma~\ref{l:one_generation_per_island} and
  Lemma~\ref{l:concentration}.
  After inserting~\eqref{eq:summed_over_k}
  into~\eqref{eq:an_estimate_as},
  we see that
  \begin{equation}  \label{eq:X_Z_difference}
    \limsupN\,\E\eckbb{\sum_{i=1}^N\sum_{k=0}^{\infty}
        \ruB{\Delta_t^{k}(i)}^{+}}
    \leq\int_0^t(1+L)\,\limsupN\,\E\eckbb{\sum_{i=1}^N\sum_{k=0}^{\infty}
        \absB{\Delta_s^{k}(i)}}\,ds
  \end{equation}
  for all $t\in[0,T]$.
  Note that the right-hand side is finite by
  Lemma~\ref{l:essentially_finitely_many_generations}.
  Similarly we obtain~\eqref{eq:X_Z_difference} with
  $(\cdot)^+$ replaced by $\abs{\cdot}$.
  Finally apply Gronwall's inequality to arrive at~\eqref{eq:X_close_to_Z}.
  
  In the second step, we consider functions $\sigma_N^2$ and $\mu_N$ which are not globally
  Lip\-schitz-con\-ti\-nu\-ous.
  For each $K>0$ choose function $\sigma_{N,K}^2$ and $\mu_{N,K}$
  which agree with $\sigma_N$ and $\mu_N$, respectively,
  on $[0,K]$ and which
  are uniformly globally Lipschitz continuous and uniformly bounded.
  Existence of such functions follows from the uniform local
  Lipschitz continuity of $\sigma_N^2$ and $\mu_N$.
  Recall the stopping time $\taut_K^N$, $K,N\in\N$, from~\eqref{eq:tautKN}.
  The process
  $\curl{\rub{Z_s^{N,m}(i)}_{s\leq t}\colon i\leq N, m\in\N_0}$
  agrees with the loop-free
  $(N,\sigma^2_{N,K},\mu_{N,K})$-process
  on the event $\{\taut_K^N>t\}$.
  Furthermore the process
  $\curl{\rub{X_s^{N,m}(i)}_{s\leq t}\colon i\leq N, m\in\N_0}$ agrees
  with an $(N,\sigma^2_{N,K},\mu_{N,K})$-process with migration levels
  on $\{\taut_K^N>t\}$.
  The\-re\-fore
  \begin{equation}
    \limsupN\;\E\eckbb{\sum_{i=1}^N\sum_{k=0}^{\infty}
        \absB{X_{t}^{N,k}(i)-Z_{t}^{N,k}(i)}\1_{\taut_K^N>t}}=0
  \end{equation}
  for every $K>0$
  by the preceding step.
  Lemma~\ref{l:tau_theta} handles the event $\{\taut_K^N\leq t\}$.
  This completes the proof of Lemma~\ref{l:X_close_to_Z}.
\end{proof}
\subsection{Proof of Theorem~\ref{thm:convergence}}
\label{ssec:convergence_to_VIM}
\begin{proof}[Proof of Theorem~\ref{thm:convergence}]
  First we prove Theorem~\ref{thm:convergence} under the additional
  Assumption~\ref{a:second_moments}. This will be relaxed later.

  We begin with convergence of fi\-nite-di\-men\-sio\-nal distributions.
  Recall $\MCE_{s,T}$ from~\eqref{eq:MCET}.
  Let $F(\eta)=\prod_{j=1}^n f_j(\eta_{t_j})\in\MCE_{s,T}$ satisfy the
  Lipschitz condition~\eqref{eq:F_is_globally_Lipschitz}
  with Lipschitz constant $L_F\in(0,\infty)$
  and let $F$ be bounded by $C_F<\infty$.
  Furthermore let the function $f\colon I\to\R$ have compact support in $(0,\abs{I})$,
  let $f$ be bounded by $C_f<\infty$ and let $f$
  be globally Lipschitz continuous
  with Lipschitz constant $L_f\in(0,\infty)$.
  Recall the $(N,\mu_N,\sigma_N)$-process with migration levels
  from~\eqref{eq:XN_k}. 
  We will exploit below that all individuals
  on one island have the same migration level in order to show that
  \begin{equation}  \label{eq:same_migration_level_thm1}
    \limN \E F\biggl(\Bigl(\sum_{i=1}^N f\rub{X_t^N(i)}\Bigr)_{t\leq T}\biggr)
    =\limN \E F\rubb{\Bigl(\sum_{i=1}^N
      \sum_{k=0}^\infty f\rub{X_t^{N,k}(i)}\Bigr)_{t\leq T}}.
  \end{equation}
  Assuming~\eqref{eq:same_migration_level_thm1}
  we now prove
  convergence of fi\-nite-di\-men\-sio\-nal distributions.
  We may replace the $(N,\mu_N,\sigma_N^2)$-process with migration levels
  in~\eqref{eq:same_migration_level_thm1} by the loop-free process
  because of Lemma \ref{l:X_close_to_Z} and the Lipschitz continuity
  of $F$ and $f$.
  Hence~\eqref{eq:same_migration_level_thm1} and Lemma~\ref{l:X_close_to_Z}
  imply that
  \begin{equation}  \begin{split}
    \limN \E F\biggl(\Bigl(\sum_{i=1}^N f\rub{X_t^N(i)}\Bigr)_{t\leq T}\biggr)
    &=\limN \E F\rubb{\Bigl(\sum_{i=1}^N
      \sum_{k=0}^\infty f\rub{X_t^{N,k}(i)}\Bigr)_{t\leq T}}\\
    =\limN \E F\rubb{\Bigl(\sum_{i=1}^N
      \sum_{k=0}^\infty f\rub{Z_t^{N,k}(i)}\Bigr)_{t\leq T}}
    &= \E F\rubb{\Bigl(\sum_{(s,\eta)\in\VIM}
          f\rob{\eta_{t-s}}\Bigr)_{t\leq T}}. 
  \end{split}     \end{equation}
  The last equality is the convergence of the loop-free $(N,\mu_N,\sigma_N^2)$-process
  to the virgin island model and has been established in
  Lemma~\ref{l:convergence_of_the_loop_free_process}.
  
  Next we prove~\eqref{eq:same_migration_level_thm1}.
  According to Lemma~\ref{l:decomposition} if we ignore the migration levels in
  the $(N,\mu_N,\sigma_N^2)$-process with migration levels,
  then we obtain a version of the
  $(N,\mu_N,\sigma_N^2)$-process, that is,
  \begin{equation}  \label{eq:decomposition_thm1}
    \limN \E F\rubb{\Bigl(\sum_{i=1}^N f\rub{X_t^N(i)}\Bigr)_{t\leq T}}
    =\limN \E F\rubb{\Bigl(\sum_{i=1}^N f\ruB{\sum_{m=0}^\infty X_t^{N,m}(i)}\Bigr)_{t\leq T}}
  \end{equation}
  For proving~\eqref{eq:same_migration_level_thm1} we observe that
  \begin{equation}  \begin{split}
    \absbb{1-\sum_{k=0}^\infty\1_{x^k\geq\dl}}
    &\leq\1_{x^m\leq\dl\fa m\in\N_0}
       +\1_{\exists m\neq l\colon x^m\geq\dl,x^l\geq\dl}\mal\sum_{k=0}^\infty\1_{x^k\geq\dl}\\
    &\leq\1_{x^m\leq\dl\fa m\in\N_0}
       +\frac{1}{\dl^2}\sum_{k=0}^\infty x^k\sum_{l\neq k}x^l
  \end{split}     \end{equation}
  for every sequence $(x^k)_{k\in\N_0}\subseteq[0,\infty)$ and every $\dl>0$.
  Thus we get that
  \begin{equation}  \begin{split} \label{eq:abab}
      \E&\absbb{ \sum_{i=1}^N f\ruB{\sum_{m=0}^\infty X_t^{N,m}(i)}
          \roB{1-\sum_{k=0}^\infty\1_{X_{t}^{N,k}(i)\geq\dl}}
             }\\
      &\leq L_f\sum_{i=1}^N\E\eckbb{\sum_{m=0}^\infty X_{t}^{N,m}(i)\wedge\dl}
        +C_f\sum_{i=1}^N\frac{1}{\dl^2}
          \E\eckbb{\sum_{k=0}^\infty X_{t}^{N,k}(i)
             \sum_{l\neq k} X_{t}^{N,l}(i)}\\
      &=:C(N,\dl,t)
  \end{split}     \end{equation}
  for all $N\in\N$, $t\geq0$ and $\dl>0$.
  The second summand on the right-hand side converges to zero as $N\to\infty$
  according to Lemma~\ref{l:one_generation_per_island}.
  The first summand on the right-hand side converges to zero as $\dl\to0$
  uniformly in $N\in\N$
  according to Lemma~\ref{l:concentration}.
  Using~\eqref{eq:abab} we obtain that
  \begin{equation}  \begin{split}  \label{eq:proof_of_same_migration_level_thm1}
    \frac{1}{L_F}&\absB{
      \E F\rubb{\Bigl(\sum_{i=1}^N f\ruB{\sum_{m=0}^\infty X_t^{N,m}(i)}\Bigr)_{t\leq T}}
      -\E F\rubb{\Bigl(\sum_{i=1}^N \sum_{k=0}^\infty f\rub{ X_t^{N,k}(i)}\Bigr)_{t\leq T}}}\\
    &\leq\sum_{j=1}^n
      \E\eckbb{\sum_{i=1}^N\sum_{k=0}^\infty\1_{X_{t_j}^{N,k}(i)\geq\dl}
          \Bigl| f\ruB{\sum_{m=0}^\infty X_{t_j}^{N,m}(i)}
               - f\ruB{X_{t_j}^{N,k}(i)}\Bigr|}\\
    &\qquad+\sum_{j=1}^n C(N,\dl,t_j)
      +\sum_{j=1}^n\E\eckbb{\sum_{i=1}^N\sum_{k=0}^\infty\1_{X_{t_j}^{N,k}(i)<\dl}
          \Bigl| f\ruB{X_{t_j}^{N,k}(i)}\Bigr|}\\
    &\leq \frac{L_f}{\dl}\sum_{j=1}^n
      \E\eckbb{\sum_{i=1}^N\sum_{k=0}^\infty X_{t_j}^{N,k}(i)
          \sum_{m\neq k}X_{t_j}^{N,m}(i)}\\
    &\qquad+\sum_{j=1}^n C(N,\dl,t_j)
      +L_f\sum_{j=1}^n \E\eckbb{\sum_{i=1}^N\sum_{k=0}^\infty
          X_{t_j}^{N,k}(i)\wedge\dl}
  \end{split}     \end{equation}
  for all $N\in\N$ and $\dl>0$.
  Letting first $N\to\infty$ and then $\dl\to0$,
  the right-hand side converges to zero according to
  Lemmas~\ref{l:one_generation_per_island}
  and~\ref{l:concentration} and according to the preceding step.
  Inserting this into~\eqref{eq:decomposition_thm1}
  proves~\eqref{eq:same_migration_level_thm1}.

  The next step is to prove tightness. This is analogous to the 
  tightness proof in Lemma~\ref{l:vanishing_immigration_weak_process}.
  Use Lemmas~\ref{l:tau_theta}
  and \ref{l:second_moment_estimate_X} instead of Lemma~\ref{l:Y_2_supsum}.
  So we omit this step. 

  It remains to prove Theorem~\ref{thm:convergence} in the case when
  Assumption~\ref{a:second_moments} fails to hold.
  Fix $T\in[0,\infty)$.
  Let $H\colon D_{\MCM_F(I)}\left([0,T]\right)\to\R$
  be a bounded continuous function on measure-valued \cadlag-paths.
  It follows from Assumption~\ref{a:initial} that $\sum_{i=1}^N X_0^N(i)$
  converges in $L^1$ and thus also in distribution to $\sum_{i\in G}X_0(i)$.
  By the Skorokhod representation theorem
  (e.g.\ Theorem II.86.1 in~\cite{RogersWilliams2000a}),
  there exists a version
  of $\{X_0^N(\cdot)\colon N\in\N\}$ such that $\sum_{i=1}^N X_0^N(i)$
  converges almost surely to $\sum_{i\in G}X_0(i)$ as $N\to\infty$.
  Now the previous step implies that
  \begin{equation}
    \E\eckbb{H\biggl(\biggl(\sum_{i=1}^N \dl_{X_{t}^N(i)}
                     \biggr)_{t\leq T}\biggr)|X_0^N(\cdot)}
    \lraN
    \E\eckbb{H\biggl(\biggl(\sum_{(u,\eta)\in\VIM}\dl_{\eta_{t-u}}
                     \biggr)_{t\leq T}\biggr)|X_0(\cdot)}
  \end{equation}
  almost surely.
  Taking expectations and applying the dominated convergence theorem results in
  \begin{equation}
    \E\eckbb{H\biggl(\biggl(\sum_{i=1}^N \dl_{X_{t}^N(i)}\biggr)_{t\leq T}\biggr)}
    \lraN
    \E\eckbb{H\biggl(\biggl(\sum_{(u,\eta)\in\VIM}\dl_{\eta_{t-u}}\biggr)_{t\leq T}\biggr)}
  \end{equation}
  almost surely. This finishes the proof of Theorem~\ref{thm:convergence}.
\end{proof}

\newpage
\subsection{McKean-Vlasov limit of the $N$-island model as $N\to\infty$}
Proposition~\ref{p:MVL} below establishes convergence of the $N$-island model as $N\to\infty$
in the case of exchangeable initial configurations.
First we prove estimate~\eqref{eq:M.monotone} which implies pathwise uniqueness and monotonicity of the solution.
The following lemma is a special case of Lemma~\ref{l:first_moment_estimate}.
Define $x^+=\max(x,0)$ for all $x\in\R$.
\begin{lemma}  \label{l:Yzeta.estimate}
  Let Assumption~\ref{a:A1} be fulfilled, let $x,y\in I$, $s\in[0,\infty)$
  and let $\zeta,\bar{\zeta}\colon[s,\infty)\to I$ be locally square Lebesgue integrable.
  Then there exists a unique strong solution
  $(Y_{t,s}^{\zeta,x},\bar{Y}_{t,s}^{\bar{\zeta},y})_{t\in[s,\infty)}$
  of
  \begin{align}
    \label{eq:Yzeta}
    dY_{t,s}^{\zeta,x}&=\zeta(t)\,dt-Y_{t,s}^{\zeta,x}dt
                           +\mu\rub{Y_{t,s}^{\zeta,x}}\,dt
         +\sqrt{\sigma^2\rub{Y_{t,s}^{\zeta,x}}}dB_t
    \\
    d\Yb_{t,s}^{\zetab,x}&=\zetab(t)\,dt-\Yb_{t,s}^{\zetab,x}dt
                           +\mu\rub{\Yb_{t,s}^{\zetab,x}}\,dt
         +\sqrt{\sigma^2\rub{\Yb_{t,s}^{\zetab,x}}}dB_t
  \end{align}
  having almost surely continuous paths
  and satisfying $Y_{s,s}^{\zeta,x}=x$, $\bar{Y}_{s,s}^{\bar{\zeta},y}=y$ for all $s\geq0$.
  Moreover
  \begin{equation}
    \E\left[(Y_{t,s}^{\zeta,x}-\bar{Y}_{t,s}^{\bar{\zeta},y})^+\right]
    \leq e^{L_\mu (t-s)}\rubb{\int_s^t\big(\zeta(r)-\zetab(r)\big)^+\,dr+\big(x-y\big)^+}
  \end{equation}
  for all $t\in[s,\infty)$.
\end{lemma}
\begin{lemma}   \label{l:M.monotone}
  Let Assumption~\ref{a:A1} be fulfilled.
  Moreover let $(M_t)_{t\geq0}$ and $(\Mb_t)_{t\geq0}$ be two solutions
  of the McKean-Vlasov equation~\eqref{eq:M} with respect to the same Brownian motion
  having almost surely continuous paths and satisfying $\E[|M_0|+|\bar{M}_0|]<\infty$.
  Then
  \begin{equation}  \label{eq:M.monotone}
    \sup_{t\in[0,T]}\E\left[(M_t-\Mb_t)^+\right]
    \leq e^{L_\mu T+Te^{L_\mu T}}\E\left[(M_0-\Mb_0)^+\right]
  \end{equation}
  for all $T\in[0,\infty)$.
\end{lemma}
\begin{proof}
  Applying Lemma~\ref{l:Yzeta.estimate} and Jensen's inequality yields that
  \begin{equation}  \begin{split}
    \E\left[(M_t-\Mb_t)^+\right]
    &=
    \int\int\E\left[\Big(Y_{t,0}^{E[M],x}-Y_{t,0}^{E[\bar{M}],y}\Big)^+\right]\P(M_0\in dx, \bar{M}_0\in dy)
    \\
    &\leq
       e^{L_\mu T}\rubb{\int_0^t\big(\E[M_r]-\E[\Mb_r]\big)^+\,dr+\E\big[\big(M_0-\bar{M}_0\big)^+\big]}
   \\
    &\leq
       e^{L_\mu T}\int_0^t\E\big[\big(M_r-\Mb_r\big)^+\big]\,dr+e^{L_\mu T}\E\big[\big(M_0-\bar{M}_0\big)^+\big]
  \end{split}     \end{equation}
  for all $t\in[0,T]$ and all $T\in[0,\infty)$.
  Therefore, Gronwall's inequality implies inequality~\eqref{eq:M.monotone}.
\end{proof}
\begin{lemma} \label{l:MVL.existence}
  Let Assumption~\ref{a:A1}
  be fulfilled
  and let $M_0$ be an $I$-valued random variable
  with $\E[|M_0|]<\infty$.
  Then the McKean-Vlasov equation~\eqref{eq:M}
  has a unique strong solution $(M_t)_{t\geq 0}$.
\end{lemma}
\begin{proof}
  Fix a standard Brownian motion $(B_t)_{t\geq0}$ which is independent of $M_0$.
  Let the processes $(Z^{(k)})_{t\geq0}$, $k\in\N_0$, be the unique strong solutions of
  \begin{equation}  \label{eq:Zk}
    dZ_t^{(k)}= \1_{k\geq 1}\E\big[Z_t^{(k-1)}\big]\,dt-Z_t^{(k)} dt +\mu\rub{Z_t^{(k)} }\,dt +\sqrt{\sigma^2\rub{Z_t^{(k)} }}dB_t
  \end{equation}
  and $Z_0^{(k)}=M_0$ for $k\in\N_0$.
  We show by induction on $k\in\N_0$ that $Z_t^{(k)}\leq Z_t^{(k+1)}$ for all $t\in[0,\infty)$ and $k\in\N_0$
  almost surely.
  The base case $k=0$ follows from $\E\big[Z_t^{(0)}\big]\geq0$ for all $t\in[0,\infty)$
  and from a time inhomogeneous version of the monotonicity result in Lemma 3.3 in~\cite{HutzenthalerWakolbinger2007}.
  For the induction step $k\to k+1$, note that the induction hypothesis implies that
  $\E\big[Z_t^{(k)}\big]\leq \E\big[Z_t^{(k+1)}\big]$ for all $t\in[0,\infty)$.
  Thus the induction step follows from
  a time inhomogeneous version of the monotonicity result in Lemma 3.3 in~\cite{HutzenthalerWakolbinger2007}.
  Consequently the process $(M_t)_{t\geq0}$ defined through
  $M_t=\uparrow\lim_{k\to\infty}Z_t^{(k)}$ for $t\in[0,\infty)$
  is a well-defined progressively measurable stochastic process with values in $I\cup\{\infty\}$.
  Note that  Lemma~\ref{l:Yzeta.estimate} implies that
  \begin{equation}
    \E\big[Z_t^{(k)}\big]
    \leq e^{L_\mu t}\Big(\E\big[M_0\big]+\int_0^t\1_{k\geq 1}\E\big[Z_s^{(k-1)}\big]\,ds\Big)
    \leq e^{L_\mu t}\Big(\E\big[M_0\big]+\int_0^t\E\big[Z_s^{(k)}\big]\,ds\Big)
  \end{equation}
  for all $k\in\N_0$ and all $t\in[0,\infty)$.
  By induction on $k\in\N_0$,
  we get that $\E[Z_t^{(k)}]$ is finite for all $k\in\N_0$ and all $t\in[0,\infty)$.
  Therefore the monotone convergence theorem and Gronwall's lemma result in
  \begin{equation}
    \E\big[M_t\big]
    =\E\big[\lim_{k\to\infty}Z_t^{(k)}\big]
    =\lim_{k\to\infty} \E\big[Z_t^{(k)}\big]
    \leq e^{L_\mu T+Te^{L_\mu T}}\E\left[M_0\right]
  \end{equation}
  for all $t\in[0,T]$ and all $T\in[0,\infty)$.
  Next we show that $(M_t)_{t\geq0}$ solves
  the McKean-Vlasov equation~\eqref{eq:M}.
  Define a stopping time $\tau_K:=\inf(\{t\in[0,\infty)\colon M_t\geq K\}\cup\{\infty\})\in[0,\infty]$
  for every $K\in\N$.
  Doob's $L^2$ inequality
  (e.g.\ Theorem II.70.2 in~\cite{RogersWilliams2000a}) implies that
  \begin{equation}  \begin{split}  \label{eq:limit.k}
  &
    \E\left[\int_0^{T\wedge\tau_K}\left|\mu(M_s)-\mu\big(Z_s^{(k)}\big)\right|ds\right]
    +
    \E\left[\sup_{t\in[0,T]}
    \left(\int_0^{t\wedge\tau_K} \sigma(M_s)dB_s
               -\int_0^{t\wedge\tau_K}\sigma\big(Z_s^{(k)}\big) dB_s\right)^2\right]
  \\
  &
    \leq 
    \int_0^T\E\left[\left|\mu(M_{s\wedge\tau_K})-\mu\big(Z_{s\wedge\tau_K}^{(k)}\big)\right|ds\right]
    +4\int_0^T \E\left[\left(\sigma\big(M_{s\wedge\tau_K}\big)-\sigma\big(Z_{s\wedge\tau_K}^{(k)}\big)\right)^2\right]ds
  \\
  &
  \leq
  \left(\sup_{x\neq y\in[0,K]}\tfrac{|\mu(x)-\mu(y)|}{|x-y|}
        +4\sup_{x\neq y\in[0,K]}\tfrac{|\sigma^2(x)-\sigma^2(y)|}{|x-y|}
  \right)
    \int_0^T \E\left[\big|M_{s\wedge\tau_K}-Z_{s\wedge\tau_K}^{(k)}\big|\right]ds
  \end{split}     \end{equation}
  for all $T\in[0,\infty)$, $K\in\N$ and all $k\in\N_0$.
  The right-hand side converges to $0$ as $k\to\infty$ for every $K\in\N$ and every  $T\in[0,\infty)$
  by the dominated convergence theorem.
  Thus, letting $k\to\infty$ in~\eqref{eq:Zk} and using~\eqref{eq:limit.k},
  we conclude that $(M_t)_{t\geq0}$ solves
  the McKean-Vlasov equation~\eqref{eq:M}
  for all $t\in[0,\tau_K]$ almost surely for all $K\in\N$, that is,
  $(M_t)_{t\geq0}$ is a solution of
  the McKean-Vlasov equation~\eqref{eq:M}. This proves existence of a solution.

  Applying Lemma~\ref{l:M.monotone} twice implies that two solutions $(M_t)_{t\geq0}$ and $(\Mb_t)_{t\geq0}$
  of~\eqref{eq:M} with respect to the same Brownian motion and with the same initial point have the same
  finite-dimensional distributions. Together with path continuity this yields pathwise uniqueness
  and -- together with weak existence -- that the SDE~\eqref{eq:M} is exact (Definition V.9.3 in~\cite{RogersWilliams2000b}).
  Moreover Lemma~\ref{l:M.monotone} implies almost sure monotonicity and continuity of the solution
  in the initial point. More precisely if $M_0\leq\Mb_0$ almost surely then $M_t\leq\Mb_t$ almost surely
  for all $t\in[0,\infty)$ and, due to path continuity, $M_t\leq\Mb_t$ for all $t\in[0,\infty)$ almost surely.
  The proof of Theorem V.13.1 in~\cite{RogersWilliams2000b} shows that if the solution of an exact SDE is
  almost surely continuous in the initial point, then there exists a unique strong solution.
  This finishes the proof of Lemma~\ref{l:MVL.existence}.
\end{proof}

Recall the $(N,\mu,\sigma^2)$-island process $(X_t^N)_{t\geq 0}$ from~\eqref{eq:XN}.
Let $\MCM_1(I)$ be the set of probability measures on $I$
equipped with the Prohorov metric, which induces weak convergence
(see e.g.~Theorem 3.3.1 in~\cite{EthierKurtz1986}).
Moreover let
$\MCM_1^2(I)=\{\nu\in\MCM_1(I)\colon \langle v,x^2\rangle=\int x^2\nu(dx)<\infty$ be the set of probability measures on $I$
with finite second moments on which a sequence $(\nu_n)_{n\in\N}$ converges to $\nu\in\MCM_1^2(I)$
if and only if $\nu_n\to\nu$ as $n\to\infty$ in $\MCM_1(I)$
and $\sup_{n\in\N}\langle\nu_n,x^2\rangle<\infty$.
The space $\C([0,\infty),\MCM_1(I))$ is endowed with the topology of uniform convergence
and the subspace $\C([0,\infty),\MCM_1^2(I))$ is endowed with a topology such that
$\nu_n\to\nu$ as $n\to\infty$ in  $\C([0,\infty),\MCM_1^2(I))$
if and only if
$\nu_n\to\nu$ as $n\to\infty$ in  $\C([0,\infty),\MCM_1(I))$
and $\sup_{n\in\N}\sup_{t\in[0,T]}\langle\nu_n(t),x^2\rangle<\infty$ for every $T\in[0,\infty)$,
for details see Appendix B in~\cite{Gaertner1988}.

\begin{proposition} \label{p:MVL}
  Let Assumption~\ref{a:A1}
  be fulfilled
  and let $M_0$ be an $I$-valued random variable.
  Moreover let $X_0^N(j)$, $j\leq N$, be exchangeable random variables with values in $I$
  for every $N\in\N$
  such that $\sup_{N\in\N}\E[(X_0^N(1))^2]<\infty$,
  such that $X_0^N(1)\to M_0$ in distribution as $N\to\infty$
  and
  such that $\tfrac{1}{N}\sum_{j=1}^N \dl_{X_0^N(j)}\to \E[\dl_{M_0}]$
  in distribution in $\MCM_1^2(I)$ as $N\to\infty$.
  Let $(X_t^N)_{t\geq 0}$ be the unique solution of~\eqref{eq:XN}
  and let $(M_t)_{t\geq0}$ be the unique solution of~\eqref{eq:M}.
  Then
  $(X_t^N(i))_{t\geq0}\to (M_t)_{t\geq0}$ in distribution as $N\to\infty$ for every $i\in\N$
  and
  $\big(\tfrac{1}{N}\sum_{i=1}^N \dl_{X_t^N(i)}\big)_{t\in[0,\infty)}
  \to
  \big(\E\big[ \dl_{M_t}\big]\big)_{t\in[0,\infty)}$
  in distribution in $\C([0,\infty),\MCM_1^2)$
  as $N\to\infty$.
  Moreover we have that
  \begin{equation}  \label{eq:L1.XM}
    \sqrt{N}\E\Big[\big|X_t^N(1)-M_t\big|\Big]
    \leq e^{(1+L_\mu)t}
    \left(
    \sqrt{N}
    \E\Big[\big|X_0^N(1)-M_0\big|\Big]
    +\int_0^t\left(\Var\big(M_s)\right)^{\frac{1}{2}}ds
    \right)
    \in[0,\infty)
  \end{equation}
  for all $N\in\N$ and all $t\in[0,\infty)$.
\end{proposition}
\begin{proof}
  Theorem 4.1 of G\"artner (1988)\nocite{Gaertner1988} establishes an
  analogous assertion
  under general assumptions including the ellipticity assumption that (in our notation) $\sigma^2(x)>0$ for all $x\in I$.
  This assumption is not satisfied in our situation.
  Nevertheless,
  parts of the proof of
  Theorem 4.1 of G\"artner (1988)\nocite{Gaertner1988,HutzenthalerTaylor2010}
  carry over to our situation.

  We only prove the assertion for $i=1$, the general case being analogous.
  First we show that
  the distributions of the sequence $(\tfrac{1}{N}\sum_{k=1}^N \dl_{X_t^N(k)},X_t^N(1))_{t\in[0,\infty)}$, $N\in\N$,
  are relatively sequentially weakly compact.
  For this, we introduce more notation.
  Let  $\GenA(\nu)\colon \C^2(I)\to \C(I)$, $\nu\in\MCM_1^2(I)$, be operators defined through
  \begin{equation}  \begin{split}
    \GenA(\nu)f(x)&=f^{'}(x)\left(\int_I z\nu(dz)-x+\mu(x)\right)+\tfrac{1}{2}f^{''}(x)\sigma^2(x)
  \end{split}     \end{equation}
  for all $x\in I$, $f\in \C^2(I)$ and all $\nu\in\MCM_1^2(I)$
  and let the operators $\Gen^N\colon \C^2(I^N)\to \C(I^N)$, $N\in\N$, be defined through
  \begin{equation*}  \begin{split}
    \Gen^N f(x)&=\sum_{k=1}^N \tfrac{\del}{\del x_k}f(x)
                   \bigg(\tfrac{1}{N}\sum_{j=1}^N x_j-x_k+\mu(x_k)\bigg)
                   +\tfrac{1}{2}\sum_{k=1}^N\tfrac{\del^2}{\del x_k^2}f(x)\sigma^2(x_k)
               = \sum_{k=1}^N \GenA_k\bigg( \tfrac{1}{N}\sum_{j=1}^N\dl_{x_j}\bigg)f(x)
  \end{split}     \end{equation*}
  for all $x\in I^N$, $f\in \C^2(I^N)$
  and all $N\in\N$
   where $\GenA_k$ is the operator $\GenA$ acting on the $k$-th variable.
  Note that the $N$-island process~\eqref{eq:XN} solves the well-posed martingale problem for $\Gen^N$ for every $N\in\N$.
  The functions $\mu$ and $\sigma$ are continuous so that
  Assumption (A) of~\cite{Gaertner1988} is satisfied except for $\sigma(0)=0$
  and $\sigma(I)=0$ if $|I|<\infty$.
  Define functions $\ld\colon [0,\infty)\ni x\mapsto 2(1+L_\mu+2\theta+L_\sigma)(1+x)\in(0,\infty)$ and
  $\varphi\colon I\ni x\mapsto 1+x^2\in(0,\infty)$.
  Then $\int_1^\infty dx/\ld(x)=\infty$ and
  $\langle \nu,\GenA(\nu)\varphi\rangle\leq \ld\big(\langle \nu,\varphi\rangle\big)$
  for every $\nu\in\MCM_1^2(I)$ due to Assumption~\ref{a:A1}
  so that
  Assumption (B) of~\cite{Gaertner1988} is satisfied.
  Lemma~\ref{l:average.square.Nisland}
  and the inequality $a\leq 1+a^2$ for $a\in\R$
  yield that
  \begin{equation*}  \begin{split}
    &\lim_{r\to\infty}\sup_{N\in\N}\P\left[\sup_{t\in[0,T]}\tfrac{1}{N}\sum_{k=1}^N \big(X_t^N(k)\big)^2\geq r\right]
    \leq
    \lim_{r\to\infty}\frac{1}{r}\sup_{N\in\N}\E\left[\sup_{t\in[0,T]}\tfrac{1}{N}\sum_{k=1}^N \big(X_t^N(k)\big)^2
           \right]
    \\
    &
    \leq
    \lim_{r\to\infty}\frac{2}{r}\left(
       24T\theta^2
       +8L_\sigma T \Big(2\theta T+1\Big)e^{L_\mu T}+(1+8L_\sigma Te^{L_\mu T})\sup_{N\in\N}\E\left[\big(X_0^N(1)\big)^2\right]
    \right)
    e^{ 40(1+T)(1+L_\mu+L_\sigma)^2T}
    \\&
    =0
  \end{split}     \end{equation*}
  for every $T\in[0,\infty)$.
  This implies assumption (i) of Lemma 1.4 of~\cite{Gaertner1988}.
  The proof of assertion (ii) in the proof of Theorem 1.5 of~\cite{Gaertner1988}
  only requires that $\sup_{x\in I\cap [0,r]}(|\mu(x)|+|\sigma(x)|)<\infty$
  for every $r\in[0,\infty)$ (which follows from Assumption~\ref{a:A1})
  and carries over to our situation without further changes.
  This
  implies assumption (ii) of Lemma 1.4 of~\cite{Gaertner1988}.
  Lemma 1.4 of~\cite{Gaertner1988} thus yields that
  the distributions of
  the sequence
  $\big(\tfrac{1}{N}\sum_{k=1}^N \dl_{X_t^N(k)}\big)_{t\in[0,\infty)}$, $N\in\N$,
  are relatively sequentially weakly compact 
  in $C([0,\infty),\MCM_1^2(I))$.
  In addition,
  the proof of Theorem 4.1 of~\cite{Gaertner1988}
  shows
  relative weak compactness of $(X_t^N(1))_{t\geq0}$, $N\in\N$, in $\C([0,\infty),I)$
  (without using Assumptions (C) or (E) in~\cite{Gaertner1988}).
  This proves that
  the distributions of the sequence
  $\big(X_t^N(1),\tfrac{1}{N}\sum_{k=1}^N \dl_{X_t^N(k)}\big)_{t\in[0,\infty)}$, $N\in\N$,
  are relatively sequentially weakly compact
  in $C([0,\infty),I\times \MCM_1^2(I))$.

  Next we identify the weak limit of the sequence
  $\big(X_t^N(1),\tfrac{1}{N}\sum_{k=1}^N \dl_{X_t^N(k)}\big)_{t\in[0,\infty)}$, $N\in\N$.
  Let $(M_t(j))_{t\geq0}$ be the unique strong solution of~\eqref{eq:M} with
  respect to the Brownian motion $(B_t(j))_{t\geq0}$ for every $j\in\N$
  such that $M_0(j)$, $j\in\N$, are independent copies of $M_0$.
  As in Theorem 1 of Yamada and Watanabe (1971)~\cite{YamadaWatanabe1971}, an approximation
  of the function
  $\R\ni x\to |x|\in [0,\infty)$ with $\C^2$-functions
  (see also the proof of Lemma~\ref{l:X_close_to_Z} for this approximation)
  results in
  \begin{equation}  \begin{split}
    d\big|X_t^N(j)-M_t(j)\big|
    &=\sgn\rob{X_t^N(j)-M_t(j)}
       d\rub{X_t^N(j)-M_t(j)}
  \end{split}     \end{equation}
  for all $j\in\{1,\ldots,N\}$ and all $N\in\N$.
  Taking expectations and using Assumption~\ref{a:A1} yields that
  \begin{equation}  \begin{split}
  \lefteqn{
    \E\bigg[\frac{1}{N}\sum_{j=1}^N\big|X_t^N(j)-M_t(j)\big|\bigg]
    -
    \E\bigg[\frac{1}{N}\sum_{j=1}^N\big|X_0^N(j)-M_0(j)\big|\bigg]
  }\\
    &\leq \frac{1}{N}\sum_{j=1}^N\int_0^t\E\bigg[
        \Big|\frac{1}{N}\sum_{k=1}^N X_s^N(k)-\E\big[M_s(j)\big]\Big|
        +\sgn\rob{X_s^N(j)-M_s(j)}\Big(\mu\big(X_s^N(j)\big) -\mu\big(M_s(j)\big)\Big)
       \bigg]ds
    \\
    &\leq
    \int_0^t\E\Big[ \Big|\frac{1}{N}\sum_{k=1}^N M_s(k)-\E\big[M_s(1)\big]\Big| \Big]ds
    +(1+L_\mu)\int_0^t\E\bigg[\frac{1}{N}\sum_{j=1}^N\Big| X_s^N(j)-M_s(j)\Big|\bigg]ds
  \end{split}     \end{equation}
  for all $N\in\N$ and all $t\in[0,\infty)$.
  Therefore, Gronwall's inequality implies that
  \begin{equation}  \begin{split}\label{eq:Gronwall.L1.XM}
  \lefteqn{
    \E\bigg[\big|X_t^N(1)-M_t(1)\big|\bigg]
    =
    \E\bigg[\frac{1}{N}\sum_{j=1}^N\big|X_t^N(j)-M_t(j)\big|\bigg]
  }
  \\
  &
    \leq e^{(1+L_\mu)t}
    \bigg(
    \E\bigg[\frac{1}{N}\sum_{j=1}^N\big|X_0^N(j)-M_0(j)\big|\bigg]
    +
    \int_0^t\E\Big[ \Big|\frac{1}{N}\sum_{k=1}^N M_s(k)-\E\big[M_s(1)\big]\Big| \Big]ds
    \bigg)
  \\
  &
    \leq e^{(1+L_\mu)t}
    \bigg(
    \E\Big[\big|X_0^N(1)-M_0(1)\big|\Big]
    +
    \int_0^t
    \bigg(\Var\Big(\frac{1}{N}\sum_{k=1}^N M_s(k)\Big)\bigg)^{\frac{1}{2}}ds
    \bigg)
  \\
  &
    = e^{(1+L_\mu)t}
    \bigg(
    \E\Big[\big|X_0^N(1)-M_0(1)\big|\Big]
    +\frac{1}{\sqrt{N}}
    \int_0^t
    \Big(\Var\big(M_s(1)\big)\Big)^{\frac{1}{2}}ds
    \bigg)
  \end{split}     \end{equation}
  for all $N\in\N$ and all $t\in[0,\infty)$.
  This proves inequality~\eqref{eq:L1.XM}.
  Moreover we infer that
  \begin{equation}  \begin{split}
    &
    \left|
      \E\left[e^{
          -\sum_{j=1}^n \ld_j X_{t_j}^N(1)
          -\sum_{j=1}^n\ldt_j \frac{1}{N}\sum_{k=1}^N f_j\big( X_{t_j}^N(k)\big)
       }
       \right]
      -
      \E\left[e^{
          -\sum_{j=1}^n \ld_j M_{t_j}(1)
          -\sum_{j=1}^n\ldt_j \E\big[f_j\big( M_{t_j}(1) \big)\big]
       }
       \right]
    \right|
     \\
     &\leq
       \sum_{j=1}^n \ld_j \E\left[\left|X_{t_j}^N(1)-M_{t_j}(1)\right|\right]
       +
       \sum_{j=1}^n \ldt_j \E\left[
                \Big|\frac{1}{N}\sum_{k=1}^N f_j\big(X_{t_j}^N(k)\big)-\E\Big[f_j\big(M_{t_j}(1)\big)\Big]\Big|\right]
     \\
     &\leq
       \sum_{j=1}^n \Big(\ld_j +\ldt_j \sup_{x\neq y\in I}\tfrac{|f_j(x)-f_j(y)|}{|x-y|}
       \Big)
         \E\left[\left|X_{t_j}^N(1)-M_{t_j}(1)\right|\right]
       +
       \sum_{j=1}^n \ldt_j \tfrac{1}{\sqrt{N}}
       \left(\Var\Big( f_j\big(M_{t_j}(1)\big)\Big)\right)^{\frac{1}{2}}
  \end{split}     \end{equation}
  for all $N\in\N$,
  $0\leq t_1<t_2<\ldots<t_n<\infty$,
  $\ld_1,\ldots,\ld_n,\ldt_1,\ldots\ldt_n\in[0,\infty)$,
  all globally Lipschitz continuous functions $f_1,\ldots,f_n\colon I\to[0,\infty)$
  and all
  $n\in\N$.
  The right-hand side converges to $0$ as $N\to\infty$ due to inequality~\eqref{eq:Gronwall.L1.XM}.
  This identifies the limit and, together with tightness, implies 
  that the sequence
  $\big(X_t^N(1),\tfrac{1}{N}\sum_{k=1}^N \dl_{X_t^N(k)}\big)_{t\in[0,\infty)}$, $N\in\N$,
  converges to $\big(M_t,\E\big[ \dl_{M_t}\big]\big)_{t\in[0,\infty)}$
  in distribution in $C([0,\infty),I\times \MCM_1^2(I))$.
  This finishes the proof of Proposition~\ref{p:MVL}.
\end{proof}
\section{Comparison with the virgin island model}
\label{sec:comparison_with_the_VIM}
As in Section~\ref{sec:convergence_to_VIM}, we define a loop-free process.
Let 
$\curlb{(Z_t^{(k)}(i))_{t\geq0}\colon k\in\N_0, i\in G}$
be the solution of
\begin{equation}  \begin{split}   \label{eq:ZG}
  dZ_t^{(k)}(i)
    =&\rubb{\sum_{j\in G} Z_t^{(k-1)}(j)m(j,i)-Z_t^{(k)}(i)
         +\mu\rub{Z_t^{(k)}(i)}}\,dt\\
     & +\sqrt{\sigma^2\rub{Z_t^{(k)}(i)}}\,dB_t^k(i),\quad
       Z_0^{(k)}(i)=\1_{k=0}X_0(i),\qquad i\in G,k\in\N_0,
\end  {split}     \end  {equation}
where  we agree on $Z_t^{(-1)}(i):=0$ for $t\geq0$ and $i\in G$.
We will refer to this process as the loop-free $(G,m,\mu,\sigma^2)$-process.
The main two steps in our proof of
Theorem~\ref{thm:comparison} are as follows.
Lemma~\ref{l:Z_dominates_X} below shows that the total mass of the
$(G,m,\mu,\sigma^2)$-process is dominated by the total mass of the loop-free
$(G,m,\mu,\sigma^2)$-process. 
Lemma~\ref{l:V_dominates_Z} then proves that the total mass of the
loop-free $(G,m,\mu,\sigma^2)$-process is dominated by the total mass of
the virgin island model.
Our proof of Lemma~\ref{l:V_dominates_Z} exploits the hierarchical
structure of the loop-free process. Note that conditioned on migration level $k-1$,
the islands with migration level $k$ are independent one-dimensional
diffusions. We prepare this in
Subsection~\ref{ssec:Decomposition of a one-dimensional diffusion with immigration into subfamilies}
by studying the one-dimensional time-inhomogeneous diffusion
\begin{equation}  \label{eq:Y_immi}
  dY_{t,s}^{\zeta,x}(i)=\zeta(t)\,dt-Y_{t,s}^{\zeta,x}(i)dt
                         +\mu\rub{Y_{t,s}^{\zeta,x}(i)}\,dt
       +\sqrt{\sigma^2\rub{Y_{t,s}^{\zeta,x}(i)}}dB_t(i)
\end{equation}
where $Y_{s,s}^{\zeta,x}=x\in I$ and $s\geq0$.
The path $\zeta\in\C\rob{[0,\infty),I}$
will later represent the mass immigrating from lower migration levels.

The core of the comparison result is the following generator calculation
which manifests the intuition that separating mass onto different
islands increases the total mass.
If $\zeta\equiv c$ is constant, then a formal generator of
$(Y_{t,s}^{c,\cdot})_{t\geq s}$ is
\begin{equation}
  \Gen^{(c)}f(x):=\rob{c-x+\mu(x)}f^{'}(x)+\frac12\sigma^2(x)f^{''}(x)\quad x\in I
\end{equation}
where $f\in\C^2(I)$, see e.g.\ Section 5.3 in \cite{EthierKurtz1986}.
Recall $\MCF_{++}^{(1)}$ from~\eqref{eq:F++}.
%
\begin{lemma}  \label{l:generator_estimate}
  Assume~\ref{a:A1}.
  Suppose that $c,c_1,c_2\in I$ satisfy $c\leq c_1+c_2$.
  Assume $\mu$ to be subadditive, that is,
  $\mu(x+y)\leq \mu(x)+\mu(y)$ for all $x,y\in I$ with
  $x+y\in I$.
  Let $x,y,x+y\in I$.
  If $\sigma^2$ is superadditive, then
  \begin{equation}  \label{eq:generator_estimate_superadditive}
    \Gen^{(c)} f(x+y)\leq \rub{\Gen^{(c_1)} f(\cdot+y)}(x)+\rub{\Gen^{(c_2)} f(x+\cdot)}(y)
    \ \fa f\in\MCF_{+-}^{(1)}\cap\C^2.
  \end{equation}
  If $\sigma^2$ is subadditive, then
  \begin{equation}  \label{eq:generator_estimate_subadditive}
    \Gen^{(c)} f(x+y)\leq \rub{\Gen^{(c_1)} f(\cdot+y)}(x)+\rub{\Gen^{(c_2)} f(x+\cdot)}(y)
    \ \fa f\in\MCF_{++}^{(1)}\cap\C^2.
  \end{equation}
  If $\sigma^2$ is additive, then
  \begin{equation}  \label{eq:generator_estimate_additive}
    \Gen^{(c)} f(x+y)\leq \rub{\Gen^{(c_1)} f(\cdot+y)}(x)+\rub{\Gen^{(c_2)} f(x+\cdot)}(y)
    \ \fa f\in\MCF_{+\pm}^{(1)}\cap\C^2.
  \end{equation}
\end{lemma}
\begin{proof}
  For $f\in\MCF_{+-}^{(1)}\cap\C^2$, the first derivative is non-negative and the
  second derivative is nonpositive. Thus
  \begin{equation}  \begin{split}
    \Gen^{(c)} f(x+y)
    &=f^{'}(x+y)\ruB{c-(x+y)+\mu(x+y)} +\frac12 f^{''}(x+y)\sigma^2(x+y)\\
   &\leq f^{'}(x+y)
     \ruB{c_1-x+\mu(x)}+\frac12 f^{''}(x+y)\sigma^2(x)\\
   &\quad+ f^{'}(x+y) \ruB{c_2-y+\mu(y)} +\frac12 f^{''}(x+y)\sigma^2(y)\\
   &= \rub{\Gen^{(c_1)} f(\cdot+y)}(x)+\rub{\Gen^{(c_2)} f(x+\cdot)}(y).
  \end{split}     \end{equation}
  This is inequality~\eqref{eq:generator_estimate_superadditive}.
  The proof of inequality~\eqref{eq:generator_estimate_subadditive}
  is analogous.
  If $\sigma^2$ is additive, then $\sigma^2(x+y)=\sigma^2(x)+\sigma^2(y)$
  and no property of $f^{''}$ is needed in the above calculation.
\end{proof}
\noindent
As a remark, we observe that the operator on the right-hand side
of~\eqref{eq:generator_estimate_superadditive} is a formal generator
of the superposition $\rob{Y_{t,s}^{c_1,x}+\Yt_{t,s}^{c_2,y}}_{t\geq s}$
of two independent solutions of~\eqref{eq:Y_immi}.
This follows from Theorem 4.10.1 in~\cite{EthierKurtz1986}.

We will lift inequality~\eqref{eq:generator_estimate_superadditive} between
formal generators to an inequality between the associated semigroups.
For this we use an integration by parts formula.
For its formulation, let $\Gen_S$ and $\Gen_T$ be two
generators associated with the semigroups $(S_t)_{t\geq0}$
and $(T_t)_{t\geq0}$, respectively.
Then, for $t\in[0,\infty)$, we have that
\begin{equation}       \label{eq:comp_int_by_parts}
  S_t f-T_t f
  =\int_0^t T_{t-s}\rub{\Gen_S-\Gen_T}S_s f\,ds
\end{equation}
if $S_s f\in\MCD(\Gen_S)\cap\MCD(\Gen_T)$ for all $s\leq t$,
see p.\ 367 in Liggett \nocite{Liggett1985}(1985).
The idea of using~\eqref{eq:comp_int_by_parts} for a comparison
is borrowed from Cox et al.\ \nocite{CoxEtAl1996}(1996).
As the generator
inequality~\eqref{eq:generator_estimate_superadditive} 
holds for functions in $\MCF_{+-}^{(1)}$, we need to show that the semigroup
of $(Y_{t,s}^{x,\cdot})_{t\geq s}$ preserves $\MCF_{+-}^{(1)}$.
This is subject of the following subsection.
\subsection{Preservation of convexity}%
\label{ssec:Preservation of convexity}
We write $\uline{x}_n:=(x_1,\ldots,x_n)$
for $x_1,\ldots,x_n\in\R$ and $n\in\N$.
The $i$-th unit row vector is denoted as
$e_i$ for every $i\in\N$.
Recall $\MCF_{++}^{(n)}$ from~\eqref{eq:F++}.
\begin{lemma}  \label{l:delete}
  For every $n\in\N$ and $f\in\MCF_{++}^{(n)}\rob{[0,\infty)}$, we have that
  \begin{equation} \label{eq:delete}
     f\ruB{\uline{x}_n+\sum_{j=1}^k h_j e_{i_j}}
    -f\ruB{\uline{x}_n+\sum_{j=2}^k h_j e_{i_j}}
    \geq
     f\ruB{\uline{x}_n+h_1 e_{i_1}}
    -f\rub{\uline{x}_n}
  \end{equation}
  for all $\uline{x}_n\in [0,\infty)^n$, all $h_1,\ldots,h_k\geq0$,
  all $(i_1,\ldots,i_k)\in\{1,\ldots,n\}^k$ and all $k\in\N$.
  The reverse inequality holds if $\MCF_{++}^{(n)}$ is replaced by
  $\MCF_{+-}^{(n)}$.
\end{lemma}
\begin{proof}
  The proof is by induction on $k\in\N$.
  The base case $k=1$ is trivial.
  Now assume that~\eqref{eq:delete} holds for some $k\geq 1$.
  Fix $\uline{x}_n\in[0,\infty)^n$,
  $h_1,\ldots,h_k\geq0$ and 
  $(i_1,\ldots,i_k)\in\{1,\ldots,n\}^k$.
  Applying the induction hypothesis at location
  $\uline{x}_n+\sum_{j=2}^k h_j e_{i_j}$ to the index tuple $(i_1,i_{k+1})$,
  we obtain that
  \begin{equation}  \begin{split}
    &f\ruB{\uline{x}_n+\sum_{j=2}^k h_j e_{i_j}+h_1 e_{i_1}+h_{k+1}e_{i_{k+1}}}
    -f\ruB{\uline{x}_n+\sum_{j=2}^k h_j e_{i_j}+h_{k+1}e_{i_{k+1}}}\\
    &\geq
     f\ruB{\uline{x}_n+\sum_{j=2}^k h_j e_{i_j}+h_1 e_{i_1}}
    -f\ruB{\uline{x}_n+\sum_{j=2}^k h_j e_{i_j}}
    \geq
     f\ruB{\uline{x}_n+h_1 e_{i_1}}
    -f\rub{\uline{x}_n}.
  \end{split}     \end{equation}
  The last step is again the induction hypothesis.
\end{proof}
\begin{lemma}   \label{l:double_argument}
  Let $n\in\N$, $c\in[0,\infty)$
  and $f\in\MCF_{++}^{(n+1)}\rob{[0,\infty)}$.
  Then the two functions
  \begin{equation}
    \ft\colon[0,\infty)^n\to\R, \uline{x}_n\to f(x_1,\ldots,x_n,x_n)
  \end{equation}
  and
  \begin{equation}
    \fb\colon[0,\infty)^n\to\R,  \uline{x}_n\to f(x_1,\ldots,x_n,c)
  \end{equation}
  are elements of $\MCF_{++}^{(n)}\rob{[0,\infty)}$.
  This is also true if $\MCF_{++}$ is replaced by $\MCF_{+-}$ and
  $\MCF_{+\pm}$, respectively.
\end{lemma}
\begin{proof}
  The functions $\ft$ and $\fb$ are non-decreasing and either bounded
  or non-negative.
  It is clear that $\fb$ is again $(i,j)$-convex
  for $1\leq i,j\leq n$ and that $\ft$ is $(i,j)$-convex for $1\leq i,j\leq n-1$.
  It remains to prove $(i,n)$-convexity of $\ft$ for $1\leq i\leq n$.
  Applying Lemma~\ref{l:delete} at location $\uline{x}_{n+1}$ to the index
  tuple $(i,n,n+1)$, we obtain for all $h_1,h_2\geq0$ that
  \begin{equation}  \begin{split}
    &f\rob{\uline{x}_n+h_1 e_i+h_2 e_n, x_n+h_2}
      =f\rob{\uline{x}_{n+1}+h_1 e_i+h_2 e_n+h_2 e_{n+1}}|_{x_{n+1}=x_n}\\
    &\geq
    \roB{ 
       f\rob{\uline{x}_{n+1}+h_2 e_n+h_2 e_{n+1}}
      +f\rob{\uline{x}_{n+1}+h_1 e_i}
      -f\rob{\uline{x}_{n+1}}
        }|_{x_{n+1}=x_n}\\
    &=
      f\rob{\uline{x}_n+h_2 e_n, x_n+h_2}
     +f\rob{\uline{x}_n+h_1 e_i, x_n}
     -f\rob{\uline{x}_n, x_n},
  \end{split}     \end{equation}
  that is, $\ft$ is $(i,n)$-convex.
\end{proof}
\begin{lemma}   \label{l:preservation}
  Assume~\ref{a:A1}.
  Let $n\in\N$ and $f\in\MCF_{+\pm}^{(n+1)}\rob{[0,\infty)}$.
  Then the function
  \begin{equation}
    (x_1,\ldots,x_n)\mapsto 
    \E f\rob{x_1,\ldots,x_n,Y_{t,s}^{\zeta,x_n}}
  \end{equation}
  is an  element of $\MCF_{+\pm}^{(n)}\rob{[0,\infty)}$
  for every $0\leq s\leq t$. 
  If $\mu$ is concave, then this property still holds if 
  $\MCF_{+\pm}^{(n+1)}$ is replaced by $\MCF_{++}^{(n+1)}$
  and if $\MCF_{+\pm}^{(n)}$ is replaced by $\MCF_{++}^{(n)}$, respectively.
\end{lemma}
\begin{proof}
  Fix $0\leq s\leq t$ and $n\in\N$.
  We only prove the case of $\mu$ being concave and $f\in\MCF_{++}^{(n+1)}$
  as the remaining cases are similar.
  According to Lemma~\ref{l:double_argument}, it suffices to prove that
  \begin{equation}
    \ft\colon[0,\infty)^{n+1}\to\R, (x_1,\ldots,x_{n+1})\mapsto 
    \E f\rob{x_1,\ldots,x_n,Y_{t,s}^{\zeta,x_{n+1}}}
  \end{equation}
  is an  element of $\MCF_{++}^{(n+1)}\rob{[0,\infty)}$.
  Let $(Y_{t,s}^{\zeta,x})_{t\geq s}$, $x\in I$, be solutions
  of~\eqref{eq:Y_immi} with respect to the same Brownian motion.
  It is known that $Y_{t,s}^{\zeta,x}\leq Y_{t,s}^{\zeta,x+h}$
  holds
  almost surely for all $x\leq x+h\in I$,
  see e.g. Theorem V.43.1 in~\cite{RogersWilliams2000b} for the time-homogeneous case.
  Thus the function $\ft$ is again non-decreasing.
  Moreover $\ft$ inherits $(i,j)$-convexity from $f$ for
  every $1\leq i,j\leq n$.
  It remains to show  that $\ft$ is $(i,n+1)$-convex for
  $1\leq i\leq n+1$.
  If $i\leq n$, then $(i,n+1)$-convexity of $f$ at the point
  $\rob{x_1,\ldots,x_n,Y_{t,s}^{\zeta,x_{n+1}}}$ implies that
  \begin{equation}  \begin{split}
    &\ft\rob{\uline{x}_{n+1}+h_1 e_i+h_2 e_{n+1}}
    =\E f\rob{\uline{x}_{n}+h_1 e_i,Y_{t,s}^{\zeta,x_{n+1}}
    +Y_{t,s}^{\zeta,x_{n+1}+h_2}
    -Y_{t,s}^{\zeta,x_{n+1}}}\\
    &\geq
     \E f\rob{\uline{x}_{n}+h_1 e_i, Y_{t,s}^{\zeta,x_{n+1}}}
    +\E f\rob{\uline{x}_{n}, Y_{t,s}^{\zeta,x_{n+1}+h_2}}
    -\E f\rob{\uline{x}_{n}, Y_{t,s}^{\zeta,x_{n+1}}}\\
   &=
     \ft\rob{\uline{x}_{n+1}+h_1 e_i}
     +\ft\rob{\uline{x}_{n+1}+h_2 e_{n+1}}
     -\ft\rob{\uline{x}_{n+1}}
  \end{split}     \end{equation}
  for every $h_1,h_2\geq0$ and $(x_1,\ldots,x_n)\in[0,\infty)^n$,
  that is,
  $(i,n+1)$-convexity of $\ft$ in the case $i\leq n$.
  One can establish convexity of
  \begin{equation}
    y\mapsto \ft(x_1,\ldots,x_n,y)
  \end{equation}
  as in Lemma 6.1 of~\cite{HutzenthalerWakolbinger2007} (this Lemma 6.1 shows
  concavity if $f$ is $(n+1,n+1)$-concave and smooth; for the general case, approximate
  $f$, $\mu$ and $\sigma$ with smooth functions and exploit that convexity is preserved under pointwise limits).
  This step uses concavity of $\mu$.
  Consequently, $\ft$ is $(n+1,n+1)$-convex.
  This completes the proof of $\ft\in\MCF_{++}^{(n+1)}\rob{[0,\infty)}$.
\end{proof}
Lemma~\ref{l:preservation} extends Proposition 16 of Cox et al.\ (1996) \cite{CoxEtAl1996}.
This Proposition 16 is used in~\cite{CoxEtAl1996} to establish a comparison result
between diffusions with different diffusion functions, see Theorem 1 
in~\cite{CoxEtAl1996}.
Using the above Lemma~\ref{l:preservation}, this comparison result can be extended
to more general test functions.
\subsection{Decomposition of a one-dimensional diffusion with immigration into subfamilies}
\label{ssec:Decomposition of a one-dimensional diffusion with immigration into subfamilies}

Feller's branching diffusion
with immigration can be decomposed into independent families which originate either from an individual
at time zero or from an immigrant,
see e.g.\ Theorem 1.3 in Li and Shiga (1995)\nocite{LiShiga1995}.
A diffusion does in general not agree with its family decomposition
if individuals interact with each other,
e.g.\ if the branching rate depends on the population size.
If the drift function is subadditive and if the branching function is superadditive,
however, then
we get at least a comparison result.
In that situation,
the diffusion is dominated by its family decomposition.
More precisely, the total mass increases in the order $\leqFix$
if we let all subfamilies evolve independently, see Lemma~\ref{l:family_decomposition}
below. The following lemma 
is a first step in this direction.
%
\begin{lemma}    \label{l:semigroup_estimate}
  Assume~\ref{a:A1}.
  Let $x,y,x+y\in I$ and let $\zeta,\zetat\colon[0,\infty)\to[0,\infty)$
  be locally Lebesgue integrable.
  If $\mu$ is concave and $\sigma^2$ is superadditive, then
  \begin{equation}   \label{eq:semigroup_estimate_ix}
    \ruB{Y_{t,s}^{\zeta+\zetat,x+y}}_{t\geq s}
    \leq_{\SetF_{+-}\ro{[0,\infty)}}
    \ruB{Y_{t,s}^{\zeta,x}+\Yt_{t,s}^{\zetat,y}}_{t\geq s},
    \qquad\fa s\geq0,
  \end{equation}
  where $\big(Y_{t,s}^{\zeta,x}\big)_{t\geq s}$
  and
  $\big(\tilde{Y}_{t,s}^{\tilde{\zeta},x}\big)_{t\geq s}$
  are independent processes.
  If $\mu$ is concave and $\sigma^2$ is subadditive, then
  inequality~\eqref{eq:semigroup_estimate_ix}
  holds with $\SetF_{+-}$ replaced by $\SetF_{++}$.
  If $\mu$ is subadditive and $\sigma^2$ is additive, then
  inequality~\eqref{eq:semigroup_estimate_ix}
  holds with $\SetF_{+-}$ replaced by $\SetF_{+\pm}$.
\end{lemma}
\begin{proof}
  Let $F_n(\eta)=f_n(\eta_{t_1},\ldots,\eta_{t_n})\in\SetF_{+-}([0,\infty))$
  where $f_n\in\MCF_{+-}^{(n)}([0,\infty))$.
  We begin with the case of $\zeta,\zetat$ being simple functions.
  W.l.o.g.\ we consider
  $\zeta(t)=\sum_{i=1}^n c_i\1_{[t_{i-1},t_i)}(t)$
  and
  $\zetat(t)=\sum_{i=1}^{n} \ct_i\1_{[t_{i-1},t_i)}(t)$
  where $c_1,...,c_n,\ct_1,\ldots,\ct_n\geq0$, $t_0=s$ and $t_{n+1}=\infty$
  as we may let $F_n$ depend trivially on further time points.
  We will prove by induction on $n\in\N$ that
  \begin{equation}  \label{eq:semigroup_inequality_F}
    \E F_n\ruB{\rub{Y_{t,s}^{\zeta+\zetat,x+y}}_{t\geq s}}
    \leq
    \E F_n\ruB{\rub{Y_{t,s}^{\zeta,x}+\Yt_{t,s}^{\zetat,y}}_{t\geq s}}.
  \end{equation}
  For the base case $n=1$ additionally assume $f_1\in\mathbf{C}^2$.
  Approximate ${\sigma}$ and $\mu$ with functions
  ${\sigma_l},\mu_l\in\mathbf{C}^{\infty}(\R)$ having the following properties.
  All derivatives $\sigma_l^{(k)}$,$\mu_l^{(k)}$, $k\in\N_0$,
  are bounded, $\mu_l$ is concave and $\sigma_l^2$ is superadditive.
  Both functions vanish at zero.
  If $\abs{I}<\infty$, then $\mu_l(\abs{I})\leq0=\sigma_l^2(\abs{I})$.
  Moreover $\mu_l(x)\to \mu(x)$ and
  $\sigma^2_l(x)\to \sigma^2(x)$ as $l\to\infty$ for all $x\in I$.
  Let $(Y_{t,s}^{\zeta,x,l})_{t\geq s}$, $x\in I$, be solutions
  of~\eqref{eq:Y_immi} with $\sigma^2$ and $\mu$ replaced by
  $\sigma_l^2$ and $\mu_l$, respectively,
  and let 
  $(\Yt_{t,s}^{\zetat,y,l})_{t\geq s}$, $y\in I$, be an independent version hereof.
  Then $x\mapsto S_{t} f_1(x):=\E f_1\rub{Y_{t,s}^{c_1,x,l}}$
  is twice continuously differentiable for every $t\geq s$,
  see Theorem 8.4.3 in Gikhman and Skorokhod \nocite{GikhmanSkorokhod1969}(1969).
  In addition,
  Lemma~\ref{l:preservation} proves
  $S_{t} f_1\in\MCF_{+-}^{(1)}$ for all $t\in[s, t_1]$.
  Consequently, we may apply Lemma~\ref{l:generator_estimate} to
  $S_{t} f_1\in\MCF_{+-}^{(1)}\cap\C^2$
  for every $t\in[s,t_1]$
  and the integration by parts formula~\eqref{eq:comp_int_by_parts}
  yields that
  \begin{equation}   \label{eq:semigroup_estimate_n}
    \E\eckB{f_1\rub{Y_{t_1}^{c_1+\ct_1,x+y,l}}}
    \leq
    \E\eckB{f_1\rub{Y_{t_1}^{c_1,x,l}+\Yt_{t_1}^{\ct_1,y,l}}}.
  \end{equation}
  Now as $l\to\infty$, $(Y_t^{c_1,x,l})_{t\geq0}$ converges weakly
  to $(Y_t^{c_1,x})_{t\geq0}$ for every $x\in I$,
  see Lemma 19 in Cox et al.\ \nocite{CoxEtAl1996}(1996) for a sketch of the proof.
  Therefore letting $l\to\infty$ in~\eqref{eq:semigroup_estimate_n}
  proves~\eqref{eq:semigroup_inequality_F} for $n=1$ if
  $f_1\in\C^2$. The case of general $f_1\in\MCF_{+-}^{(1)}$ follows by
  approximating $f_1$ with smooth functions in $\MCF_{+-}^{(1)}$.
  For the induction step $n\to n+1$, define
  \begin{equation}
    \ft_{n}(y_1,\ldots,y_n)
    :=\E f_{n+1}\roB{y_1,\ldots,y_n,Y_{t_{n+1},t_n}^{\zeta+\zetat,y_n}}
    \quad\fa (y_1,\ldots,y_n)\in I^n.
  \end{equation}
  Note that the induction hypothesis implies that
  \begin{equation}
    \ft_{n}(x_1+y_1,\ldots,x_{n}+y_{n})
    \leq\E f_{n+1}\roB{x_1+y_1,\ldots,x_n+y_n,
            Y_{t_{n+1},t_n}^{\zeta,x_n}+ \Yt_{t_{n+1},t_n}^{\zetat,y_n}}
  \end{equation}
  and that Lemma~\ref{l:preservation} implies that $\ft_{n}\in\MCF_{+-}^{(n)}$.
  Therefore, using the Markov property and the induction hypothesis,
  we get that
  \begin{equation}  \begin{split}
    &\E F_{n+1}\ruB{Y_{\cdot,s}^{\zeta+\zetat,x+y}}
    =\E \ft_{n}\ruB{
            Y_{t_1,s}^{\zeta+\zetat,x+y},\ldots,
            Y_{t_n,s}^{\zeta+\zetat,x+y}}\\
    &\leq
    \E \ft_{n}\ruB{
            Y_{t_1,s}^{\zeta,x} +\Yt_{t_1,s}^{\zetat,y},\ldots,
            Y_{t_n,s}^{\zeta,x}+\Yt_{t_n,s}^{\zetat,y}}\\
    &\leq\E\eckbb{
      \E \eckbb{f_{n+1}\Bigl.\roB{x_1+y_1,\ldots,x_n+y_n,
            Y_{t_{n+1},t_n}^{\zeta,x_n}+ \Yt_{t_{n+1},t_n}^{\zetat,y_n}}
            \Bigr|_{x_i=Y_{t_i,s}^{\zeta,x},\;y_i=\Yt_{t_i,s}^{\zetat,y}}}}\\
    &=\E F_{n+1}\rubb{\ruB{Y_{t,s}^{\zeta,x}+\Yt_{t,s}^{\zetat,y}}_{t\geq s}}
  \end{split}     \end{equation}
  for all $x,y\in I$ satisfying $x+y\in I$.
  The last step follows from the Markov property and from independence of the
  two processes. This proves~\eqref{eq:semigroup_inequality_F}.

  In case of general functions $\zeta$ and $\zetat$, approximate $\zeta$ and $\zetat$
  with simple functions $\zeta_l$ and $\zetat_l$, $l\in\N$, respectively.
  The process $Y_{\cdot,s}^{\zeta_l,x}$ converges in the sense of
  fi\-nite-di\-men\-sio\-nal distributions in $L^1$,
  see Lemma~\ref{l:second_moment_estimate_Y},
  and due to tightness also weakly to
  the process $Y_{\cdot,s}^{\zeta,x}$.
  This completes the proof as the remaining cases are analogous.
\end{proof}
%
\begin{lemma}   \label{l:small_initial_mass_weak}
  Assume~\ref{a:A1} and~\ref{a:Hutzenthaler2009EJP}.
  Then we have that
  \begin{equation}   \label{eq:small_initial_mass_weak}
     \robb{\sum_{i=1}^{N}
      \dl_{Y_{t,s}^{0,\frac{x}{N}}(i)}}_{s\leq t\leq T}\wlra
      \robb{\int_0^x\int \dl_{\eta_{t-s}}   \Pi(dy,d\eta)}_{s\leq t\leq T}   \qqasN
  \end{equation}
  for all $x\geq0$ and all $s\leq T$
  where $\Pi$ is a Poisson point process on $[0,\infty)\times U$
  with intensity measure
  $\Leb\otimes Q$.
\end{lemma}
\begin{proof}
  The proof is analogous to the proof of
  Lemma~\ref{l:vanishing_immigration_weak_process}.
  For convergence of fi\-nite-di\-men\-sio\-nal distributions use the
  convergence~\eqref{eq:Q} instead of
  Lemma~\ref{l:vanishing_immigration_weak_process}.
  Tightness follows from an estimate as in~\eqref{eq:esti_tightiY}
  together with boundedness (see Lemma 9.9 in~\cite{Hutzenthaler2009EJP}) of second moments.
\end{proof}
Finally we prove the main result of this subsection.
The following lemma shows that the total mass increases if we let all
subfamilies evolve independently.
In the special case of $\mu$ and $\sigma^2$ being linear,
inequality~\eqref{eq:family_decomposition_ix} is actually an equality
according to the classical family decomposition of Feller's branching
diffusion with immigration.
%
\begin{lemma}   \label{l:family_decomposition}
  Assume~\ref{a:A1}.
  Let $\zeta\colon[0,\infty)\to[0,\infty)$
  be locally Lebesgue integrable
  and let $x\in I$.
  If the drift function $\mu$ is concave and the diffusion function
  $\sigma^2$ is superadditive, then
  \begin{equation}   \label{eq:family_decomposition_ix}
    \ruB{Y_{t,s}^{\zeta,x}}_{t\geq s}
    \leq_{\SetF_{+-}\ro{[0,\infty)}}
    \Big(\int_0^x \int\eta_{t-s}\Pi(dy,d\eta)
        +\int_s^\infty\int\eta_{t-u}\Pit(du,d\eta)\Big)_{t\geq s}
  \end{equation}
  for every $s\geq0$
  where $\Pi$ is a Poisson point process on $[0,\infty)\times U$
  with intensity measure $\Leb\otimes Q$ and where $\Pit$ is an independent
  Poisson point process on $[0,\infty)\times U$ with intensity measure
  $\zeta(s)\,ds\otimes Q$.
  If $\mu$ is concave and $\sigma^2$ is subadditive,
  then~\eqref{eq:family_decomposition_ix}
  holds with $\SetF_{+-}$ replaced by $\SetF_{++}$.
  If $\mu$ is subadditive and $\sigma^2$ is additive, then~\eqref{eq:family_decomposition_ix}
  holds with $\SetF_{+-}$ replaced by $\SetF_{+\pm}$.
\end{lemma}
\begin{proof}
  The idea is to split the initial mass and the immigrating mass
  into smaller and smaller pieces.
  Fix $s\geq0$.
  Let $\mu$ be concave and let $\sigma^2$ be superadditive.
  According to Lemma~\ref{l:semigroup_estimate}
  \begin{equation}  \label{eq:smaller_pieces}
    \ruB{Y_{t,s}^{\zeta,x}}_{t\geq s}
    \leq_{\SetF_{+-}\ro{[0,\infty)}}
    \biggl(\sum_{i=1}^N   Y_{t,s}^{0,\frac{x}{N}}(i)
        +\sum_{i=1}^N\Yt_{t,s}^{\frac{\zeta}{N},0}(i)\biggr)_{t\geq s}
  \end{equation}
  for every $N\in\N$
  where all processes are independent of each other.
  Letting $N\to\infty$ in~\eqref{eq:smaller_pieces},
  the right-hand side of~\eqref{eq:smaller_pieces}
  converges to the right-hand side of~\eqref{eq:family_decomposition_ix},
  see Lemma~\ref{l:vanishing_immigration_weak_process}
  and Lemma~\ref{l:small_initial_mass_weak}.
  The remaining cases are analogous.
\end{proof}
\subsection{The $(G,m,\mu,\sigma^2)$-process is dominated by the
loop-free $(G,m,\mu,\sigma^2)$-process}
\label{ssec:The process is dominated by the loop-free process}
\begin{lemma}  \label{l:Z_dominates_X}
  Assume~\ref{a:A1}.
  If $\mu$ is concave and $\sigma^2$ is superadditive, then
  \begin{equation}  \label{eq:semigroup_estimate_loop_free}
    \ruB{X_t}_{t\geq0}
    \leq_{\SetF_{+-}\ro{G,[0,\infty)}}
    \biggl(\sum_{k=0}^\infty Z_t^{(k)}\biggr)_{t\geq0}.
  \end{equation}
  If $\mu$ is concave and $\sigma^2$ is subadditive, then
  inequality~\eqref{eq:semigroup_estimate_loop_free} holds with
  $\SetF_{+-}$ replaced by $\SetF_{++}$.
  If $\mu$ is subadditive and $\sigma^2$ is additive, then
  inequality~\eqref{eq:semigroup_estimate_loop_free} holds with
  $\SetF_{+-}$ replaced by $\SetF_{+\pm}$.
\end{lemma}
\begin{proof}
  Assume that $\mu$ is subadditive and that
  $\sigma^2$ is superadditive.
  We follow the proof of Lemma~\ref{l:semigroup_estimate}
  and begin with a generator calculation similar to
  Lemma~\ref{l:generator_estimate}.
  Let $\Gen^X$ and $\Gen^Z$ denote the formal generators of
  $(G,m,\mu,\sigma^2)$-process and of the
  loop-free $(G,m,\mu,\sigma^2)$-process, respectively.
  Assume that $f\in\C^2\rob{[0,\infty)^G}\cap\MCF_{+-}^{(1)}\rob{[0,\infty)^G}$
  depends only on finitely many coordinates.
  Associated with this test function is
  \begin{equation}
    \ft\roB{(x_i^{(k)})_{i\in G,k\in\N_0}}:= f\robb{\sum_{k=0}^\infty x_\cdot^{(k)}}
  \end{equation}
  where $(x_i^{(k)})_{i\in G,k\in\N_0}\in I^{G\times\N_0}$.
  Note that $\ft\in\Dom(\Gen^Z)=\C^2\rob{I^{G\times\N_0}}$.
  The first partial derivatives
  $f_i:=\rob{\tfrac{\del}{\del x_i}}f\rob{\sum_{k=0}^\infty x_{\cdot}^{(k)}}$, $i\in G$,
  are non-negative and the second partial derivatives
  \begin{equation}
    f_{ii}:=\roB{\tfrac{\del^2}{\del x_i^2}}f\robb{\sum_{k=0}^\infty x_{\cdot}^{(k)}},
    \qquad i\in G,
  \end{equation}
  are nonpositive.
  Thus we see that
  \begin{equation}  \begin{split} \label{eq:generator_estimate_loop_free}
    \lefteqn{\ruB{\Gen^X f}\ruB{\sum_{k=0}^\infty x_\cdot^{(k)}}}\\
    &=\sum_{i\in G}f_i\eckbb
     {
      \sum_{j\in G}\ruB{\sum_{k=0}^\infty x_j^{(k)}m(j,i)-\sum_{k=0}^\infty x_i^{(k)}
      +\mu\ruB{\sum_{k=0}^\infty x_i^{(k)}} }
     }
     +\frac12\sum_{i\in G}f_{ii}\mal \sigma^2\ruB{\sum_{k=0}^\infty x_i^{(k)}}\\
    &\leq \sum_{k=0}^\infty\eckbb
     {
       \sum_{i\in G}f_i\ruB
       {
        \sum_{j\in G}x_j^{(k)}m(j,i)-x_i^{(k)}
        +\mu\rub{x_i^{(k)}} 
       }
       +\frac12\sum_{i\in G}f_{ii}\mal \sigma^2\rub{ x_i^{(k)} }
     }\\
    &= \sum_{k=0}^\infty\eckbb
     {
       \sum_{i\in G}f_i\ruB
       {
        \sum_{j\in G}x_j^{(k-1)}m(j,i)-x_i^{(k)}
        +\mu\rub{x_i^{(k)}} 
       }
       +\frac12\sum_{i\in G}f_{ii}\mal \sigma^2\rub{ x_i^{(k)} }
     }\\
    &= \Gen^Z \ft\roB{\rob{x_i^{(k)}}_{i\in G,k\in\N_0}}
  \end{split}     \end{equation}
  for every
  $(x_i^{(k)})_{i\in G,k\in\N_0}\in I^{G\times\N_0}$.
  Now we wish to apply the integration by parts formula~\eqref{eq:comp_int_by_parts}.
  In order to guarantee $x\mapsto \E^x f(X_s)\in\Dom(\Gen^X)$,
  approximate $f$ with smooth functions in $\MCF_{+-}^{(1)}$,
  approximate $\mu$ and $\sigma$ as
  in the proof of Lemma~\ref{l:semigroup_estimate}
  and approximate $G$ with finite sets.
  Moreover in order to exploit the generator
  inequality~\eqref{eq:generator_estimate_loop_free}
  in the integration by parts formula~\eqref{eq:comp_int_by_parts},
  we note that $x\mapsto \E^x f(X_s)\in\MCF_{+-}^{(1)}$
  (see Lemma 6.1 in~\cite{HutzenthalerWakolbinger2007}).
  Therefore the integration by parts formula~\eqref{eq:comp_int_by_parts}
  together with inequality~\eqref{eq:generator_estimate_loop_free} implies that
  \begin{equation}
    \E^{\sum_{k=0}^\infty x_\cdot^{(k)}}\ruB{f\rub{X_t}}
    \leq
    \E^{ x_\cdot^{(\cdot)}}\rubb{f\ruB{\sum_{k=0}^\infty Z_t^{(k)}}}
    \qquad\fa f\in\MCF_{+-}^{(1)}\rob{[0,\infty)^G}
  \end{equation}
  for all
  $(x_i^{(k)})_{i\in G,k\in\N_0}\in I^{G\times\N_0}$.
  In addition note that $(X_t)_{t\geq0}$ is stochastically non-decreasing
  in its initial configuration, see Lemma 3.3 in~\cite{HutzenthalerWakolbinger2007}.
  As stochastic monotonicity is the only input into the proof of
  Lemma~\ref{l:preservation},
  the assertion of Lemma~\ref{l:preservation} holds also for $(X_t)_{t\geq0}$.
  Thus, for all $f\in\MCF_{+-}^{(n+1)}\rob{[0,\infty)^G}$, we have that
  \begin{equation}
    \rob{I^G}^{n}\ni (x_1,\ldots,x_{n})\mapsto 
    \E^{x_n} f\rob{x_1,\ldots,x_n,X_t}
    \in\MCF_{+-}^{(n)}\rob{I^G}.
  \end{equation}
  Using this, the assertion follows as in
  the proof of Lemma~\ref{l:semigroup_estimate}
  by induction on the number of arguments of $F\in\SetF_{+-}\rob{[0,\infty)^G}$.
\end{proof}

\subsection{The  loop-free process is dominated by the virgin island model}
\label{sec:V_dominates_Z}
We show in Lemma~\ref{l:V_dominates_Z} below that the total mass of the
loop-free process is dominated by the total mass of the virgin island model.
In the proof of this lemma, we use that the Poisson point processes appearing
in the definition of the virgin island model preserve convexity in a suitable way.
This is subject of the following lemma.
\begin{lemma}  \label{l:preservation_PPP}
  For every vector $\uline{z}=(z_1,\ldots,z_m)\in[0,\infty)^m$,
  $m\in\N$, let
  $\Pi^{(\uline{z})}$ be a Poisson point process on a Polish space $S$
  with intensity measure $\sum_{i=1}^{m}z_i\mu_i$
  where
  $\mu_1,\ldots,\mu_m$
  are fixed measures on $S$.
  If $f\in\MCF_{+-}^{(n+1)}\rob{[0,\infty)}$, $n\in\N$, then
  the function
  \begin{equation}
    \ft\colon[0,\infty)^{n+m}\to\R,
    \;(\uline{x},\uline{z})\mapsto \E f\rob{\uline{x},\scal{g}{\Pi^{(\uline{z})}}}
  \end{equation}
  is an element of $\MCF_{+-}^{(n+m)}\rob{[0,\infty)}$ for every
  measurable test function $g\colon S\to[0,\infty)$
  satisfying $\scal{g}{\mu_i}=\int g\,d\mu_i<\infty$ for $i=1,\ldots,m$.
  Analogous
  results hold if $\MCF_{+-}$ is replaced by $\MCF_{++}$
  and $\MCF_{+\pm}$, respectively.
\end{lemma}
\begin{proof}
  The function $\ft$ is non-decreasing
  in the first $n$
  variables and $(i,j)$-concave
  for all $1\leq i,j\leq n$.
  Furthermore $\ft$ is non-decreasing in the last $m$ variables
  as $\Pi^{(\uline{z})}$ is stochastically non-decreasing in $\uline{z}$.
  Fix $1\leq i,j\leq m$ and $h_1,h_2\geq0$.
  Let the Poisson point processes 
  $\Pi^{(\uline{z})}$, $\Pi^{(h_1 e_i)}$ and $\Pi^{(h_2 e_j)}$
  be independent.
  Fix $\uline{z}\in[0,\infty)^m$, $\uline{x}\in[0,\infty)^n$
  and a measurable test function $g\colon S\to[0,\infty)$.
  Note that
  $\Pi^{(\uline{z})}+\Pi^{(h_1 e_i)}+\Pi^{(h_2 e_j)}$
  has the same distribution as $\Pi^{(\uline{z}+h_1 e_i+h_2 e_j)}$.
  Therefore, using $(n+1,n+1)$-concavity of $f$ in the point
  $(\uline{x},\scal{g}{\Pi^{(\uline{z})}})$, we obtain that
  \begin{equation}  \begin{split}
    \lefteqn{\E f\roB{\uline{x},\scalb{g}{\Pi^{(\uline{z}+h_1 e_i+h_2 e_j)}}}
        +\E f\roB{\uline{x},\scalb{g}{\Pi^{(\uline{z})}}
             }
            }\\
    &=\E f\roB{\uline{x},\scalb{g}{\Pi^{(\uline{z})}}
                 +\scalb{g}{\Pi^{(h_1 e_i)}}
                 +\scalb{g}{\Pi^{(h_2 e_j)}}
             }
        +\E f\roB{\uline{x},\scalb{g}{\Pi^{(\uline{z})}}
             }
             \\
    &\leq\E f\roB{\uline{x},\scalb{g}{\Pi^{(\uline{z})}}
                 +\scalb{g}{\Pi^{(h_1 e_i)}}
             }
        +\E f\roB{\uline{x},\scalb{g}{\Pi^{(\uline{z})}}
                 +\scalb{g}{\Pi^{(h_2 e_j)}}
             } \\
    &= \ft(\uline{x},\uline{z}+h_1 e_i)
      +\ft(\uline{x},\uline{z}+h_2 e_j).
  \end{split}     \end{equation}
  For the last step, note that
  $\Pi^{(\uline{z})}+\Pi^{(h_1 e_i)}$
  has the same distribution as $\Pi^{(\uline{z}+h_1 e_i)}$
  and that
  $\Pi^{(\uline{z})}+\Pi^{(h_2 e_j)}$
  has the same distribution as $\Pi^{(\uline{z}+h_2 e_j)}$.
  This proves $(n+i,n+j)$-concavity of $\ft$ for all
  $1\leq i,j\leq m$.
  Similar arguments prove $(i,n+j)$-concavity for $1\leq i\leq n$
  and $1\leq j\leq m$.
\end{proof}

Recall the total mass process $(V_t^{(n)})_{t\geq0}$ of
the $n$-th generation of the virgin island model
from~\eqref{eq:def:Vn} for every $n\in\N_0$.
Define $|Z_t^{(k)}|:=\sum_{i\in G}Z_t^{(k)}$ for $t\in[0,\infty)$ and $k\in\N_0$.
\begin{lemma}   \label{l:V_dominates_Z}
  Assume~\ref{a:A1} and~\ref{a:migration}.
  If $\mu$ is concave and $\sigma^2$ is superadditive, then
  \begin{equation}  \label{eq:V_dominates_Z}
    \ruB{\rub{\absb{Z_t^{(k)}}}_{k=0,\ldots,k_0}}_{t\geq0}
    \leq_{\SetF_{+-}\rob{\{0,\ldots,k_0\},[0,\infty)}}
    \ruB{\rub{V_t^{(k)}}_{k=0,\ldots,k_0}}_{t\geq0}
  \end{equation}
  for every $k_0\in\N_0$.
  If $\mu$ is concave and $\sigma^2$ is subadditive, then
  inequality~\eqref{eq:V_dominates_Z} holds with
  $\SetF_{+-}$ replaced by $\SetF_{++}$.
  If $\mu$ is subadditive and $\sigma^2$ is additive, then
  inequality~\eqref{eq:V_dominates_Z} holds with
  $\SetF_{+-}$ replaced by $\SetF_{+\pm}$.
\end{lemma}
\begin{proof}
  We  prove~\eqref{eq:V_dominates_Z} by induction on $k_0\in\N_0$.
  The base case $k_0=0$ follows from $\abs{Z_\cdot^{(0)}}\eqd V_{\cdot}^{(0)}$.
  We apply Lemma~\ref{l:family_decomposition}
  for the induction step $k_0\to k_0+1$.
  Fix
  \begin{equation}
    F_{n+1}(\chi)=f_{n+1}(\chi_{s_1}^{(k_1)},\ldots,\chi_{s_{n+1}}^{(k_{n+1})})
    \in\SetF_{+-}\rob{\{0,\ldots,k_0+1\},[0,\infty)}
  \end{equation}
  where $k_1,\ldots,k_{n+1}\in\{0,\ldots,k_0+1\}$
  and $0\leq s_1\leq\cdots\leq s_{n+1}$.
  Let $\Pi^{i,\zeta}$, $i\in G,\zeta\in\C\rob{[0,\infty),[0,\infty)}$,
  be independent Poisson point processes on $[0,\infty)\times U$, independent
  of $\{Z^{(0)},\ldots,Z^{(k_0)}\}$ and with
  $\Pi^{i,\zeta}\eqd \Pi^{\zeta}$ where $\Pi^\zeta$ has intensity measure
  \begin{equation}
    \E \Pi^\zeta(du,d\eta)=\zeta(u)du\otimes Q(d\eta).
  \end{equation}
  Note that conditioned on $\zeta_i^{(k_0)}(\cdot):=\sum_{j\in G}Z_\cdot^{(k_0)}(j)m(j,i)$,
  the law of $Z_{\cdot}^{(k_0+1)}(i)$ is equal to the law of
  $Y_{\cdot,0}^{\zeta_i^{(k_0)},0}(i)$ (defined in~\eqref{eq:Y_immi})
  where $Y_{\cdot,0}^{\zeta_i^{(k_0)},0}(i)$, $i\in G$, are independent
  of each other.
  Thus conditioning on $\{Z^{(0)},\ldots,Z^{(k_0)}\}$ and applying
  Lemma~\ref{l:family_decomposition} with $x=0$
  and $s=0$, we obtain that
  \begin{equation}  \begin{split}   \label{eq:estimate_Z_V_1}
    \lefteqn{
      \E \eckbb{F_{n+1}\roB{\rob{\absb{Z^{(k)}}}_{k=0,\ldots,k_0+1}}}
    }\\
    &= \E\curlbb{\E \eckbb{F_{n+1}\roB{\rob{\absb{Z^{(k)}}}_{k=0,\ldots,k_0},
        \sum_{i\in G} Y_{\cdot,0}^{\zeta_i^{(k_0)},0}(i)}
           \Bigm|\rob{Z^{(k)}}_{k=0\ldots k_0}
                }         }\\
    &\leq \E \eckbb{ F_{n+1}\roB{\rob{\absb{Z^{(k)}}}_{k=0,\ldots,k_0},
        \sum_{i\in G} \int_0^\infty \int\eta_{\cdot-u}\Pi^{i,\zeta_i^{(k_0)}}(du,d\eta)
                   }       }\\
    &=    \E \eckbb{ F_{n+1}\roB{\rob{\absb{Z^{(k)}}}_{k=0,\ldots,k_0},
         \int_0^\infty \int\eta_{\cdot-u}\Pi^{\sum_{i\in G} \zeta_i^{(k_0)}}(du,d\eta)
                   }       }\\
    &\leq    \E \eckbb{ F_{n+1}\roB{\rob{\absb{Z^{(k)}}}_{k=0,\ldots,k_0},
         \int_0^\infty \int\eta_{\cdot-u}\Pi^{\abs{Z^{(k_0)}}}(du,d\eta)
                   }       }.
  \end{split}     \end{equation}
  The last inequality follows from $\sum_{i\in G}\zeta_i^{(k_0)}(\cdot)
  =\sum_{j\in G}Z_\cdot^{(k_0)}\sum_{i\in G}m(j,i)\leq\abs{Z_\cdot^{(k_0)}}$,
  where we used the inequality $\sum_{i\in G}m(j,i)\leq 1$ from Assumption~\ref{a:migration}, and
  where we used that $F_{n+1}$ is non-decreasing.
  Next we would like to apply the induction hypothesis.
  However the right-hand side of~\eqref{eq:estimate_Z_V_1} depends
  on $(\abs{Z_t^{(k_0)}})_{t\geq0}$ through a continuum of time
  points and not only through finitely many
  time points.
  To remedy this, we approximate the Poisson point process on the
  right-hand side of~\eqref{eq:estimate_Z_V_1} by approximating
  $(\abs{Z_t^{(k_0)}})_{t\geq0}$ with simple functions.
  For each $\N\ni m\geq\max(s_n,2n)$,
  choose a discretization $\{t_0,\ldots,t_{m^2}\}$
  of $[0,m]$ of maximal width $\tfrac{2}{m}$ such that
  $t_0=0$, $t_{m^2}=m$ and
  $\{s_1,\ldots,s_n\}\subset\{t_0,\ldots,t_{m^2}\}$.
  Define $l_i:=k_i$ if $t_i\in\{s_1,\ldots,s_n\}$
  and $l_i:=0$ otherwise.
  For a path $(\chi_t)_{t\geq0}\in\C\rob{[0,\infty),[0,\infty)}$,
  define
  \begin{equation}
    \rob{D_m\chi}(t):=\sum_{i=1}^{m^2}\chi_{t_{i-1}}\1_{[t_{i-1},t_i)}(t)
    \quad\text{for }t\geq0.
  \end{equation}
  Note that $(D_m\chi)(t)\to\chi_t$ for every $t\geq0$ as $m\to\infty$.
  Thus the intensity measure $\E\Pi^{D_m\chi^{(k_0)}}(du,d\eta)$
  converges to $\E\Pi^{\chi^{(k_0)}}(du,d\eta)$ as $m\to\infty$.
  This convergence of the intensity measures
  implies weak convergence of the Poisson point process $\Pi^{D_m\chi^{(k_0)}}$
  to the Poisson point process $\Pi^{\chi^{(k_0)}}$.
  Due to Lemma~\ref{l:preservation_PPP},
  the function $\fb\colon [0,\infty)^{m^2+1}\to\R$ defined through
  \begin{equation}
    \fb\roB{\chi_{t_0}^{(l_0)},\ldots,\chi_{t_{m^2}}^{(l_{m^2})}}
    =\E f_{n+1}\roB{\chi_{s_1}^{(k_1)},\ldots,\chi_{s_n}^{(k_n)},
       \int_0^\infty \int\eta_{\cdot-u}\Pi^{D_m\chi^{(k_0)}}(du,d\eta)
              }
  \end{equation}
  is an element of $\MCF_{+-}^{m^2+1}\rob{[0,\infty)}$.
  Now we apply the induction hypothesis and obtain that
  \begin{equation}  \begin{split}  \label{eq:last_eqn}
    &\E \eckbb{ F\roB{\rob{\absb{Z^{(k)}}}_{k=0,\ldots,k_0},
            \int_0^\infty \int\eta_{\cdot-u}\Pi^{\abs{Z^{(k_0)}}}(du,d\eta)
                   }       }\\
    &=\limm\E \eckbb{ F\roB{\rob{\absb{Z^{(k)}}}_{k=0,\ldots,k_0},
            \int_0^\infty \int\eta_{\cdot-u}\Pi^{D_m\abs{Z^{(k_0)}}}(du,d\eta)
                   }       }\\
    &\leq\limm\E \eckbb{ F\roB{\rob{{V^{(k)}}}_{k=0,\ldots,k_0},
            \int_0^\infty \int\eta_{\cdot-u}\Pi^{D_m{V^{(k_0)}}}(du,d\eta)
                   }       }\\
    &=\E \eckbb{ F\roB{\rob{{V^{(k)}}}_{k=0,\ldots,k_0},
            \int_0^\infty \int\eta_{\cdot-u}\Pi^{{V^{(k_0)}}}(du,d\eta)
                   }       }\\
    &=\E \eckbb{ F\roB{\rob{{V^{(k)}}}_{k=0,\ldots,k_0},
            V^{(k_0+1)}
                   }       }.
  \end{split}     \end{equation}
  Putting~\eqref{eq:estimate_Z_V_1} and~\eqref{eq:last_eqn} together
  completes the induction step.
\end{proof}
\subsection{Proof of the comparison result of Theorem~\ref{thm:comparison}}
\label{ssec:Proof of the comparison result of Theorem 2}

\begin{proof}[Proof of Theorem~\ref{thm:comparison}]
  We prove the case of $\mu$ being concave and $\sigma^2$ being
  superadditive. The remaining two cases are analogous.
  According to Lemma~\ref{l:Z_dominates_X} we have that
  \begin{equation}
    \biggl(\sum_{i\in\Lambda}X_t(i)\biggr)_{t\geq0}
    \leq_{\SetF_{+-}\ro{[0,\infty)}}
    \biggl(\sum_{i\in\Lambda}\sum_{k=0}^\infty Z_t^{(k)}(i)\biggr)_{t\geq0}
  \end{equation}
  for every finite subset $\Lambda\subseteq G$. Letting $\Lambda\nearrow G$,
  we see that the total mass of the
  $(G,m,\mu,\sigma^2)$-process is dominated by the total mass of the loop-free
  $(G,m,\mu,\sigma^2)$-process.
  Now we get from Lemma~\ref{l:V_dominates_Z} that
  \begin{equation}
    \biggl(\sum_{k=0}^{k_0}\absb{Z_t^{(k)}}\biggr)_{t\geq0}
    \leq_{\SetF_{+-}\ro{[0,\infty)}}
    \biggl(\sum_{k=0}^{k_0}V_t^{(k)}\biggr)_{t\geq0}
  \end{equation}
  for every $k_0\in\N$. Letting $k_0\to\infty$, we obtain that
  the total mass of the loop-free
  $(G,m,\mu,\sigma^2)$-process is
  dominated by the total mass of the virgin island model.
  Therefore, the total mass of the
  $(G,m,\mu,\sigma^2)$-process is dominated by the total mass of the
  virgin island model.
\end{proof}




\def\cprime{$'$}
  \hyphenation{Sprin-ger}\def\cftil$1{\ifmmode\setbox7\hbox{$\accent"5E$1$}\el%
se \setbox7\hbox{\accent"5E$1}\penalty 10000\relax\fi\raise 1\ht7
  \hbox{\lower1.15ex\hbox to 1\wd7{\hss\accent"7E\hss}}\penalty 10000
  \hskip-1\wd7\penalty 10000\box7} \def\cprime{$'$}


\ACKNO{
I thank Anton Wakolbinger and Don Dawson for inspiring discussions and valuable remarks.
Moreover I am indebted to two anonymous referees for very helpful comments
and suggestions which led to a considerable improvement of the presentation
and of the proofs.
}


\end{document}